\documentclass[10pt,centertags]{amsart}
\usepackage{amsmath,amssymb}
\usepackage{hyperref}
\newtheorem*{main*}{Main Theorem}

\newtheorem{theorem}{Theorem}[section]
\newtheorem*{theorem*}{Theorem}
\newtheorem{proposition}[theorem]{Proposition}
\newtheorem{lemma}[theorem]{Lemma}
\newtheorem{corollary}[theorem]{Corollary}

\newtheorem*{question*}{Question}

\newtheorem*{conjecture*}{Conjecture}

\theoremstyle{definition}

\newtheorem{definition}[theorem]{Definition}
\newtheorem{definitions}[theorem]{Definitions}
\newtheorem*{definition*}{Definition}
\newtheorem{example}[theorem]{Example}

\theoremstyle{remark}
\newtheorem{remark}[theorem]{Remark}
\newtheorem{remarks}[theorem]{Remarks}

\numberwithin{equation}{section}

\def\ds{\displaystyle}

\newcommand{\C}{\mathbb{C}}
\newcommand{\R}{\mathbb{R}}
\renewcommand{\H}{\mathbb{H}}
\newcommand{\HH}{\mathbb{H}}
\newcommand{\Z}{\mathbb{Z}}
\newcommand{\PP}{\mathbb{P}}

\newcommand{\eps}{\epsilon }

\newcommand{\set}[1]{\left\{ #1 \right\}}
\newcommand{\op}[1]{\operatorname{#1}}

\newcommand{\id}{{\mathbf e}}

\newcommand{\Isom}{\operatorname{Isom}}

\newcommand{\D}{\partial }

\newcommand{\til}{\widetilde}

\newcommand{\of}{\circ }
\providecommand{\to}{\rightarrow }

\renewcommand{\bar}{\overline}

\newcommand{\abs}[1]{\left\lvert #1 \right\rvert }
\newcommand{\norm}[1]{\left\lVert #1 \right\rVert }

\newcommand{\conv}{\star}
\newcommand{\group}{\Gamma}

\renewcommand{\(}{\left(}
\renewcommand{\)}{\right)}
\newcommand{\numeq}[1]{\begin{align}\begin{split} #1
\end{split}\end{align}}

\DeclareMathOperator*{\ess-sup}{ess\ sup}
\def\[#1\]{\begin{align*}\begin{split} #1 \end{split}\end{align*} }

\def\ssu{\subset}

\def\pr{\prime}

\def\nn{ \mathbb N}
\def\rr{\mathbb R}

\def\Ga{ \Gamma}
\def\ga{\gamma}

\def\la{ \lambda}
\def\La{ \Lambda}
\def\til{ \tilde}
\def\<{\langle}
\def\>{\rangle}
\def\pr{\prime}
\def\pa{\partial}
\def\ep{\epsilon}

\begin{document}
\address{Indiana University} %
\email{cconnell@math.indiana.edu} %

\address{University of Chicago} %
\email{roma@math.uchicago.edu} %

\author[C. Connell]{Chris Connell$^\dagger$}
\thanks{$\dagger$ The author was
  supported in part by an NSF postdoctoral fellowship and DMS-0306594.}
\author[R. Muchnik]{Roman Muchnik$^\ddagger$}
\thanks{$\ddagger$ The author was
  supported in part by an NSF postdoctoral fellowship.}

\title[Stationary Measures]{Harmonicity of quasiconformal measures and
Poisson boundaries of hyperbolic spaces}

\begin{abstract}
  We consider a group $\group$ of isometries acting on a (not necessarily
  geodesic) $\delta$-hyperbolic space and possessing a radial limit set of
  full measure within its limit set. For any continuous
  $\alpha$-quasiconformal measure $\nu$ supported on the limit set, we
  produce a stationary measure $\mu$ on $\group$. Moreover the limit set
  together with $\nu$ forms a $\mu$-boundary and $\nu$ is harmonic with
  respect to the random walk induced by $\mu$. In the case when $X$ is a
  CAT$(-\kappa)$ space and $\group$ acts cocompactly, for instance, we
  show that $\mu$ has finite first moment. This implies that $(\pa X,
  \nu)$ is the unique Poisson boundary for $\mu$. In the course of the
  proofs, we establish sufficient conditions for a set of continuous
  functions to form a positive basis, either in the $L^1$ or $L^\infty$
  norm, for the space of uniformly positive lower-semicontinuous functions
  on a metric measure space.
\end{abstract}
\maketitle

\thispagestyle{empty} %

On the hyperbolic plane $\HH^2$ we can represent any bounded harmonic
function $h$
by the formula
$$h(x)=h_f(x)=\int_{\partial \HH^2}f d\nu_x$$
for some $f\in
L^\infty(\partial \HH^2)$
where $\partial \HH^2$ is the circle at infinity and $\nu_x$ are the
harmonic measures. Representing $\HH^2$ by the unit disk in $\C$, the
harmonic measure corresponding to the origin, $\nu_0$, is just the
unit Lebesgue measure on $S^1$. The others are given by
$$\frac{d\nu_x}{d\nu_0}(z)=\frac{1-|x|^2}{|z-x|^2}.$$
The measures $\nu_x$ also arise from $\nu_0$ by image measures under the
transitive isometry group: $\nu_x=g_*\nu_0$ for any
$g\in\Isom(\HH^2)$ such that $g(0)=x$.

Since $\nu_x$ tends to the dirac measure at $z\in \D \HH^2$ as
$x\to z$, we obtain
$$\norm{h_f-h_g}_{L^\infty(\HH^2)}=\ess-sup_{x\in X} \abs{\int_{\D \HH^2}
  f(z)-g(z)\ d\nu_x(z)}=\norm{f-g}_{L^\infty(\D \HH^2)}.$$
In other words, we have a maximum principle so that the map
$f\mapsto h_f$ determines an isometry, with respect to the
$L^\infty$ norms, between the Banach spaces $H^\infty(\HH^2)$ of
all bounded harmonic functions and $L^\infty(\partial
\HH^2,\nu_0)$. Another consequence of the maximum principle is
that harmonic functions satisfy the averaging condition. If we
identify $\HH^2$ with any transitive Lie subgroup of
$\Isom(\HH^2)$, then at each point $x\in\HH^2$ any harmonic
function $h_f$ satisfies
$$h_f(x)=\int_{S} h_f(x y) d\mu(y),$$
where $\mu$ is the uniform measure on the unit distance circle
$S\subset \HH^2$ around the identity element with respect to the
hyperbolic metric.

We can generalize these concepts greatly to any measure space $X$
with a Markov operator $P$. We say that a function $h$ is
$P$-harmonic if $P h=h$. A space $B$ together with a family of
mutually absolutely continuous measures $\set{\nu_x}_{x\in X}$ is a
{\em Poisson Boundary} if the map $f\mapsto h_f$ given by the {\em
Poisson formula}
\begin{gather}\label{eq:Poisson_form}
   h_f(x):=\int_B f(y)d\nu_x(y)
\end{gather}
is an isometry between $L^\infty(B,\set{\nu_x})$ and bounded
$P$-harmonic functions on $X$. In this event, $\nu_x$ is called the
{\em harmonic measure at $x$}. The Poisson boundary is a purely
measurable object and is unique up to measurable isomorphism. As
such, even when $X$ is a manifold, the underlying space $B$ could be
quite different from that of the topological Martin boundary, which
gives a representation space for $\mu$-harmonic functions. The
Poisson boundary always exists and has many equivalent descriptions.
For instance, it can be identified with the space of ergodic
components of the shift map $T$ acting on the space of sample paths
of the Markov chain on $X$ associated with the operator $P$. The
harmonic measures $\nu_x$ are the images under the quotient map of
the measures $\PP_x$ in the path space corresponding to starting the
Markov process from state $x\in X.$ For other characterizations of
the Poisson boundary see \cite{Kaimanovich03}.

Before moving to random walks on groups, we recall the notion of convolution
measures. Suppose $(X,\nu)$ is a measure space and $\Gamma$ is any
set of $\nu$-measurable transformations of $X$. For any measure
$\mu$ on $\Gamma$ we define the convolution of the two measures
$\mu\conv \nu$ to be the measure on $X$ given by $$\mu\conv
\nu :=\int_{\Gamma} \gamma_*\nu\ \  d\mu(\gamma).$$

As a special case of the above constructions, we can define the (right
sided) random walk determined by a finite measure $\mu$ on a group $G$
as follows. Let $G^{\Z_+}=\prod_{i=1}^\infty G$ and denote by $\PP$
the measure obtained as the image of $\mu\times\mu\times\cdots$ under
the map $(x_1,x_2,x_3,\dots)\mapsto (e,x_1,x_1x_2,x_1x_2x_3,\dots)$.
The conditional measure for $\PP$ on the $n-th$ coordinate of
$G^{\Z_+}$ is the $n$-fold convolution measure $\mu_n=\mu^{\conv
  n}$ on $G$. The measure of all paths starting from an initial
distribution $\theta$ is $\theta\conv\PP(G^{\Z_+})$. In this context,
the natural Markov operator associated to the random walk is $P_\mu:
L^\infty(G,\mu)\to L^\infty(G,\mu)$ defined by
$$P_\mu(f)(g)=\int_G f(gh)d\mu(h).$$
In this setting we call $P_\mu$-harmonic functions simply
$\mu$-harmonic. Moreover, for a Poisson boundary
$(B,\set{\nu_g}_{g\in G})$ it follows that $\nu_g=g_*\nu_e$ so we
may simply write the boundary as $(B,\nu)$ where $\nu=\nu_e$. It
follows from the definition of $\mu$-harmonicity and the Poisson
formula \eqref{eq:Poisson_form} that for all $f\in L^\infty(B,\nu)$
at $e$ we have,
$$\nu(f)=h_f(e)=\int_G h_f(g)d\mu(g)=\int_G g_*\nu(f) d\mu(g).$$
In short, $\mu\conv \nu=\nu$, in which case we also say $\nu$ is
{\em $\mu$-stationary}. The importance of $\mu$-stationary measures
is that the Poisson formula \eqref{eq:Poisson_form} yields
$\mu$-harmonic functions: for any $x\in G$,

\[
\int_G h_f(x g)\mu(g)&=\int_G\int_B f(z)
d(xg)_*\nu(z)d\mu(g)\\&=\int_B f(z)dx_*\left(\int_G g_*\nu
d\mu(g)\right)(z)=\int_B f(z)dx_*\nu(z)=h_f(x).
\]

A measured $G$-space $(B^\pr,\nu^\pr)$ is called a {\em
  $\mu$-boundary of $G$} if the corresponding Poisson formula, $f\mapsto
\(g\mapsto h_f(g):=g_*\nu^\pr(g)\)$, defines an isometric
embedding from $L^\infty(B^\pr,\nu^\pr)$ into the space of bounded
$\mu$-harmonic functions $H^\infty(G,\mu)$.  Any $\mu$-boundary
arises as a $G$-equivariant measurable quotient $\pi:(B,\nu)\to
(B^\pr,\nu^\pr)$ of the Poisson boundary $(B,\nu)$ since the
induced lift map $\pi_*:L^\infty(B^\pr,\nu^\pr)\to
L^\infty(B,\nu)$ is an isometric embedding. In particular, the
Poisson boundary can be characterized as the maximal
$\mu$-boundary (see the unpublished survey \cite{Kaimanovich03}).
Since the Poisson formula still holds, we also have
$\mu\conv\nu^\pr=\nu^\pr$ for any $\mu$-boundary
$(B^\pr,\nu^\pr)$. On the other hand, it does not follow that any
$G$-space with a $\mu$-stationary measure is isomorphic to a
$\mu$-boundary. However, if $\nu^\pr$ is $\mu$-stationary and all
of the Dirac measures $\set{\delta_z}_{z\in B^\pr}$ occur in the
weak-$*$ closure of $\set{g\nu^\pr}_{g\in \group}$, then as in the
case of $\HH^2$ above, we have
$\norm{h_f-h_g}_{L^\infty(G,\mu)}=\norm{f-g}_{L^\infty(B^\pr,\nu^\pr)}$
which implies that $(B^\pr,\nu^\pr)$ is a $\mu$-boundary. This
last condition will hold in the case that the action of $G$ on
$(B^\pr,\nu^\pr)$ is {\em $\mu$-proximal} in the sense of
Furstenberg (see \cite{Furstenberg73}).  This construction can be
useful for identifying $\mu$-boundaries as subsets of
geometrically defined boundaries for $G$.

For example, consider the case $G=\Isom(\HH^2)=SL(2,\R)$ with a
maximal compact subgroup $K=\op{Stab}(o)\equiv SO(2)$ for a fixed
basepoint $o\in \HH^2$. Let $m_K$ denote the (bi-invariant) Haar
measure on $K$ which we will think of as a measure on $G$ supported
on $K$, and choose an element $g\in G$ such that $d(o,go)=1$. If
$\mu_0=m_K\conv g_*m_K$, then $\mu_0$-harmonic functions on $G$ are
right $K$ invariant and their quotients in $G/K=\H^2$ are ordinairy
harmonic functions. Moreover every $\mu_0$ harmonic function is the
lift of one on $\H^2$. Hence $(\partial\H^2,\nu_0)$ is the Poisson
boundary of $(G,\mu_0)$.

In fact this same correspondence was established by Furstenberg in
\cite{Furstenberg63} for any symmetric space $G/K$ where $G$ is
semisimple Lie group of noncompact type and $K$ is a maximal compact
subgroup. If we again take the bi-$K$ invariant measure
$\mu_0=m_K\conv g_*m_K$, then the Poisson boundary of $(G,\mu_0)$ is
the Furstenberg boundary $G/P$ together with the unique
$K$-invariant measure $\nu_0$. Later, Furstenberg extended this in
\cite{Furstenberg67} and \cite{Furstenberg71} to show that for any
lattice $\Gamma < G$, one can build a measure $\mu$ on $\Gamma$ for
which $(G/P,\nu_0)$ is the Poisson boundary. Geometric intuition
makes it tempting to believe that passing from a Lie Group to a
lattice, at least a uniform lattice, should be a simple operation
when it comes to garnering asymptotic information of any kind.
However, the measure $\mu$ constructed by Furstenberg is quite
different from the measure $\mu_0$ on $G$. For instance, there are
no $K$-invariant measures on the lattice and the measure $\mu$ need
not be compactly supported. At least in the rank one case, we will
show that there is an infinite dimensional space of measures $\mu$
on $\Gamma$ for which $(G/P,\nu_0)$ is a Poisson boundary.

Furstenberg proved this result for $\Gamma$ in two steps. First he
constructed a $\mu$ for which $\nu_0$ was $\mu$-stationary. By the
argument mentioned above he concluded that this was a
$\mu$-boundary. Next he showed that any $\mu$-boundary which has
finite first moment is the Poisson boundary. Here finite first
moment means,
$$\sum_{\gamma\in\Gamma} \mu(\gamma)d(e,\gamma)<\infty.$$
He used this geometric characterization of the Poisson boundary to
distinguish envelopes of discrete groups.  Namely, a discrete
group $G$ cannot be a lattice in $SL(n,\R)$ for two different
values of $n$.

Kaimanovich and Vershik in \cite{KaimanovichVershik83} gave sufficient and
necessary conditions for the Poisson boundary of a random walk on a
locally compact group to be trivial. In \cite{Kaimanovich00}, Kaimanovich
generalized this to a criterion to decide when a certain geometric
boundary for a group together with a family of exit measures could be a
Poisson boundary for a given random walk on the group. For $(B,\nu)$ to be
a Poisson boundary for $(G,\mu)$ they showed that in addition to
stationarity, $\mu\conv
\nu=\nu$, it is sufficient for $\mu$ to have both finite first
log-moment and finite entropy:
$$\sum_{\gamma\in\Gamma} \mu(\gamma)\log d(e,\gamma)<\infty\quad\text{and}\quad
\op{h}(\mu):=-\sum_{\gamma\in\group} \mu(\gamma)\log
\mu(\gamma)<\infty.$$

The goal of the present paper is to generalize Furstenberg's results
to Gromov hyperbolic groups and to a general class of boundary
measures. By so doing, we can partially answer a converse question
to that answered by Kaimanovich and Vershik's results stated above.
Namely, starting with a measure $\nu$ on $\partial G$ can we find a
measure $\mu$ on $G$ such that $(\D G,\nu)$ is the Poisson boundary
of $G$? In fact, not every measure $\nu$ can arise as a Poisson
boundary, and there are examples of measures $\nu$ which do arise
but $f\nu$ does not for certain positive measurable functions $f$
(see Remark \ref{rem:after_thm}). Nevertheless, we shall give an
affirmative answer for any measure Lipschitz equivalent to a
Patterson-Sullivan measure on a CAT$(-1)$ group. If one asks the
same question for $\mu$-boundaries instead, then we show existence
for continuous measures in this class on a large family of groups
which includes the Gromov hyperbolic groups.

In a second paper, we will broaden some of these results to multiple
measure classes within the family of Gibbs streams.

To achieve the stated goal we will extend the approach of
Furstenberg's original work (\cite{Furstenberg63}) to our wider
context. We restrict our attention here to spaces which are negatively
curved in a general sense. We hold out the hope that in the future
some of these techniques should also be able to address similar
problems for certain nonamenable nonpositively curved groups.

Consider a locally compact $\delta$-hyperbolic metric space $(X,d)$,
and let $\group$ be an arbitrary group of isometries of $X$. We will
assume $(X,d,\group)$ satisfies two mild conditions which we call
``Gromov product bounded'' and ``uniformly quasigeodesic'' which
hold whenever $X$ is a complete geodesic space. The first,
guarantees that the metric is well behaved near infinity and the
second guarantees that points are sufficiently well distributed in
$X$. Since $\group$ acts by isometries on $\partial X$, we can
consider its ideal limit set $\Lambda\subset\D X$ and radial (or
conical) limit set $\Lambda_r\subset \Lambda$ (see Section
\ref{sec:delta_hyp}). To avoid dealing with uninteresting cases, we
assume $\group$ is nonelementary. This means that $\Lambda$ has at
least $3$ points in it. Note that $\group$ itself need not be a
hyperbolic group. In fact, many relatively hyperbolic convergence
groups (in the sense of Bowditch \cite{Bowditch99}) will satisfy the
hypotheses of our first theorem below. For instance, if $X=\HH^3$
and $\group$ is a geometrically finite discrete group of isometries
then the complement of $\Lambda_r$ in $\Lambda$ is a countable dense
set of points, and $\Lambda$ may have topological codimension $1$ or
$2$ in $\partial X$.

The metric structure $(X,d)$ induces a natural 1-parameter family of
metrics  on $\partial X$ parameterized by $\eps$. If $\Gamma$ is
discrete and acts properly discontinuously, then the metric
corresponding to $\eps$ gives $\Lambda$ Hausdorff dimension
$\delta(\group)/\eps$, where $\delta(\group)$ is the critical
exponent for the action of $\group$. The corresponding Hausdorff
measure is in the same measure class as the Patterson-Sullivan
measures. If the Gromov product extends continuously to $\partial
X$, then these measures are examples of continuous
$\alpha$-quasiconformal densities for $\alpha=\delta(\group)$ (see
Section \ref{sec:delta_hyp} for all definitions and details).

\begin{theorem}[Stationarity]\label{thm:main1}
  Let $X$ be a Gromov product bounded, uniformly quasigeodesic,
  $\delta$-hyperbolic space and choose a nonelementary $\group<\op{Isom}(X)$ with
  limit set $\Lambda\subset \pa X$. Let $\nu$ be a continuous
  $\alpha$-quasiconformal measure on $\Lambda\subset\partial X$ for
  any $\alpha>0$. Suppose $\nu^{\pr}$ is an equivalent measure whose
  Radon-Nikodym derivative $\frac{d\nu^{\pr}}{d\nu}$ is a uniformly
  positive lower semicontinuous function. Assume $\group$ posesses a
  radial limit set $\Lambda_r$ of full $\nu$-measure in $\Lambda$.
  If $(X,\group)$ is quasiconvex cobounded or else $\nu$ belongs to
  a bounded quasiconformal density, then there exists a (nontrivial)
  measure $\mu$ on $\group$ such that $\mu\conv\nu=\nu^{\pr}$.
\end{theorem}

\begin{remark}
  The assumption that $X$ is uniformly quasiconvex and that $\group$ is
  quasiconvex cobounded in the case that $\nu$ does not arise from a
  density are only used to establish a decay condition for $\nu$
  \eqref{eq:integ_cond} which akin to a weakened form of upper Ahlfors
  regularity. In particular, these conditions can be dropped if this
  regularity can be established by some other means. Note that the
  uniformly quasiconvex assumption is much weaker than being geodesic.
\end{remark}

Coornaert showed in \cite{Coornaert93} that when $X$ is geodesic,
$\alpha$-quasiconformal densities only exist for $\alpha\geq
\delta(\group)$. When $X$ is not geodesic this follows from our
version of Sullivan's Shadow Lemma (\ref{lem:Shadow_lemma}). The
main significance of the above theorem is the following.

\begin{corollary}\label{cor:mu_boundary}
  If $\group$, $X$ and $\nu$ are as in the above theorem,
  then $(\Lambda,\nu)$ is a $\mu$-boundary of $\group$.
\end{corollary}

Unfortunately, even in the case of negatively curved manifolds, it
does not follow that $(\Lambda,\nu)$ is a Poisson boundary for $\mu$
despite the fact that $\nu$ is supported on all of $\Lambda$.  The
problem is twofold: two divergent sequences may actually
asymptotically represent the same $\mu$ walk, and two asymptotically
metrically convergent sequences may asymptotically represent
different $\mu$ walks. To guarantee the maximality of the above
boundary we need to connect the large scale behavior of $\mu$ to the
large scale behavior of the metric. Kaimanovich \cite{Kaimanovich00}
has formulated very general criterion for establishing maximality.
When $X$ is a CAT($-\kappa$) space, we were able to establish these
criteria in certain cases resulting in our second main result.

\begin{theorem}[Poisson Boundary]\label{thm:main2}
  Suppose that in addition to the hypotheses of Corollary
  \ref{cor:mu_boundary} we assume $X$ is a CAT($-\kappa$) space,
  $\group$ is locally compact and that a uniform neighborhood of a
  $\group$ orbit contains the convex hull of $\La$. If $\nu$ is a
  Lipschitz $\alpha$-quasiconformal measure, then there is a measure
  $\mu$ on $\group$ such that $(\Lambda,\nu)$ is a Poisson boundary
  for $(\group,\mu)$.
\end{theorem}

The most important examples to which we apply this theorem is given
by the following result which we prove in Section
\ref{sec:derivs_are_Q}.
\begin{corollary}\label{cor:Lip_examp}
  If $(X,\group)$ is as in the above theorem and $\nu$ is bounded
  Lipschitz equivalent to either the Hausdorff measure of a Busemann
  metric or a Patterson-Sullivan measure, then there is a measure $\mu$ on
  $\group$ such that $(\Lambda,\nu)$ is the Poisson boundary for
  $(\group,\mu)$.
\end{corollary}

\begin{remark}
  In each of the results above, we find solution measures $\mu$
  which are an infinite sum of atomic measures. However in
  Corollary \ref{cor:param_space} we show that, in each case, the family of
  stationizing measures $\mu$ is infinite dimensional and has
  members in any Borel measure class supported on all of $\Gamma$.
  Unfortunately, we were unable to determine whether a symmetric
  example always exists; i.e. one with $\mu(\ga)=\mu(\ga^{-1})$ for all
  $\ga\in \group$.
\end{remark}

We shall prove Theorem \ref{thm:main2} by showing that both the
finite entropy and the finite first log-moment condition are
satisfied by a carefully constructed measure given by Theorem
\ref{thm:main1}. The results of Kaimanovich mentioned above will
finish the proof. Since we do not assume that $\group$ acts
cocompactly, there is a strong constraint on the log-moment
property, but it will be satisfied whenever the limit set contains a
conical limit set of full measure. In fact, we will give criteria in
terms of the limit set, for when the usual first moment is finite.
This implies some properties of the random walk.

Before concluding the introduction, we mention an explicit
application in the following example; a setting which has seen
considerable recent interest.
\begin{example}[Fuchsian buildings]\label{exam:buildings} Throughout this example we refer to \cite{BourdonPajot00} for
details and proofs of stated facts. Let $R$ be a right angled
regular $r$-gon in the hyperbolic plane $\HH^2$. Given an $r$-tuple
of integers $(q_1, \dots, q_r)$  with $q_i\geq 2$ we assign the
cyclic group of order $q_i+1$ to the $i$-th edge and the trivial
group to the face of $R$. This gives an "orbihedron" structure to
$R$ which is developable in the sense of \cite{Haefliger91}. Its
universal developing cover $\Delta$ is a two dimensional cell
complex called a right-angle Fuchsian building. Moreover, $\Delta$
with the induced path metric is a CAT(-1) space. Let $\Ga \leqslant
Isom(\Delta)$ be the fundamental group of the orbihedron $R$, so
that $\group\backslash\Delta = R$. Since $R$ is a compact
orbihedron, the limit set and conical limit set coincide, so we may
apply Corollary \ref{cor:Lip_examp}. In particular, for a
Patterson-Sullivan measure $\nu_p$ on the boundary $\pa \Delta$ we
can find a probability measure $\mu \in P(\Ga)$ such that $(\pa
\Delta,\nu_p)$ is the Poisson boundary for $\mu$. The flexibility of
our approach allows us to handle other measures as well. Let
$G(\Ga)$ be the dual graph to the 1-skeleton of $\Delta$. If we
adjust the weight of edges by requiring the length of each edge that
crosses an edge associated to $q_i$ to be $\log(q_i)$, then we
obtain a $\delta$-hyperbolic space for some $\delta>0$. Moreover,
the Patterson-Sullivan measures for this metric graph satisfy our
conditions in Theorem \ref{thm:main1}. The associated Busemann
metric based at each chamber $p$ is called combinatorial metric
$\delta_p$. Denote the Hausdorff measure of $\delta_p$ by
$\mathcal{H}_p^{comb}$. The most important thing is that the
Hausdorff measures form a conformal density
$\frac{d(\ga^{\star}\mathcal{H}_p^{comb})}{d\mathcal{H}_p^{comb}}(z)=e^{(\tau+1)
N_{z}(p,\ga^{-1}p)}$ where $\tau$ is the critical exponent of
$\group$ with respect to $\delta_p$ and $N_z$ is the associated
Busemann function. Moreoever, this is a locally constant, hence
continuous, function on the boundary $\pa \Delta$ outside of a set
of measure $0$ (this set is the complement of the so called
tree-wall ends)\footnote{Personal communication with Marc Bourdon}.
So we can apply Corollary \ref{cor:mu_boundary} to conclude that
there exists a measure $\mu^{\pr} \in P(\Ga)$ such that
$\mu^{\pr}\star \mathcal{H}_p^{comb}=\mathcal{H}_p^{comb}$ and $(\pa
\Delta,\mathcal{H}_p^{comb})$ is a $\mu$-boundary. Since the
Patterson-Sullivan measures have the same Radon-Nikodym derivatives,
the ergodicity of the action of $\Gamma$ on $\Lambda$ implies that
the Patterson-Sullivan measures of the combinatorial metric coincide
up to a constant multiple with the Hausdorff measure with respect to
the same base chamber (see Proposition \ref{prop:coincide}).
\end{example}

We conclude with a brief outline of the paper. In Section
\ref{sec:delta_hyp} we introduce the basic tools used in working
with nongeodesic Gromov hyperbolic spaces, and we generalize some
well known estimates to this setting for later use. Section
\ref{sec:metric_measure} sets up the presentation of
Patterson-Sullivan theory in this context. Sections
\ref{sec:regularity_covers} and \ref{sec:spikes} present the
notation and background for the conditions which will be needed in
order to guarantee that a family of functions can form a positive
basis. In Section \ref{sec:basis} we present a general theorem of
independent interest which establishes when lower semicontinuous
functions can be approximated by positive sums of basis functions.
In Section \ref{sec:first_moment} this theorem is extended to give
conditions for which the functional approximations can be done with
finite first moment or log-moment. These theorems are applied to the
case of Patterson-Sullivan measures in Sections
\ref{sec:derivs_are_spikes} and \ref{sec:derivs_are_Q} where the two
main theorems are also proved. Finally we demonstrate the theorems
in the simple example of a free group in Section
\ref{sec:free_group}.

\section*{Acknowledgements}
Both authors would like to thank Alex Eskin, Alex Furman and Vadim
Kaimanovich for helpful communications, and Yehuda Shalom and Marc
Bourdon for pointing out the application in Example
\ref{exam:buildings}. The second author would also like to thank
Gregory Margulis and Carlos Kenig for insightful discussions.

\section{Background on $\delta$-hyperbolic spaces}\label{sec:delta_hyp}

We first recall the definition of a geodesic Gromov
$\delta$-hyperbolic metric space. After doing so, we recall the
general definition which may be found in \cite{BridsonHaefliger99}.
The remainder of the section is devoted to reconstructing the basic
facts about such spaces which we will need later.

\subsection{Geodesic $\delta$-hyperbolic spaces}

\begin{definition}
 We say that a geodesic metric space
$(X,d)$ is $\delta$-hyperbolic if for every geodesic triangle
$\Delta\subset X$, each side of $\Delta$ has Hausdorff distance at most
$\delta$ to the union of the other two sides. Equivalently, each side is
contained in the uniform $\delta$ neighborhood of the other two sides.
\end{definition}

\subsection{General definition}

First we define one of the basic quantities used in our asymptotic
analysis.

\begin{definition}
Let $(X,d)$ be a metric space and $x \in X$. The {\it Gromov product} of
$y,z \in X$ with respect to $x$ is defined to be
$$(y\cdot z)_x= \frac{1}{2}(d(x,y)+d(x,z)-d(y,z)).$$

It is easy to observe that $(y \cdot z)_x \leq \min\set{d(x, y), d(x,z)}$.
Now we present a definition of $\delta$-hyperbolic space without resorting
to geodesics.
\end{definition}
\begin{definition}
Let $\delta\geq 0$.  A metric space $(X,d)$ is said to be {\it
$\delta$-hyperbolic} (or {\it Gromov hyperbolic}) if
$$(x\cdot y)_w \geq \min\set{(x\cdot z)_w, (y\cdot z)_w} - \delta,$$
for all $w,x,y,z \in X$.
\end{definition}

For geodesic spaces, we could allow $w$ in the second definition
to vary along one side of the triangle $\Delta(x,y,z)$. A simple
application of the triangle inequality then shows that a
$\delta$-hyperbolic space in the second sense is
$2 \delta$-hyperbolic in the first sense. The converse is
proved in Chapter III.H of \cite{BridsonHaefliger99}.

\begin{definitions}
We say a sequence $(x_i)$ of points in $X$ {\em converges at
infinity}
\\ if $(x_i~\cdot~ x_j)_p\to~\infty$ as $i,j\to\infty$. Two sequences
converging at infinity, $(x_i)$ and $(y_i)$, are equivalent if
$(x_i\cdot y_j)_p\to\infty$ as $i,j\to\infty$. The space of all
equivalence classes of sequences converging to infinity is denoted
by $\partial X$, and is called the {\em boundary of $X$}. For a
sequence $(x_i)$ converging at infinity and an equivalence class
$x\in \pa X$, we write $x=\lim_{i\to \infty} x_i$ if $(x_i)\in x$.

Note that if $X$ is geodesic, then there is a natural bijection from this
boundary to the geodesic boundary. Now we can extend the definition of
Gromov product to $\pa X$.
\end{definitions}
\begin{definition}
Let $(X,d)$ be a $\delta$-hyperbolic space with base point $p\in X$.
Let $\bar{X} =X \cup \pa X$. We extend the Gromov product to
$x,y\in\bar{X}$ by
$$(x \cdot y)_p = \sup \liminf_{i,j \to \infty} (x_i\cdot  y_j)_p,$$
where the supremum is taken over all sequences $(x_i)$ and $(y_j)$
in $X$ such that $x =\lim_{i \to \infty} x_i$ and $y =\lim_{j \to
\infty} y_j$. (As a consequence, $\ds{(x \cdot y)_p = \sup
\limsup_{i,j \to \infty} (x_i\cdot  y_j)_p.}$)
\end{definition}

Since the metric $d$ is by definition continuous, the extended
Gromov product and the original Gromov product agree on $X\times X$.
The extended product allows us to define a topology on $\bar{X}$.
Namely, a set is closed if and only if it contains all of its limit
points.

\begin{proposition}\label{prop:product}
(\cite{BridsonHaefliger99}) Let $X$ be a $\delta$-hyperbolic space
and fix $p \in X$.
\begin{enumerate}
\item The extended product $(\  \cdot\  )_p$ is continuous on $X\times X$,
 but not necessarily on $\bar{X}\times\bar{X}$.

\item In the definition of $(x\cdot y)_p$, if we have $x\in X$
(resp. $y\in X$), then we may always take the respective sequence
to be the constant value $x_i=x$ (resp. $y_j=y$).

\item\label{item:seqexist} For all $x, y \in \overline{X}$ there
exist sequences $(x_n)$ and $(y_n)$ such that $\ds{x =\lim_{n \to
\infty} x_n}$ and $\ds{y =\lim_{n \to \infty} y_n}$ and $\ds{(x
\cdot y)_p = \lim_{n \to \infty}(x_n \cdot y_n)_p}$.

\item\label{item:still_ultrametric} For all $x,y,z \in \overline{X}$ by taking limits we still
have $$(x \cdot y)_p \geq \min\set{(x \cdot z)_p, (y, \cdot
z)_p}-2\delta.$$

\item\label{item:limits_pinched} For all $x, y \in \pa X$ and all sequences $(x_i^{\pr})$ and
$(y_j^{\pr})$ in $X$ with $\ds{x =\lim_{i \to \infty} x_i^{\pr}}$
and $\ds{y =\lim_{j \to \infty} y_i^{\pr}}$,
$$(x \cdot y)_p - 2\delta \leq \liminf_{i, j \to \infty}(x_i \cdot
y_j)_p \leq (x \cdot y)_p.$$
\end{enumerate}
\end{proposition}

Recall that a $(\lambda,C)$-quasi-isometric embedding between two
metric spaces $(X,d_X)$ and $(Y,d_Y)$ is a map $f:X\to Y$ such
that for all $x,y\in X$,
$$\frac{1}{\la} d_X(x,y)-C\leq d_Y(fx,fy)\leq \la d_X(x,y)+C.$$
A $C$-quasigeodesic in $X$ is the image of a
$(1,C)$-quasi-isometric embedding of $\R$ into $X$.

\begin{definition}
We say that a metric space $X$ is {\em ($C$-)uniformly
quasigeodesic} if there is a $C\geq 0$ such that for any two
distinct points in $X$ there is a $C$-quasi-geodesic joining them.
In particular, a $0$-uniformly quasigeodesic space is geodesic.
\end{definition}

\begin{definition}
Let $(X,d)$ be a $\delta$-hyperbolic space with base point $p\in X$.  Fix
$\ep>0$. We consider the following measure of separation of the points in
$\overline{X}$
$$d_p^{\ep} ( x,y) = e^{-\ep(x \cdot y)_p},$$
for $x,y \in \overline{X}$. Denote $d_p^1 =d_p$.
\end{definition}

\subsection{Metrics on $\partial X$}

Now we recall the existence of a compatible metric on $\partial X$.
(See \cite{BridsonHaefliger99} or \cite{GhysDeLaHarpe88} for the
nongeodesic case.)
\begin{proposition}\label{prop:metric}
\cite{Gromov87} Let $(X,d)$ be a $\delta$-hyperbolic space.  If
$0<\ep \leq \frac{\log 2}{4\delta}$, then there exists a metric
$\delta_p^\eps$ on $\pa X$ so that
$$(3-2e^{2\delta\ep})d_p^{\ep}(z, y) \leq \delta_p^\eps(z, y) \leq
d_p^{\ep}(z, y),$$ for all $z, y \in \pa X$.
\end{proposition}

Since the Gromov product is nonnegative, the above proposition
implies $\op{diam}(\pa X)\leq 1$ in the metric $\delta_p^\eps$ for
any $p\in X$ and $0<\ep \leq \frac{\log 2}{4\delta}$. Moreover, we
have

\begin{proposition}[6.2 of \cite{BonkSchramm00}]
With respect to $\delta_p^\ep$, the boundary $\pa X$ is complete.
\end{proposition}

Note that $\pa X$ need not be compact. For instance, the boundary of
$\HH^\infty$ is the Hilbert sphere. By taking a tree whose level $n$
leaves form a $1/n$ net in the $n$-dimensional sphere, one obtains a
locally compact example with the same boundary. On the other hand,
if all of the closed balls in $X$ are compact, then we say $X$ is
{\em proper}. This condition is much stronger, and we omit the
(straightforward) proof of the following.

\begin{proposition}
  If $X$ is a proper $\delta$-hyperbolic space, then $\pa X$ is
  compact and $\op{Isom}(X)$ is locally compact.
\end{proposition}
With some difficulty we avoid using this assumption on $X$. In fact,
we don't even assume $\group$ is locally compact until Theorem
\ref{thm:main2} and its corollaries. The point is to allow groups in
the first theorem which come from infinite dimensional
constructions.

Recall that a {\it $C$-quasimetric} is a function $d:X\times X\to
\R$ satisfying all of the properties of a metric except for the
triangle inequality which is substituted by the condition that
$d(x,y)\leq C(d(x,z)+d(z,y))$ for all $x,y,z\in X$ and some $C\geq
1$. Proposition \ref{prop:metric} implies that $d_p^{\ep}(z, y)$
is a $\frac{1}{(3-2e^{2\delta\ep})}$-quasimetric which can be used
in most computations instead of the more complicated metric.
However, the $s$ power of a $C$-quasimetric is a
$C2^{s-1}$-quasimetric, so for any $s>0$ and any $0<\ep<
\frac{\log 2}{4\delta}$, $d_p^s$ is a
$\frac{2^{\frac{s}{\ep}-1}}{3-2 e^{2\delta\ep}}$-quasimetric.
Whenever $\delta>\frac{\left( 4 + 3\,{\sqrt{2}} \right) \,\log
(2)}{4 s}$, the choice $\ep=\frac{\log 2}{4\delta}$ is optimal in
the valid interval giving that $d_p^s$ is a $\left( \frac{3}{2} +
{\sqrt{2}} \right) e^{4\,\delta s}$-quasimetric.

\begin{lemma} \label{lem:value_on_balls}Let $(X,d)$ be a $\delta$-hyperbolic space.  Fix two points $p, q
\in X$.  Then for all $z, y \in \pa X$  and $C>0$ such that
$$d_p(z, y) \leq e^{C}e^{- d(p,q)}\leq e^{C}\min\set{e^{- (z\cdot q)_p},e^{- (y\cdot q)_p}},$$
we have $|(z\cdot q)_p - (y\cdot q)_p|\leq 2\delta+C$. (Recall that
$\max\set{(y\cdot q)_p,(z\cdot q)_p}\leq d(p,q)$.)
\end{lemma}

\medskip
\begin{proof}
 Without loss of generality assume that $(y\cdot q)_p
\geq (z\cdot q)_p$. Since $d_p(z, y) =
e^{-(z\cdot y)_p}$, we have
$$(z\cdot y)_p \geq d(p,q) -C \geq \max\set{(z\cdot
q)_p,(y\cdot q)_p}-C = (y \cdot q)_p -C.$$ By Item
\ref{item:still_ultrametric} in Proposition \ref{prop:product} we
have
$$\begin{aligned}
(z\cdot q)_p &\geq \min\set{(y\cdot q)_p, (z\cdot
y)_p}-2\delta = \min\set{(y\cdot q)_p, (y \cdot q)_p -C}-2\delta \\
& =(y \cdot q)_p -2\delta C.\end{aligned}$$ This proves the
lemma.

\end{proof}

Now we present some estimates we shall need later for how $d_p(x,y)$
varies on $\overline{X}$.
\begin{lemma} \label{lem:compare_to_distance} Let $(X,d)$ be a $\delta$-hyperbolic space and $p, q \in X$.
 Let $U_p(q) = \sup_{z \in \pa X} \set{(z \cdot q)_p}$. Let $z_{p,q}$ be any point in $\pa X$ such that
$(z_{p,q}\cdot q)_p \geq U_p(q)-\delta$.
\begin{enumerate}
\item[a)] For all $ y \in \pa X$, we have
$$ e^{(y\cdot q)_p} \leq \frac{e^{3\delta}}{d_p(z_{p,q}, y)}.$$

\item[b)] If $d_p(z_{p,q}, y) \ge e^{- U_p(q)}$, then
$$e^{(y\cdot q)_p} \geq \frac{e^{-3\delta}}{d_p(z_{p,q},
y)}.$$
\end{enumerate}
 \end{lemma}

\medskip
\begin{proof}
For part a), recall that $$(z_{p,q}\cdot y)_p \ge
\min\set{(z_{p,q}\cdot q)_p, (y\cdot q)_p}-2\delta. $$

Since $(z_{p,q}\cdot q)_p \geq U_p(q)-\delta$ and $U_p(q)\geq (y\cdot
q)_p$, we obtain that
$$(z_{p,q}\cdot y)_p \ge
\min\set{(y\cdot q)_p -\delta, (y\cdot q)_p }-2\delta =(y\cdot q)_p
-3\delta. $$

For b), the condition $d_p(z_{p,q}, y) \ge e^{- U_p(q)}$ implies that
$(z_{p,q}\cdot y)_p \leq U_p(q)$. Since $(z_{p,q}\cdot q)_p \geq
U_p(q)-\delta$, we have

$$\begin{aligned} (y \cdot q)_p &\geq \min\set{(z_{p,q}\cdot q)_p,
(z_{p,q}\cdot y)_p}  -2\delta \geq \\&\geq \min\set{(z_{p,q}\cdot
y)_p-\delta, (z_{p,q}\cdot y)_p} -2\delta = (z_{p,q}\cdot y)_p-
3\delta. \end{aligned}$$ This proves the lemma.
\end{proof}

\subsection{Isometric actions on $\bar{X}$.}\label{subsec:isometric}
Assume that a group $\group$ acts by isometries on $X$.  Fix $p \in
X$.

One easily observes that $(\ga x \cdot \ga y)_{\ga p} = (x \cdot
y)_p$ for all $x, y \in X$ and $\ga \in \group$. The same is true for
all $x,y\in\overline{X}$. For if $(x_n)$ and $(y_n)$ are two
sequences such that
$$(x\cdot y)_p = \lim_{n \to \infty}(x_n, y_n)_p$$ then $$ \lim_{n
\to \infty}(\ga x_n, \ga y_n)_{\ga p} = (\ga x, \ga y)_{\ga p}.$$

For a similar proof of the following in the manifold case, see
\cite{Yue96}.
\begin{lemma}\label{lem:ratio_in_delta_space}
For all $x,y\in \pa X$ and all $p,q\in X$ we have

$$e^{-d(p, q) -2\delta} e^{((x\cdot q)_p+
(y\cdot q)_p )} \leq \frac{d_q(x, y)}{d_p(x,y)} \leq e^{-d(p, q)+2\delta}
e^{((x\cdot q)_p+ (y\cdot q)_p )}.$$
\end{lemma}

\medskip
\begin{proof}
  Let $x, y \in \pa X$.  Assume that $(x_n)$ and $(y_n)$
are two sequences in $X$ such that $\lim_{n \to \infty}x_n = x$ and
$\lim_{n \to \infty}y_n = y$ and the limit $\lim_{n\to \infty}
(x_n\cdot y_n)_p$ exists. By item \ref{item:seqexist} in Proposition
\ref{prop:product}, the sequence $(x_n)$ and $(y_n)$ may be chosen
so that $\lim_{n\to \infty}(x_n\cdot q)_p = (x\cdot q)_p$ and
$\lim_{n\to \infty}(y_n\cdot q)_p = (y\cdot q)_p$.

By hyperbolicity, we have $$(x\cdot y)_p - 2\delta \leq \lim_{n\to
\infty} (x_n\cdot y_n)_p \leq (x\cdot y)_p,$$ and
$$( x\cdot y)_q - 2\delta \leq \liminf_{n\to \infty}
( x_n\cdot y_n)_q \leq (x\cdot y)_q.$$

Therefore
$$ \liminf_{n\to \infty}((x_n\cdot y_n)_q
-( x_n\cdot y_n)_p)- 2\delta \leq ( x\cdot y)_p-(\ga
x\cdot \ga y)_p \leq \liminf_{n\to \infty}((x_n\cdot  y_n)_q -( x_n \cdot y_n)_p) + 2\delta.$$

Now observe that $$( x_n\cdot y_n)_q -(  x_n\cdot  y_n)_p = -(x_n\cdot q )_{p}- (y_n\cdot q)_{p} + d(p, q).$$

So we obtain, $$-(x\cdot q )_{p}- (y\cdot q)_{p} + d(p, q)-2\delta\leq ( x\cdot y)_q-(x\cdot y)_p \leq
-(x\cdot q )_{p}- (y\cdot q)_{p} + d(p, q)+
2\delta.$$

Therefore,
$$\frac{d_q( x, y)}{d_p(x,y)}  =
\frac{e^{-(x\cdot y)_q}}{e^{-(x\cdot y)_p}} =
e^{-((x \cdot y)_q + (x\cdot y)_p)} \leq e^{
2\delta}e^{- d(p, q)} e^{((x\cdot q)_p+
(y\cdot q)_p )},$$ and

$$\frac{d_q(x, y)}{d_p(x,y)} \geq e^{-
2\delta}e^{- d(p, q)} e^{((x\cdot q)_p+
(y\cdot q)_p )}.$$

\end{proof}

\begin{corollary}
For all $x,y\in \pa X$ and any $p,q \in X$ we have

$$e^{(d(p, q) - 2\delta)} e^{-((x\cdot p)_q+
(y\cdot p)_q )} \leq \frac{ d_q(x, y) }{ d_p(x,y) } \leq e^{ (d(p,
q)+2\delta)} e^{((x\cdot p)_q + (y\cdot p)_q )}.
$$
\end{corollary}

\medskip
\begin{proof} Let $\set{x_n}$ be a sequence in $X$ such that
$\lim_{n \to \infty} (x_n\cdot q)_p = (x\cdot q)_p$. Since
$(x_n\cdot p)_q = d(q, p) - (x_n\cdot q)_p$ and $\liminf_{n \to
\infty} (x_n \cdot p)_q \leq (x\cdot p)_q$, we have
$$d(p,q) -
(x\cdot q)_p \leq (x \cdot p)_q.$$ Similarly, we have $d(p,q) -
(x\cdot p)_q \leq (x \cdot q)_p.$  So $d(p,q) - (x\cdot q)_p = (x
\cdot p)_q.$ This proves the corollary.
\end{proof}

\medskip
Before concluding this section, we will introduce a few more
definitions and prove the Shadow Lemma for $\delta$-hyperbolic
spaces.

\begin{definitions} For fixed $x,p\in X$, let $\Lambda=\Lambda(\group)$
  denote the subset of $\partial X$ consisting of all asymptotic equivalence
  classes, with respect to $(\ \cdot\ )_p$, of sequences of the form
  $\(\gamma_i x\)$ for $\gamma_i\in\group$. The set $\Lambda$ is called the
  {\em limit set of $\group$} and is sometimes denoted as $\partial \group$.

  The {\em radial (or conical) limit set of $\group$}, denoted by
$\Lambda_r$, is the subset of $\Lambda$ such that $z\in \Lambda_r$
if and only if there is a constant $C>0$ and a sequence $g_i\in
\group$ with $\(g_i x\)$ converging to $z$ such that $d(g_ix,p)-(g_i
x \cdot z)_p\leq C$. Note that  by Proposition \ref{prop:product},
this is equivalent to $(p\cdot z)_{g_i x}\leq C'$ for some constant
$C-2\delta<C'<C+2\delta$. In the case of a geodesic space, we may
express this condition by saying that the orbit subsequence $(g_i
x)$ must stay within distance $C$ of the geodesic between $p$ and
$z$.
\end{definitions}
It is easy to see that these definitions do not depend on the choice of
$p$ or $x$, so long as we are free to change the constant $C$.

\begin{definitions}
For a geodesic space $X$, the geodesic hull (sometimes called the
Gromov envelope) of $\Lambda$, $\op{GH}(\Lambda)$ is the union of
all geodesics in $X$ with both endpoints in $\Lambda$. The convex
hull of $\Lambda$, $\op{CH}(\Lambda)$, is the smallest subset of $X$
containing $\op{GH}(\Lambda)$ with the property that {\emph every}
geodesic segment between any pair of points $p,q\in\op{CH}(\Lambda)$
also lies in $\op{CH}(\Lambda)$. If $X$ is $C$-uniformly
quasigeodesic, then we let $\op{GH}(\Lambda)$ be the union of all
$C$-quasigeodesics with both endpoints in $\Lambda$. Here we assume
this constant $C$ has been chosen once and for all.

We say that $X$ is {\it quasiconvex cocompact (resp. quasiconvex
cobounded)} with respect to the action of $\group$ if
$\op{GH}(\Lambda)/\group$ is compact (resp. $\op{GH}(\Lambda)$
lies in the uniform neighborhood of some $\group$-orbit).
Similarly, the action of $\group$ on $X$ is {\it convex cocompact}
if $\op{CH}(\Lambda)/\group$ is compact.
\end{definitions}

It is important to note that $\op{CH}(\Lambda)$ is strictly larger than
$\op{GH}(\Lambda)$ even for (most) actions of surface groups on $\HH^3$.
Nevertheless, for CAT(-1)  spaces (see Section \ref{sec:metric_measure}),
the notion of quasiconvex cocompactness is equivalent to convex
cocompactness. The quasiconvex cocompact actions are the simplest families
of examples where $\Lambda_r$ and $\Lambda$ coincide.

\begin{definitions}
We say that a $\delta$-hyperbolic space $X$ is {\em upper (resp.
lower) Gromov product bounded from above (resp. below)} if there
exists a $p\in X$ and a constant $C$ such that
$$\sup_{z \in \pa X}(z \cdot q)_p \geq d(p,q)-C, \quad\left(\inf_{z
\in \pa X}(z \cdot q)_p \leq C,\right)$$ for all $q\in X$. The space
$X$ is {\em Gromov product bounded} if it is both upper and lower
Gromov product bounded. The $\delta$-hyperbolic space $X$ is {\it
weakly Gromov product bounded} if there exists a constant $C$ and
$p\in X$ such that
$$\inf_{z \in \pa X}(z \cdot q)_p+ \sup_{z \in \pa X}(z \cdot q)_p
-d(p, q) \leq C,$$ for all $q\in X$. Lastly we say $X$ is
(weakly,upper,lower) Gromov bounded with respect to $\group$ if we
only require the corresponding condition to hold for $q\in
\group\cdot p$. A straightforward application of the triangle
inequality shows that if any of these conditions holds for one $p\in
X$, then it also holds for every $p\in X$ for an appropriate choice
of $C$.
\end{definitions}

If $\Gamma$ acts cocompactly, then any $p$ would do, and we could
choose $C$ independently of $p$. Also, if $X$ is complete and
geodesic, then we can set $C=0$.

Being Gromov product bounded from above (respectively below) morally
means that for each $q\in X$ there is a point $z\in \pa X$ which
assumes the role of the forward (backward) endpoint of a geodesic
through $p$ and $q$, even though no such geodesic need exist. In
particular, if for every $q\in X$, $p$ and $q$ are always connected
by a discrete $C$-quasigeodesic, then $X$ is Gromov product bounded.

\subsection{Quasiconformal Measures on $\Lambda$}
We will need to consider certain classes of measures which transform
nicely under the $\group$ action.
\begin{definitions}
Let $\nu$ be a Borel measure on $\pa X$. We say that $\nu$ is {\it upper
(resp. lower) $\alpha$-quasiconformal} for $\group$ if for some $p\in X$,
$$\frac{d \ga_{\star}\nu}{d\nu}(z) \leq C e^{-\alpha d(p,
\ga^{-1} p)} e^{2\alpha(z\cdot \ga^{-1}p)_p},\quad \left(\frac{d
\ga_{\star}\nu}{d\nu}(z) \geq C^{-1} e^{-\alpha d(p, \ga^{-1} p)}
e^{2\alpha(z\cdot \ga^{-1}p)_p},\right)$$ for all $\ga \in \group$
and some $C>1$, (where $\ga_{\star}\nu = \ga^{-1}\nu$). If both
inequalities are satisfied, then we say $\nu$ is {\it
$\alpha$-quasiconformal.} If in addition for each $\ga\in\group$,
$\frac{d \ga_{\star}\nu}{d\nu}$ is continuous, then $\nu$ is said to
be a {\it continuous  $\alpha$-quasiconformal} measure.

Suppose $\nu$ is $\alpha$-quasiconformal and the bounded function
$R_\ga$ given by
$$R_\ga(z)=\frac{d \ga_{\star}\nu}{d\nu}(z)e^{\alpha d(p, \ga^{-1}
p)-2\alpha (z\cdot \ga^{-1}p)_p}$$ satisfies for some $\eps>0$ and
each $x\in \Lambda$,
$$\sup_{\set{y|0<d(x,y)<e^{-\eps
d(p,\ga^{-1}p)}}}\frac{\abs{R_\ga(x)-R_\ga(y)}}{d^\eps_p(x,y)}\leq
\frac{C}{e^{-\eps d(p,\ga^{-1}p)}}$$ for some $C$, independent of
$\ga$. Then  we say that $\nu$ is {\em Lipschitz
$\alpha$-quasiconformal.} All of the above definitions are
independent of $p$.
\end{definitions}

Note that if $\nu$ is a continuous $\alpha$-quasiconformal measure
with constant $C$ and $f:\D X \to \R$ is a uniformly positive and
bounded continuous function, then $f\nu$ is also a continuous
$\alpha$-quasiconformal measure with constant $\frac{\sup f}{\inf f}
C$.

The reason for the specific Lipschitz estimate in the definition of
Lipschitz $\alpha$-quasiconformal should become clear in Section
\ref{sec:spikes}. For our purposes, assuming that $R_\ga$ is
Lipschitz uniformly in $\ga$ would exclude the standard examples,
while just assuming the Lipschitz constant depends arbitrarily on
$\ga$ is too weak.

\begin{definitions}
A family of finite Borel measure $\set{m_x}_{x \in X}$ of $\pa X$ is
called an {\it $f$-density of $\group$} for a measurable function $f: X
\times X\times \pa X \to
\rr$ if for all points $x,y \in \pa X$ the measure $m_x, m_y$ are
equivalent with Radon-Nikodym derivatives
$$\frac{d m_y}{d m_x}( z) = f(x,y,  z),$$
for $m_x$-a.e. $ z \in \pa X$ and for all $\ga \in \group$,
$f(\ga^{-1} x, x, \ga^{-1}  z)=f(x, \ga x,  z)$. This is
equivalent to the condition that $\gamma_\star m_x=m_{\ga x}$. An
$f$-density is called {\it continuous} if $f$ is continuous. An
$f$-density is called bounded if the mass $\norm{m_x}=m_x(\pa X)$
is bounded for all $x\in X$.
\end{definitions}

In what follows it is common to express quantities in terms of Busemann
functions instead of the Gromov product, especially when $X$ is a ``nice''
space. To this end we introduce this notation here.
\begin{definition}
For $x,y \in X$ and $ z \in \pa X$, we define the {\it Busemann function}
$\rho_{x,z}$ by
$$\rho_{x, z}(y)= 2(x\cdot z)_y -d(x,y).$$
\end{definition}

Since $d$ is $\group$-invariant, we have
$$\rho_{\ga x,\ga  z}(\ga y)= \rho_{x, z}(y),$$
for all $x,y \in X$, $  z \in \pa X$ and $\ga \in \group$.

\begin{definition}
If $f(x,y, z) = e^{- \alpha \rho_{x, z}(y)}$ for all $(x,y,z)\in X
\times X\times \pa X$ and some $\alpha\in\R$, then the corresponding
$f$-density is called an {\it $\alpha$-conformal density} of
$\group$.
\end{definition}

Note that a single member $m_p$ of an $f$-density on $\Lambda(\group)$ is
an $\alpha$-quasiconformal density of $\group$ if and only if there exists
some constant $C\geq 1$ such that for all $\ga \in \group$ and $z\in \pa
X$,
$$C^{-1} e^{- \alpha \rho_{p,
 z}(\ga^{-1}p)} \leq f(p,\ga p,z) \leq C  e^{- \alpha
\rho_{p,  z}(\ga^{-1}p)}.$$ If this holds for every member of the
$f$-density for a uniform $C$, then we say that the $f$-density is an {\em
$\alpha$-quasiconformal density} of $\group$.

Now we turn our attention to two examples of $\alpha$-quasiconformal
densities.
\subsubsection{Hausdorff measures}

Let $(Z, \delta)$ be any metric space
 and $D\geq 0$ be a nonnegative constant.  Let $A$ be a subset of
 $Z$.  For each $\ep\geq 0$ we define
 $$H_{\ep}^D(A)\triangleq \inf\set{\sum_{j=1}^{\infty} r_j^D \, |\, A\ssu \cup_j B(x_j, r_j), r_j \leq \ep, x_j \in A}$$
where the infinum is taken among all coverings of $A$ by balls of radius
no more than $\ep$.  The limit measure $\mathcal{H}_{\delta}^D(A) =
\lim_{\ep
\to 0} H_{\ep}^D(A)$ is called the $D$-dimensional Hausdorff measure of
$A$. The Hausdorff dimension $\dim_H(A)$ is defined to be
$$\dim_H(A) \triangleq \inf\set{D\,:\, H_{\delta}^D(A) =0} = \sup\set{D\,:\,
H_{\delta}^D(A) = \infty}.$$

In the definition of $\mathcal{H}_{\delta}^D$ we could have used any
quasimetric bilipschitz to $\delta$ and obtained the same measure.
Now we return to the $\delta$-hyperbolic setting. As a consequence
of Item \ref{item:limits_pinched} in Proposition \ref{prop:product}
for all $x,y\in X$ and $z\in \pa X$, the Busemann functions satisfy
$$\abs{\rho_{x,z}(y)+\rho_{y,z}(x)}\leq 4\delta.$$
Using this, we can now extend Proposition 4.3 in \cite{Coornaert93}
to this more general class of $\delta$-hyperbolic spaces.

\begin{proposition}\label{prop:quasiconformal} Suppose $X$ is a
$\delta$-hyperbolic space, and suppose $A \ssu \pa X$ is a
$\group$-invariant Borel set with
$0<\mathcal{H}_{\delta_p^\ep}^{\alpha/\ep}(A)<\infty$ for some $p\in
X$. Then $\set{\mathcal{H}_{\delta_x^\ep}^{\alpha/\ep}}_{x\in X}$ is
an $\alpha$-quasiconformal density of $\group$ on $\pa X$.
\end{proposition}

\begin{proof}
The quasi-antisymmetric property of Busemann functions described
above allows us to rewrite the conclusion of Lemma
\ref{lem:value_on_balls} to obtain Lemma 2.2 of \cite{Coornaert93},
except using our definition of the Busemann function. Armed with
this lemma, the proof of the proposition is identical to that of
Proposition 4.3 in \cite{Coornaert93} in the geodesic setting.
\end{proof}
\subsubsection{Patterson-Sullivan measures}\label{subsec:Patt-Sull}

Assume that $\group$ is discrete and acts properly dicontinuously on $X$.
For two points $x, y \in X$ and for any real number $s>0$, we consider the
Poincare series $$g_s(x,y) = \sum_{\ga \in \group} e^{-sd(x, \ga y)}.$$
Let $S_k$ be the number of the orbit points  $\group y$ in annulus $B(x,
k+\frac{1}{2})\\ B(x, k-\frac{1}{2})$. Then $g_s(x,y)$ is proportional to
$\sum_{k=0}^{\infty}S_k e^{-sk}$.

We define $$\delta(\group) \triangleq \limsup_{k\to
\infty}\frac{1}{k} \log S_k$$ to be the critical exponent of
$\group$ which depends on the action of $\group$ as much as $\group$
itself. Then $g_s(x,y)$ diverges for $s <\delta(\group)$ and
converges for $s
>\delta(\group)$.  So we can consider the family of measures
$$\nu_x^s = \frac{1}{g_s(y,y)}\sum_{\ga \in \group} e^{-sd(x, \ga y)}\delta_{\ga y}, \hskip .5in  s>\delta(\group),$$
where $\delta_{\ga y}$ is the Dirac mass at $\ga y$. Since $$d(x,
\ga y)-d(x,y) \leq d(y, \ga y) \leq d(x, \ga y)+d(x,y),$$ we can
easily see that $g_s(\ga x,\ga^\pr y)=g_s(x,y)$ for all
$\ga,\ga^\pr\in\group$ and $$e^{-sd(x,y)}g_s(x,y) \leq g_s(y,y) \leq
e^{sd(x,y)}g_s(x,y).$$ Thus $\set{\nu_x^s}_{\delta(\group)<s}$ is a
family of finite measures on $X$ with total mass bounded by $$e^{-s
d(x,\group\cdot y)}\leq \norm{\nu_x^s}\leq e^{-s d(x,\group\cdot
y)}.$$ If $\liminf_{s\to \delta(\group)^+} g_s(y,y) < \infty$, then
we introduce a weighting function in front of the exponential
factors above so that $\liminf_{s\to \delta(\group)^+} g_s(y.y) =
\infty$. This can be done so that the essential properties of
$\nu_x^s$ are preserved (see \cite{Coornaert93} following
\cite{Patterson76}). Let $\nu_x = \lim_{s_j \to \delta(\group)^+}
\nu_x^{s_j}$ be a weak limit in the space of uniformly bounded
measures on $X\cup\pa X$. There may be many distinct limit measures,
however we always fix one such limit. Since $\lim_{s_i\to
\delta(\group)^+}g_s(y.y) = \infty$, the measure $\nu_x$ is
concentrated on the cluster points of the orbit $\group y$, i.e.,
the limit set $\Lambda(\group)$. From the construction it is easy to
see that for any other point $y \in X$ the limit $\nu_{y} =
\lim_{s_j \to \delta(\group)^+} \nu_{y}^{s_j}$ also exists and that
$\ga_{\star}\nu_x = \nu_{\ga x}$ for all $\ga\in\group$.

Suppose $\group$ is not discrete or does not act properly
discontinuously, but it admits a left-invariant infinite measure (or
mean) $\eta$ with the property that for any compact set $K$ and any
point $x\in X$, $\eta(\group\cdot x \cap K)<\infty$. By converting
the above sums to $\eta$-integrals we can still obtain a critical
exponent $\delta(\group)$ as well as the family of finite
Patterson-Sullivan measures $\nu_x$.

In the case when $X$ is complete and geodesic, \cite{Coornaert93}
shows that the Radon-Nikodym derivative at $ z \in \pa X$ satisfies
$$C^{-1} e^{-\delta(\group) \rho_{x, z}(y)}\leq\frac{d \nu_{y}}{d
\nu_x}( z)\leq C e^{-\delta(\group) \rho_{x, z}(y)}$$
for a constant $C\geq 1$ depending only on the hyperbolicity
constant $\delta$. In short we have the following.
\medskip

\begin{corollary}
If $X$ is complete and geodesic, then the family $\set{\nu_x}_{x\in
X}$ of Patterson-Sullivan measures is an
$\delta(\group)$-quasiconformal density.
\end{corollary}

These measures are not automatically a continuous
$\alpha$-quasiconformal density.  However, if the Gromov product
extends continuously to the boundary then both Busemann functions
and a simple computation using the above formula for $\nu^s$ shows
continuity for the Radon-Nikodym derivatives (see Coornaert
\cite{Coornaert93}). In the next section we shall see that, under
still stronger conditions, these measures become $\alpha$-conformal.

\subsection{The Shadow Lemma}

Since our $\delta$-hyperbolic space is not necessarily a geodesic space,
we cannot define a shadow for balls as in the Shadow Lemma proved by
Sullivan.  However, we can talk about the balls with respect to Busemann
distance, even though it is not in general a metric. So we denote
$$U_{p}(q)=\sup_{z \in
\pa X}((z\cdot q)_p)\quad\text{and}\quad
L_{p}(q)=\inf_{z \in
\pa X}((z\cdot q)_p).$$
We abuse the notation and denote $U_{p, \ga}= U_{p}(\ga^{-1} p)$ and
 $L_{p, \ga}= L_{p}(\ga^{-1} p)$.

For every $p$ and $q$ we
also fix a choice of points
$z_{p,q}^+,z_{p,q}^-\in\partial X$ such that
$$(z_{p, q}^+ \cdot
q)_p \geq U_{p}(q)-\delta \quad\text{ and }\quad
(z_{p,q}^- \cdot q)_p \leq L_{p}(q)+\delta.$$
And for an isometry $\ga$ we set $q=\ga^{-1}p$ and $z_{p,q}^+=z_{p,\ga}^+$ and $z_{p,q}^-=z_{p,\ga}^-$.

Lastly, we define $O_p(\gamma,D)=O_p(\gamma^{-1}p,D)$ where for
$p,q\in X$ and $D\geq 0$,
$$O_p(q, D) := \set{ y \in \pa X \, : \, d_p
(z_{p,q}^+, y) \leq e^{-U_{p}(q)+D}}.$$

\medskip

\begin{remark}\label{rem:limit_set}
In the above notation the radial limit set can be written
explicitly as
$$\La_r (\Ga) = \bigcup_{D=0}^{\infty}\bigcap_{n=1}^{\infty}
\bigcup_{\substack{\ga\in\group \\ d(p,\ga p)\geq n}}
O_p(\ga,D).$$ Again, the right hand side is independent of the
point $p$.
\end{remark}

\begin{lemma} (Shadow Lemma) \label{lem:Shadow_lemma}
Let $X$ be a $\delta$-hyperbolic space and $p\in X$.  Assume $\nu$
(or $\nu_p$) represents a finite measure that does not consist of
a single atom.

\begin{itemize}
\item[a)] If $X$ is weakly Gromov product bounded with respect to
$\group$, and $\nu$ is upper $\alpha$-quasiconformal, then there
exist numbers $\beta>0$ and $D_0\geq 0$ such that for all $D>D_0$
and $\gamma\in \group$ we have
$$\nu(O_p(\ga, D)) \ge \beta e^{-\alpha \,U_{p,\ga}}\quad
 \left(\text{in particular,}\quad \nu(O_p(\ga, D)) \ge \beta
 e^{-\alpha d(p,\gamma^{-1} p)-\alpha L_{p,\ga}}\right).$$
\item[b)] If $X$ is upper Gromov product bounded with respect to
$\group$, and $\nu$ is lower $\alpha$-quasiconformal, then there
exists a number $\beta>0$ such that for all $D\geq 0$ and $\gamma\in
\group$ we have
$$\nu(O_p(\ga, D)) \leq \beta e^{-\alpha U_{p, \ga}} e^{2\alpha
D}\quad
 \left(\text{in particular,}\quad \nu(O_p(\ga, D)) \leq \beta
 e^{-\alpha d(p,\gamma^{-1} p)+2\alpha D}\right).$$
\item[b')] If $X$ is upper Gromov product bounded and
$\set{\nu_x}_{x\in X}$ is a lower $\alpha$-quasiconformal density, then there
exists a number $\beta>0$ such that for all $D\geq 0$ and $p,q\in X$ we have
$$\nu_p(O_p(q, D)) \leq \beta\norm{\nu_q} e^{-\alpha U_{p, q}} e^{2\alpha
D}\quad
 \left(\text{in particular,}\quad \nu_p(O_p(q, D)) \leq
 \beta\norm{\nu_q}
 e^{-\alpha d(p,q)+2\alpha D}\right).$$
\end{itemize}
\end{lemma}

\begin{proof}
Without loss of generality we may assume $\nu$ is a probability
measure and let the $C$ denote the quasiconformal constant for
$\nu$.

For a), fix $\ep >0$ such that $d_p^{\ep}$ is a quasimetric. Since
$\nu$ does not consist of a single atom, there exists a number $r>0$
and $a<1$ such that
$$\nu(B(y, r))<a.$$ By Lemma \ref{lem:compare_to_distance} for all
$y \notin O_p(\ga, D)$ we have
$$e^{(y\cdot q)_p} \leq
\frac{e^{3\delta}}{d_p(z_{p,q}^+, y)} \leq
e^{3\delta}e^{(U_{p}(q)-D)}.$$

Let $K$ be the constant for which $X$ is weakly Gromov product bounded
relative to the point $p\in X$. Together with Lemma
\ref{lem:ratio_in_delta_space} we have,
$$\begin{aligned}
d_q(z_{p,q}^-, y) &\leq d_p(z_{p,q}^-,  y) e^{2\delta}e^{- d(p,
q)} e^{((z_{p,q}^-\cdot q)_p+ ( y \cdot q)_p )}\leq
\\&\leq e^{ 2\delta}e^{-d(p, q)} e^{(L_{p}(q)+\delta)}
e^{3\delta}e^{(U_{p}(q)-D)} \leq e^{6\delta}e^{- (D-K) }.
\end{aligned}$$
Thus $\op{diam}_{d_q}(O_p(q, D)^c) \leq e^{6\delta\ep}e^{- \ep (D-K) }$.
Therefore setting $q=\ga^{-1}p$, there exists $D_0$ such that for all
$D\geq D_0$,
$$\nu(\ga O_p(\ga, D)^c) \leq a \quad \text{ or equivalently}\quad
\nu(\ga O_p(\ga, D)) \geq 1-a.$$  Now we have
$$\begin{aligned}
1-a \leq & \nu(\ga O_p(\ga,D)) = \ga^{-1}_{\star}\nu(O_p(\ga,D))
\leq \int_{O_p(\ga,D)} C e^{-\alpha d(p, q)}e^{2\alpha ( y\cdot
q)_p} d\nu( y) \\ &\leq C e^{\alpha (2U_{p,\ga}-d(p,\ga^{-1}p))}
\nu(O_p(\ga,D)).
\end{aligned}$$

Thus we obtain that
$$  \nu(O_p(\ga,D)) \geq (1-a)C e^{\alpha(d(p,\ga^{-1}p)-U_{p,\ga})}
e^{-\alpha ( U_{p,\ga})}.$$

To finish the proof recall that weak Gromov product bounded implies $U_{p,
\ga} \leq d(p, \ga^{-1}p)-L_{p,\ga}+K$ for some constant $K>0$.

For part b), first note that for any $q\in X$ and every $ y \in
O_p(q, D)$ we can use Lemma \ref{lem:value_on_balls} to obtain
$$|( y\cdot q)_p - (z_{p, q}^+\cdot q)_p|
\leq 2\delta+D.$$  This implies that
$$e^{( y\cdot q)_p} \geq e^{-
(2\delta+D)}e^{\ep (U_{p,q}-\delta)} = e^{-  D} e^{-3
\delta}e^{ U_{p}(q)}.$$

For every $D\geq 0$ and setting $q=\ga^{-1} p$ we have
\begin{align*}
1 & \geq  \ga^{-1}_{\star}\nu(O_p(q,D)) =\int_{O_p(q,D)}
d\ga^{-1}_\star\nu(y)= \int_{O_p(q,D)} C^{-1} e^{-\alpha d(p,
q)}e^{2\alpha( y\cdot
q)_p} d\nu( y)\\
& \geq C^{-1} e^{\alpha(2U_{p}(q)-d(p,q))} e^{-2\alpha D}
e^{-6\alpha \delta}\nu(O_p(q,D)).
\end{align*}
Using $d(p, \ga^{-1}p) \leq U_{p,\ga}+K$ we obtain the result. For part
b') we simply replace the last estimate by
\begin{align*}
\norm{\nu_q} & \geq  \nu_q(O_p(q,D)) =\int_{O_p(q,D)} d\nu_q(y)
=\int_{O_p(q,D)} C^{-1} e^{-\alpha d(p, q)}e^{2\alpha( y\cdot
q)_p} d\nu_p( y)  \\
& \geq C^{-1} e^{\alpha(2U_{p}(q)-d(p,q))} e^{-2\alpha D}
e^{-6\alpha \delta}\nu_p(O_p(q,D)).
\end{align*}
\end{proof}

\medskip
A consequence of the generalized Shadow lemma above is that it can
be plugged into Coornaert's proof of Corollary 6.6 in
\cite{Coornaert93} to obtain the following.

\begin{corollary}
If $X$ is Gromov product bounded with respect to $\Gamma$, then any
$\alpha$-quasiconformal measure is finite and nonzero only if
$\alpha\geq \delta(\group)$.
\end{corollary}

For convex cocompact actions on a CAT(-1) space the the Shadow
Lemma implies that the Patterson-Sullivan measures
$\set{\nu_p}_{p\in X}$ are Ahlfors $Q$-regular: there exists a
$C>0$ such that for all $x\in \op{CH}(\Lambda)$ and $r<1$,
$$C^{-1} r^Q\leq \nu_p(B(x,r)) \leq C r^Q.$$
In general Gromov hyperbolic spaces this need not be the case. However, we
do not need such a strong condition.  We only need a simple decay
condition on the measure expressed in terms of a singular integral. This
will be taken care by the following lemma.
\begin{lemma}\label{lem:integ_estim}
\noindent
\begin{enumerate}
\item If $\set{\nu_x}_{x \in X}$ is a lower $\alpha$-quasiconformal density with
constant $C$and $z_{p,q}^+\in \pa X$ is any point such that
$(z_{p,q}^+,q)_p\geq U_{p}(q)-\delta$, then
$$\int_{X-O_p(q, 0)} \frac{1}{d_p(z_{p,q}^{+},  y)^{2\alpha}}
d\nu_p( y) \leq C\norm{\nu_q} e^{\alpha d(p, q) + 6\alpha\delta}.$$
\item If $\nu$ is lower $\alpha$-quasiconformal with respect to $\group$ with
constant $C$ and $z_{p,\ga}^+\in \pa X$ is any point such that
$(z_{p,\ga}^+,\ga^{-1}p)_p\geq U_{p,\ga}-\delta$, then
$$\int_{X-O_p(\ga, 0)} \frac{1}{d_p(z_{p,\ga}^{+},  y)^{2\alpha}}
d\nu( y) \leq C\norm{\nu} e^{\alpha d(p,\ga^{-1}p) + 6\alpha\delta}.$$
\end{enumerate}
\end{lemma}

\begin{proof}
By Lemma \ref{lem:compare_to_distance} and the choice of $z_{p,q}^+$ we
have
  $$e^{( y \cdot q)} \ge \frac{e^{-3 \delta
 }}{d_p(z_{p,q}^+,  y)}.$$

Therefore we have,
\begin{align*}
 C^{-1} & e^{-\alpha d(p,q)} \int_{X-O_p(q,
0)}\frac{e^{-6\alpha \delta
 }}{d_p(z_{p,q}^+,  y)^{2\alpha}} d\nu_p( y) \\ &\leq
 \int_{X-O_p(q,0)} C^{-1}e^{-\alpha d(p,q)}e^{2\alpha( y \cdot
 q)_p}d\nu_p(y) \\ &\leq  \int_{X-O_p(q,0)}
 \frac{d \nu_q}{d\nu_p}( y) d\nu_p( y) \leq
 \nu_q(O_p(q,0)^c) \leq \nu_q(X). \hskip .5in
\end{align*}
For the second statement we restrict to the case when $q=\ga^{-1}p$ and
replace $\nu_p$ and $\nu_q$ with $\nu$ and $\ga^{-1}_\star\nu$
respectively.
\end{proof}

Recall that a measure is ergodic with respect to an action of
$\group$ if for any measurable $\group$ invariant set $A$, either
$\mu(A)=0$ or $\mu(A^\complement)=0$. The following result is
Proposition 3.3.1 of \cite{Yue96} for the case of discrete group
acting freely on a pinched negatively curved Hadamard manifold. We
present his proof to show that it works in the hyperbolic setting as
well.

\begin{proposition}\label{prop:coincide}
Any two $\alpha$-conformal densities coincide up to a constant
multiple if and only if $\group$ acts ergodically on $\Lambda$.
\end{proposition}

\begin{proof}
Consider two such conformal densities $\set{\mu_x}_{x\in H}$ and
$\set{\nu_x}_{x\in H}$. If $\group$ is ergodic with respect to the
measure class of $\mu$, consider the measure $\sigma = \frac12(\mu +
\nu)$ which is clearly also an $\alpha$-conformal density. Since
$\mu_x$ and $\nu_x$ are both absolutely continuous to $\sigma_x$, their
Radon- Nikodym derivatives $\frac{d\mu_x}{d\sigma_x} ,
\frac{d\nu_x}{d\sigma_x}$ exist and are $\group$ -invariant. By the
ergodicity of these derivatives are equal to positive constants
$\mu_x,\nu_x$-almost everywhere.

If $\group$ is not ergodic with respect to the measure class of
$\mu$, then there exists a Borel $\group$-invariant subset $A
\subset \La$ such that for all $x \in H, \mu_x(A) > 0$ and
$\mu_x(A^\complement)>0$. Defne $\sigma_x(E) = \mu_x(E\cap A)$; then
$\sigma_x$ is another $\alpha$-conformal density.
\end{proof}
\section{Metrics and measures on the boundary of a CAT(-1) space}\label{sec:metric_measure}

In this section we describe some improvements that can be made on
the the previous section if we restrict our attention to the case
when $X$ belongs to the family of CAT($-\kappa$) spaces for
$\kappa>0$. These are a special class of $\delta$-hyperbolic spaces
that are uniquely geodesic, proper and for which the Gromov product
extends continuously to the boundary.  We recall definitions and few
simple lemmas related to such spaces from \cite{BridsonHaefliger99}.
It is no less general to assume $\kappa=1$ since this can be
achieved by simply rescaling the metric. For the remainder of this
section, let $H$ be a CAT(-1) space with metric $d$ and assume
$\group$ is a (nonelementary) discrete isometry group acting on $H$
properly discontinuously.

\subsection{Busemann functions revisited}
Let $\pa H$ be the ideal boundary of $H$. Again denote by
$\Lambda$ the limit set of $\group$ on $\pa H$, and
$\Lambda_r$ the radial limit set of $\group$ on $\pa H$.

For $x,y \in H$ and $ z \in \pa H$, the Busemann function
$\rho_{x,z}$ takes on the simpler expression:
$$\rho_{x, z}(y)=\lim_{z_i \to z} (d_{H}(y,z_i)-d_{H}(x,z_i)).$$

It is easy to observe that
$$\rho_{x,  z}(y)-\rho_{x^{\pr},  z}(y)= \rho_{x,  z}(x^{\pr}),$$
for all $x,x^{\pr}, y \in H$ and $ z \in \pa H$. Moreover, $\rho$ is
continuous on $H\times H\times\pa H$ (\cite{Bourdon95}). One
important consequence of these facts is that we obtain an exact
formula (\cite{Bourdon95}) for the Radon-Nikodym derivative of the
Patterson-Sullivan measures:

$$\frac{d \nu_{y}}{d \nu_x}( z) = e^{-\delta(\group) \rho_{x, z}(y)}.$$

In particular, we have the following.
\begin{corollary}\label{cor:patt-sull_conformal}
The family  $\set{\nu_x}_{x\in H}$ form a continuous
$\alpha$-conformal density on $\pa H$.
\end{corollary}
\subsection{Metrics on the boundary $\pa H$}
Throughout this subsection we will  fix a point ${p} \in H$ and $\ep
>0$. We introduce two classes of metrics on the $\pa H$.  Even
though they are equivalent, i.e. bilipschitz, it will be convenient
to use both of them for later results.

{\underline{\it Busemann metric}:}
The Busemann
metric is defined to be
$$d_{p}^{\ep}( z, y) \triangleq e^{-\frac{1}{2}\epsilon (\rho_{x,
    z}(p) + \rho_{x, y}(p))}$$
where $x$ is any point on the geodesic connecting $ z$ and $ y$.  It
is not too difficult to check that this definition is independent of
the choice of $x$. Moreover, untangling the definition of the
Busemann function shows that $d_{p}^{\ep}=e^{-\ep( z\cdot
y)_{p}}=d_{p}^\ep$ and so we will continue to denote it by
$d_{p}^{\ep}$. However, in the case of a CAT$(-\kappa)$-space, this
is a genuine metric for all $0<\ep\leq \sqrt{\kappa}$ (see
\cite{Bourdon96}).

\medskip
{\underline{\it Shadow metric}
For  $ z,  y \in \pa H$ we define
$$\ell_{p}( z, y) \triangleq \sup\set{t\,|\,, d(\gamma_{p,  z}(t), \gamma_{p,  y}(t))\leq 1}$$
$$Sh_{p}^{\ep}( z,  y) \triangleq e^{-\ep \ell_{p}( z,  y)},$$
where $\ga_{p,  z}$ and $\ga_{p,  y}$ are the geodesics starting
from ${p}$ and pointing to $ z$ and $ y$. This metric is called {\it
Shadow metric}.

The following lemma, which states that these two metrics are
equivalent, can be found in \cite{Kaimanovich90} for the pinched
negatively curved manifold case. However, the proof works the same
for CAT(-1) spaces as well.
\begin{lemma}\label{lem:equiv_metrics}
There exist $\ep_0>0$ such that $d_{p}^{\ep}(.,.)$ and
$Sh_{p}^{\ep}(.,.)$ are distances for all $ 0<\ep\leq \ep_0$ and
${p} \in H$. Moreover these metrics are bilipschitz for the same
$\eps$: there is a $C>0$ such that for all $ z,  y \in \pa H$,
$$C^{-1} \leq \frac{d_{p}^{\ep}( z,  y)}{Sh_{p}^{\ep}( z,  y)} \leq C.$$
\end{lemma}

\begin{remark}
For the sake of comparison to the existing literature on negatively
curved spaces, we point out that there exist several other natural
metrics on $\pa H$ which lie in the same bilipschitz class as the
ones above. They include  various explicit metrics under the name
"Gromov's metric," and others with names like the "geodesic metric,"
the "horospherical metric," etc.... We omit their definitions since
they will be unecessary for the discussion of this paper, but for
many constructs we could have used them instead. However, the next
lemma shows that $d_{p}^{\ep}$ is preferable in the CAT(-1)
category.
\end{remark}

\begin{lemma} \label{lem:conf}
 Each isometry $\ga$ on $H$ induces a conformal map
on $\pa H$ under the metric $d_{p}^{\ep}$.
\end{lemma}

\begin{remark}
  This lemma justifies the previous notation of an $\alpha$-(quasi)conformal
  density. We are requiring that such a measure transform in a
  (quasi)conformal way under the family $\group$ of conformal maps
  with respect to this metric.
\end{remark}

\medskip
\begin{proof}
(see \cite{Yue96}.) The proof is almost the same as in the case of
negatively curved manifolds, but we will repeat it for completion
and to help the reader to get used to this notation.

Observe that for any $x, y \in H$ and
$ z \in \pa H$ and $\ga \in \group$ we have

$$\rho_{x, \ga z}(\ga y)- \rho_{x,  z}(y) = \rho_{\ga^{-1}x,
 z}(y)- \rho_{x,  z}(y) = \rho_{\ga^{-1}x,  z}(x).$$

So for ${p}, q \in X$ and $ z,  y \in \pa H$ we have
$$\begin{aligned}
\frac{d_{q}^{\epsilon}(z,y)}{d_{p}^{\epsilon}( z, y)}= &
\frac{e^{-\ep(y\cdot z)_{q}}}{e^{-\ep(y \cdot z)_{p}}}= e^{ \ep(
(y\cdot z)_{p}-(y \cdot z)_{q})}=e^{\ep ((y \cdot q)_{p}+(z \cdot
q)_{p}-d(p,q))},\end{aligned}$$ by continuity of the Gromov product.

Now since the topology on $\pa H$ is inherited from the Gromov
product, for each $a>0$ there exists a neighborhood $V_a$ of $ z$ in
$\pa H$ such that for all $ y \in V_a$ we have
$$ (1-a) \leq e^{\ep((y \cdot q)_{p}-(z \cdot q)_{p})} \leq (1+a).$$

So we proved that for every $a>0$ there exists a neighborhood
$V_a\ssu \pa H$ of $ z$, such that for all $ y \in V_a$ we
have
$$(1-a)e^{\epsilon(2(z \cdot q)_{p}- d(p, q))} \leq \frac{d_{q}^{\epsilon}(z,
y)}{d_{p}^{\epsilon}( z, y)} \leq (1+a)e^{\epsilon(2(z \cdot q)_{p}-
d(p, q))}.$$
 This proves that $$\lim_{ y
\to  z} \frac{d_{q}^{\epsilon}(z, y)}{d_{p}^{\epsilon}( z, y)}
=e^{\epsilon(2(z \cdot q)_{p}- d(p, q)) }= e^{-\ep \rho_{p ,z}(q)} .
\hskip 1in$$
\end{proof}

The propositions and proofs for the rest of this section are taken
from \cite{Yue96} which work just as well in the CAT(-1) setting as
for negatively curved Hadamard manifolds.

First we present the analogue of Proposition
\ref{prop:quasiconformal}. It occurs as Proposition 3.1.1 of
\cite{Yue96} in the case $H$ a pinched negatively curved manifold.
However, it easily follows by the same proof using Lemma
\ref{lem:conf} in the current setting of a CAT(-1) space.

\begin{proposition}\label{prop:haus_conformal} If $A \ssu \pa H$ is a $\group$-invariant
Borel set with $0<\mathcal{H}_{d}^{\alpha/\ep}(A)<\infty$ for any
quasimetric $d$ bilipschitz to $d_{p}^{\ep}$ and $\alpha>0$, then
$\set{\mathcal{H}_{d_{p}^{\ep}}^{\alpha/\ep}}_{p\in H}$ is a
continuous $\alpha$-conformal density of $\group$.
\end{proposition}

Next we collect other pertinent results which, except where
indicated below, do not depend on the discreteness of $\group$.

\begin{proposition}[Theorem A, Corollaries 3.5.3 and 3.5.6 of
\cite{Yue96}]\label{prop:conformal_CAT_props} If
$\set{\sigma_x}_{x\in H}$ is any finite $\alpha$-conformal density
for $\group$ then $\alpha\geq \dim_H(\La_r)$ and the following are
equivalent for any $p\in H$,
\begin{enumerate}
  \item $\sigma_p(\La_r)>0$
  \item $\La_r$ has full $\sigma_p$ measure
\end{enumerate}
If $\Gamma$ is discrete and properly discontinuous then
  these are equivalent to
\begin{enumerate}
\item[(3)] The Poincar{\`e} series
$\ds{\sum_{\ga\in\group}e^{-\alpha d(p,\ga^{-1}p)}}$ diverges.
\end{enumerate}
Moreover, if any of these conditions hold then
$\alpha=\delta(\group)=\dim_H(\La_r)=\dim_H(\La)$.

\end{proposition}

If $\Gamma$ is a discrete group acting convex cocompactly, then
$\La_r=\La$ so finite nonzero $\delta(\group)$-quasiconformal
measures exist. Hence, that the Hausdorff measures
$\mathcal{H}_p^{\delta(\group)}$ and the Patterson-Sullivan measures
$\nu_p$ (regardless of the choice of weak limit) coincide up to a
constant multiple for each $p\in X$. In particular, all of the
properties listed in Proposition \ref{prop:conformal_CAT_props}
hold.

\smallskip
\section{Regularity of measures and covers}\label{sec:regularity_covers}

We first present some necessary notation. Recall that on a set $X$,
a nonnegative function $d:X\times X\to \R$ is called a
\emph{quasimetric} (or \emph{quasidistance}) if $d$ is symmetric,
zero precisely along the diagonal, and satisfies the quasitriangle
inequality:
$$d(x,y)\leq C(d(x,z)+d(z,y))$$ for some $C\geq 1$ and all $x,y,z\in X$.

\subsection{Doubling and related properties}
 Let $X$ be a space equipped with a quasidistance function
$d(\cdot,\cdot)$ and probability measure $\nu$. $\norm{\cdot}$ denotes $L^1(X,
\nu)$-norm.

\begin{definition}
\begin{enumerate}
\item We say that a measure  $\nu$ has {\it $(p,\alpha)$-decay} if
 there is a constant $D_\nu$ such that:
\begin{gather}\label{eq:integ_cond}
\int_{X-B(x, r)}\frac{1}{d(y,x)^{p+\alpha}}d\nu(y)\leq
 \frac{D_{\nu}}{r^p},
\end{gather}
 for every $x \in X$ and $1\geq r>0$. We
 replace the right hand side by $D_{\nu} \left(1+
   \left|\log\left(r\right)\right|\right)$ if $p=0$ and we replace $X$
 by $B(x,1)$ if $p+\alpha<0$.

\item We say that the measure $\nu$ has {\it  upper
$Q$-regularity} if there exists a constant $K_{\nu}>0$ such that the
following is satisfied:
$$\nu(B(x, r)) \leq K_{\nu} r^Q,$$ for every $x \in X$ and $r>0$.

\item We say that the measure $\nu$ has the {\it strong doubling
property } if there exists a constant $T_{\nu}$ such that
$$\nu(B(x, 3r)) \leq T_{\nu}\nu(B(x, r)),$$ for all $x \in X$ and
$r >0$.
\end{enumerate}

\end{definition}

\begin{remark}
  Note that if $\nu$ has the strong doubling property then to verify
  that $\nu$ has $(p,\alpha)$ decay, it is enough to verify
  condition \eqref{eq:integ_cond} for $r\in \set{1,c,c^2,c^3,\dots}$ for any
  $0<c<1$. Any other radius can be interpolated.
\end{remark}

\begin{lemma}\label{lem:regular_has_decay}  Assume that $\nu$ is a finite measure with upper
  $Q$-regularity for some $Q\geq 0$. Then $\nu$ has $(p,\alpha)$-decay
  for all $p\geq 0$ and $\alpha\in\R$ satisfying $-p\leq\alpha\leq Q$.
\end{lemma}
\medskip
\begin{proof}
  First observe that if $p+\alpha\geq 0$, then $$\int_{X-B(x,
    r)}\frac{1}{d(y,x)^{p+\alpha}}d\nu(y)\leq
  \nu(X-B(x,1))+\int_{B(x,1)-B(x,
    r)}\frac{1}{d(y,x)^{p+\alpha}}d\nu(y).$$

Without loss of generality, we may assume $r=\frac{1}{N}$ for
$N>>1$, then
$$\begin{aligned}
\int_{B(x,1)-B(x, r)}\frac{1}{d(y,x)^{p+\alpha}}d\nu(y) &\leq
\sum_{n=1}^{N} \frac{1}{(nr)^{p+\alpha}}
(\nu(B(x,(n+1)r))-\nu(B(x,nr)))  \\
& \leq \sum_{n=1}^{N} \left(\frac{1}{(nr)^{p+\alpha}} -
\frac{1}{((n+1)r)^{p+\alpha}}\right) \nu(B(x,(n+1)r))  \\ & \leq
\sum_{n=1}^{N} \left(\frac{1}{(nr)^{p+\alpha}} -
\frac{1}{((n+1)r)^{p+\alpha}}\right) K_{\nu} ((n+1)r)^Q \\ & \leq
\sum_{n=1}^{N} \left(\frac{(n+1)^{Q}}{(n)^{p+\alpha}} -
\frac{(n+1)^{Q}}{((n+1))^{p+\alpha}}\right) \frac{K_{\nu}r^{Q-\alpha}
}{r^{p}}.\end{aligned}$$ We can estimate, under the
assumption that ${p+\alpha}\geq 0$,
$$ \frac{(n+1)^{Q}}{(n)^{p+\alpha}} -
\frac{(n+1)^{Q}}{((n+1))^{p+\alpha}}= \frac{(n+1)^{p+\alpha}-n^{p+\alpha}}{(n)^{p+\alpha}
(n+1)^{p+\alpha-Q}}\leq C_{p+\alpha} \frac{n^{p+\alpha-1}}{(n)^{p+\alpha} (n+1)^{p+\alpha-Q}}\leq
\frac{C_{p+\alpha}^\pr}{n^{p+\alpha-Q+1}},$$ for some constants $C_{p+\alpha}$,$C_{p+\alpha}^\pr$
and all $n\geq 1$. (We may take $C_{p+\alpha}^\pr=C_{p+\alpha}$ when ${p+\alpha}\geq Q$.)
Therefore, if ${p+\alpha}\neq Q$, the above sum becomes
$$\begin{aligned}
  \frac{K_{\nu}r^{Q-\alpha}}{r^{p}}\sum_{n=1}^{N}\frac{C_{p+\alpha}^\pr}{n^{p+\alpha-Q+1}}&\leq
  \frac{K_{\nu}C_p^{\pr}}{r^{p}}\frac{1}{N^{Q-\alpha}}\left(\frac{N^{Q-p-\alpha}-1}{Q-p-\alpha}\right)\\
  &\leq
  \frac{K_{\nu}C_{p+\alpha}^\pr}{r^{p}}\left(\frac{N^{-p}-N^{\alpha-Q}}{Q-p-\alpha}\right).\end{aligned}$$
The quantity $\left(\frac{N^{-p}-N^{\alpha-Q}}{Q-p-\alpha}\right)$ is
bounded for $p\geq 0$ whenever $Q-\alpha\geq 0$. When ${p+\alpha}=Q$ we
may bound the sum by $1+\log(1/r)$, obtaining the same conclusion. By
using logarithmic bounds, the same conclusion is obtained in the case
${p+\alpha}=Q$. (Note we obtain the desired estimate
$K_{\nu}C_{p+\alpha}^\pr(1+\log(1/r))$ in the case ${p+\alpha}=Q$ and $p=0$.)
\end{proof}

\medskip

\begin{remarks}
The above proof also shows that if $\nu$ is any measure with
$(p,\alpha)$-decay, then it also has $(p+t,\alpha)$-decay for all $t\geq
0$.

If $X$ has Hausdorff dimension $Q$ and $\nu$ is the corresponding
Hausdorff measure, then $\nu$ has the strong doubling and
$(p,\alpha)$-decay properties for all $p>0$ and $\alpha\leq Q$. However,
$\nu$ may not be upper $Q$-regular.
\end{remarks}

\subsection{Lipschitz constants}

We shall say that a map $f:X\to Y$ between quasimetric spaces $X$
and $Y$ is locally Lipschitz if for every $r>0$ and $x\in X$ we
have,
$$\sup_{\stackrel{y \in B(x,r)}{y\neq x}}
\frac{d_Y(f(x),f(y))}{d_X(x,y)} <\infty.$$  The various extant
definitions of this notion for the most part agree when $X$ is
proper. We now recall the definition of the Lipschitz constant on
a given scale.

\begin{definition}
For a locally Lipschitz map $f:X\to Y$,
we define the {\it Lipschitz constant at $x$ of scale $r$} to be
the quantity,
$$D_rf(x) = \sup_{\stackrel{y \in B(x,r)}{y\neq x}}
\frac{d_Y(f(x),f(y))}{d_X(x,y)}.$$
\end{definition}

\begin{remark}
It is clear that if $f$ is locally Lipschitz and $s \leq r$, then $D_sf(x)
\leq D_rf(x)$.
\end{remark}

In the case of locally Lipschitz functions to $\R$, we summarize any
arithmetic relations we will need in the following lemma. These will
be mainly used in the proof of Theorem \ref{thm:first_moment}. Each
case may be verified by a simple (and omitted) computation based on
the definition.

\begin{lemma}\label{lem:Lip_const} If $F, G$ are two locally Lipschitz
  functions on $X$, then $F+G$ and $FG$ are locally Lipschitz.
  Moreover,
\begin{enumerate}
\item $D_r(F+G)(x) \leq D_r F(x) + D_r G(x)$,
\item $D_r(cF)(x) = |c|D_rF(x)$, and
\item $D_r(FG)(x) \leq \left(\sup_{d(x,y)\leq r}\abs{F(y)}\right)D_rG(x) + \left(\sup_{d(x,y)\leq r}\abs{G(y)}\right) D_r F(x)$.
\item If $G(x)\neq 0$ for all $x\in X$, then $\frac1G$ is locally Lipschitz
and
$$D_r\left(\frac{1}{G}\right)(x) \leq \frac{D_r
  G(x)}{\abs{G(x)}\left(\inf_{d(x,y)\leq r}\abs{G(y)}\right)}.$$
\item If $H:Y\to X$ is locally Lipschitz, then
$$D_r(F\circ H)(x)\leq D_{r*D_r(H)(x)}(F)(H(x))*D_r(H)(x).$$
\end{enumerate}
\end{lemma}

Because we will be using the operator $D_r$ with respect to (quasi)metrics
which differ by taking powers, it will be convenient to quantitatively
state the following well-known relationship.

\begin{lemma}\label{lem:Lip_const_power}
  Let $F$ be a locally Lipschitz function on $X$ with respect to $d^{\eps^\pr}$
  for some quasimetric $d$ and $\eps^\pr>0$. Choose any $\eps$ with $0<
  \eps\leq \eps^\pr$.  If $D^{\eps}_r F$ and $D^{\eps^\pr}_r F$ represent
  the local Lipschitz constants of $F$ at scale $r$ with respect to
  quasimetrics $d^\eps$ and $d^{\eps^\pr}$ respectively, then
\begin{align*}
D^\eps_r(F)(x)&\leq
\(\frac{\eps^\pr}{\eps}\)^\frac{\eps}{\eps^\pr}\(\sup_{B^\eps_{r}(x)}F^{1-\frac{\eps}{\eps^\pr}}\)\,
  \(D^{\eps^\pr}_{r^\frac{\eps^\pr}{\eps}}(F)(x)\)^\frac{\eps}{\eps^\pr}.\\
\intertext{In particular, if $r\leq 1$ then}
D^\eps_r(F)(x) &\leq
  2 \(\sup_{B^\eps_{r}(x)}F^{1-\frac{\eps}{\eps^\pr}}\)\,
  \(D^{\eps^\pr}_{r}(F)(x)\)^\frac{\eps}{\eps^\pr}.
\end{align*}
 Here
  $B_r^\eps(x)$ is the $r$-ball centered at $x$ with respect to
  $d^\eps$.

\end{lemma}

\begin{proof}
  The result follows by taking the supremum over $y\in B_r^\eps(x)$ of the estimate,
\begin{align*}
  \frac{\abs{F(x)-F(y)}}{d^\eps(x,y)}& \leq
  \(\frac{\abs{F(x)^\frac{\eps^\pr}{\eps}-F(y)^\frac{\eps^\pr}{\eps}}}{d^{\eps^\pr}(x,y)}
  \)^\frac{\eps}{\eps^\pr}  \\
  &\leq \(\frac{\eps^\pr}{\eps}
  \frac{\abs{F(x)-F(y)}}{d^{\eps^\pr}(x,y)}
  \max\set{F(x)^{\frac{\eps^\pr}{\eps}-1},F(y)^{\frac{\eps^\pr}{\eps}-1}}\)^\frac{\eps}{\eps^\pr}.
\end{align*}
The second estimate comes from noting that $r^{\frac{\eps^\pr}{\eps}}\leq
r$ and that $x^{-x}\leq2$ for all $x>0$.
\end{proof}
\subsection{Besicovitch covers}

Recall that a covering $\set{U_{\alpha}}_{\alpha\in\mathcal A}$ of a
space $X$ has Lebesgue number $B\in\nn$ if for every point $x \in
X$,
$$0<\#\set{\alpha\,: \, x \in U_\alpha}\leq B.$$
\begin{definition}
A collection $\set{U_{\alpha}}_{\alpha\in\mathcal A}$ of subsets of
a space $X$ is called a {\em weak cover with Lebesgue number} $B$ if
$\nu$-almost every point of $X$ lies in at least one $U_\alpha$ and
for all $x\in X$ we have
$$\#\set{\alpha\,: \,x \in U_\alpha}\leq B.$$
\end{definition}
\begin{definition}
We say that a collection $\set{U_{\alpha}}_{\alpha\in\mathcal A}$
where $U_{\alpha} \subset X$ is a {\em
  Besicovitch cover,} if there is a constant $B$ such that for every
$\epsilon>0$, there exists a countable sub-collection
 $\set{U_{\alpha_i}}_{i=1}^{\infty} \ssu
\set{U_{\alpha}}_{\alpha \in\mathcal A}$ with
$\op{diam}{(U_{\alpha_i})}\leq \epsilon$ which forms a covering with
Lebesgue number $B$. We call a Besicovitch cover {\em profinite} if
for each $\eps>0$, the corresponding subcover can be chosen to be
finite. Similarly, we define a {\em weak Besicovitch cover} and a
{\em
  weak profinite Besicovitch cover} as above, except that the
subcollections  are only expected to be weak covers with Lebesgue number
$B$.
\end{definition}

\section{Spikes}\label{sec:spikes}
To keep the discussion as general as possible, in this section we
assume that our measure is not a single atom and has support $X$.

\begin{definition}\label{def:spike}
Assume $\nu$ is a measure with $(Q,\theta)$ decay for some $Q\geq 0$
and $\theta\in \R$. A $6$-tuple $(h(x),r, a, Q,\theta, C)$ where
$h(x)$ is positive function on $X$, $r>0$, $C> 1$ and $a \in X $, is
called a {\it spike} if
\begin{enumerate}
\item $h(x)\geq \norm{h}_{L^{\infty}}/C$ on $B(a, r)$, \\
\item for each $y \in B(a, r)^c$ we have
$$0<h(y)\leq h(a) r^Q\int_{B(a,r)}
\frac{C}{d(y,x)^{Q+\theta}}d\nu(x),\quad\text{and}$$\\
\item if $y,y'\in X$ satisfy $d(y, y^\pr)\leq r$, then
$h(y^\pr) \leq C h(y)$.\\
\end{enumerate}

If $h(x)$ is a continuous function we call $(h(x),r, a, Q, \theta, C)$ a
{\it continuous spike}. Also if $\norm{h}_{L^{\infty}}=1$ we will call
$(h(x),r, a, Q, \theta,C)$ a {\it unit spike}. Lastly we will often denote
the spike by the function $h(x)$ alone with the other constants implicit.

\end{definition}

\begin{definition}
If in addition a spike $(h(x), r, a, Q, \theta, C)$ has $h(x)$  locally
Lipschitz with
$$D_r h(x) \leq \frac{C h(x)}{r}$$ for all $x \in X$, and
$$\nu(B(a,r )) \geq \frac{r^Q}{C},$$
then we call $h$ a $Q$-spike.
\end{definition}
\begin{definition}
For a spike  $(h(x),r, a, Q, \theta, C)$ we call the number
$\frac{\nu(B(a,5r))}{\nu(B(a, r))}$ the {\em local doubling
constant}.

When working with families of spikes, we will denote the supremum of
the local doubling constant over all spikes with final entry less
than $C$ by $T_{\nu,C}$. Similarly, $T_{\nu}$ indicates the supremum
of all local doubling constants. If for a given family of spikes we
have $T_{\nu,C}<\infty$ for all $C>1$, then we say that members of
this family are {\em $\nu$-spikes}.
\end{definition}

\begin{remark}
Note that we implicitly assume that $a$ must lie in the support of
$\nu$. For small families of spikes, it is much weaker to assume
that the constant $T_\nu$ is bounded than to assume that $\nu$ has
the strong doubling property.
\end{remark}

It is clear that a positive multiple of a spike is a spike. Here are few
simple Lemmas about spikes and $Q$-spikes. The first observation is
immediate.

\begin{lemma}\label{lem:spike_pinch_est}  Assume $(h, r, b, Q, \theta, C)$ is
  a spike.  If $\frac{1}{M}f \leq A g \leq Mf$ for some constants
  $A>0$ and  $M\geq 1$, then
  $(g,r,b,Q,\theta, M^2 C)$ is a spike.
\end{lemma}

\begin{lemma}
Assume $(h, r, b, Q, \theta, C)$ is a spike.  Then
$$\frac{\nu(B(b, r))}{C} \leq
\frac{\norm{h}_{L^1}}{\norm{h}_{L^{\infty}}}.$$
\end{lemma}

\medskip
\begin{proof}
Using the property 1) in the definition of a spike, we obtain
$$\norm{h}_{L^1} \geq \int_{B(b,r)}hd\nu\geq \nu(B(b,r))\frac{\norm{h}_{L^{\infty}}}{C}.$$

\end{proof}

\medskip

\begin{lemma}\label{lem:spike_power_est}  Assume that $(h, r, a, Q,\theta, C)$ is a spike.  Then
  for all $t \geq 1$, $(h^t, r, a, Q,\theta, C^t)$ is a spike.
\end{lemma}

\medskip
\begin{proof}
  Observe that conditions 1),3) hold trivially. Condition 2) holds
  because $C^t\geq C\geq 1$, so taking the power in $t$ preserves the
  inequality.
\end{proof}.

\medskip

\begin{corollary}\label{cor:spike_power_est} Assume that $(h,r,a,Q, \theta,C)$
  is a $Q$-spike, then for
  all $t \geq 0$, $(h^t,r,a,Q, \theta,\max\set{t C^t,C})$ is also a $Q$-spike.
\end{corollary}

\medskip
\begin{proof}
 We would like to use the following estimate. For all $a,
b \in \rr$ we have
$$|a^t-b^t| \leq \max\set{t,1} |a-b| \max{(|a|^{t-1}, |b|^{t-1})}.$$

Now since $h(x)$ is a $Q$-spike, $h(x) \geq 0$.  Also it is easy
to see that for all $x,y \in X$ such that $d(x,y)\leq r$ we have
$h(y)^{t-1} \leq c^{t-1} h^{t-1}(x)$.

Thus we have $$\begin{aligned} D_rh^t(x) &= \sup_{d(x,y)\leq
r}\frac{|h^t(x)-h^t(y)|}{d(x,y)} \\ &\leq \sup_{d(x,y)\leq r}
\max\set{t,1} \frac{|h(x)-h(y)|}{d(x,y)} \max{(h^{t-1}(x),
h^{t-1}(y))}\leq \\& \leq \max\set{t,z} C\frac{h(x)}{r}
C^{t-1}h^{t-1}(x) \leq C^t \max\set{t,1}\frac{h^t(x)}{r}.
\end{aligned}$$
Now if $t<1$ then $C>C^t$ so in all cases $\max\set{C,tC^t}\geq
C^t\max\set{t,1}$. Moreover, this constant is always at least as large as
$C$ so that the measure condition in the definition of the $Q$-spike $h$
persists.
\end{proof}

\begin{lemma}\label{lem:func_still_Q}
If $(h,r,a,Q, \theta,C)$ is a $Q$-spike and $f$ is any Lipschitz function
with $\frac{1}{K}\leq f\leq K$ and $D_r f\leq \frac{K}{r}$, then $\(f
h,r,a,Q,\theta,2K^2 C\)$ is also a $Q$-spike.
\end{lemma}

\begin{proof}
By Lemma \ref{lem:spike_pinch_est} it is a spike. Using the properties of
$f$ and the third property of the spike $h$, we have
\begin{align*}
D_r\(f h\)(x)&\leq \(\sup f\)D_r h+ \(\sup h\) D_r f \\
&\leq K \frac{C h(x)}{r}+C h(x) \frac{K}{r}\\
&\leq 2 K^2 C\frac{h(x) f(x)}{r}.
\end{align*}
Since $h$ was a $Q$-spike this completes the proof.
\end{proof}

\section{Basis}\label{sec:basis}

The purpose of this section is to prove the main ingredient for the
coarse version of the stationarity result. This theorem provides a
general criterion for the closure of the positive cone on a family
of continuous positive functions $\set{f_i}$ on a metric measure
space to be as large as possible in $L^1$. Here the closure is with
respect to either $L^1$, uniform, or pointwise convergence. To the
best of our knowledge, this is the first such general condition. We
also indicate some examples showing that the hypotheses are in some
aspects nearly sharp.

We call $(X,d,\nu)$ a {\em probability metric space} or {\em pm space}
if $\nu$ is a Radon probability measure (with respect to the topology
induced by $d$) of full support on $X$. We will also assume $X$ is
separable, though not necessarily complete. (N.B. the Radon condition
on $\nu$ could be replaced by the condition that $\nu$ is Borel and
for all measurable sets $S$, $\nu(S)=\sup \set{\nu(K): K\subset S \text{
  compact }}$.)

\begin{definition}
  A function $F$ defined on a probability metric space $(X,d,\nu)$ is
  {\em (almost) lower approximable} if $F$ is (resp. almost everywhere)
  equal to the pointwise limit of a nondecreasing sequence of
  continuous functions.

  Recall $F$ is {\em (almost) lower semicontinuous},  if $\liminf_{z\to
    x}F(z)\geq F(x)$ for every (resp. almost every) $x\in X$.

  Finally, $F$ is {\em (almost) uniformly positive} if there is a $c$ such
  that $F(x)\geq c>0$ for every (resp. almost every) $x\in X$. The
  greatest (essential) lower bound for such an $f$ may be expressed as
  $\norm{1/{f}}_{-\infty}$.
\end{definition}

Since $X$ is metrizable, it is perfectly normal. Hence the (almost) lower
approximable functions are exactly the (almost) lower semicontinuous
functions (\cite{Tong52}).

\begin{theorem} \label{thm:basis}
Let $(X, d, \nu)$ be a probability metric space such that $\nu$ has
$(Q,\theta)$-decay. Assume $\set{(f_\alpha, r_\alpha, b_\alpha, Q,
\theta, C_\alpha) }_{\alpha\in\mathcal{A}}$ is a family of
continuous unit $\nu$-spikes.
For any $C>1$ let $S_C = \set{\alpha\in\mathcal{A}\,:\, C_\alpha\leq
C}$ and set $$B_C(r) = \bigcup_{\substack{\alpha \in S_C,\\ r_\alpha
\leq
    r}} B(b_\alpha, r_\alpha) \quad \text{ and }\quad B_C =
    \bigcap_{r>0} B_C(r).$$
If $\lim_{C \to \infty} \nu(B_C)=1$, then there exists a countable
subset of indices $\set{\alpha_i}_{i=1}^\infty\subset \mathcal{A}$
such that for any uniformly positive almost lower semicontinuous
function $F\in L^1(X,\nu)$ there exists a sequence $\set{
\la_{\alpha_i} }_{i=1}^\infty$ of nonnegative numbers such that
$F=\sum_{i=1}^\infty \la_{\alpha_i} f_{\alpha_i}$ for $\nu$-almost
every $x\in X$ with convergence in $L^1$. Moreover, if $F$ is a
lower semicontinuous (resp. continuous) function, then the
convergence is pointwise (resp. uniform on compacta) on
$\bigcup_{C>0}B_{C}$.
\end{theorem}

\begin{remarks}\label{rem:after_thm}
  Every function $F$ which can be expressed as a positive sum of the
  continuous functions $f_\alpha$ must be lower semicontinuous and
  positive, since ignoring the tail of the series yields an
  increasing sequence of positive continuous functions whose limit
  is therefore lower semicontinuous. In this sense the conclusion of
  the theorem is sharp. In fact, the theorem asserts that
  $\sum_{i=1}^\infty \la_{\alpha_i} f_{\alpha_i}$ will be the lower
  semicontinuous hull of $F$ for a general uniformly positive $F\in
  L^1(X)$. When $X$ is noncompact and the $f_\alpha$ are not
  uniformly positive, it is conceivable that some nonuniformly
  positive $F$ could be approximated as well. However, to address
  this case we would have had to make complicated compatibility
  assumptions on how and where $F$ and the $f_\alpha$ decay to $0$
  which we felt were not worth the extra effort.

  Examples of almost lower semicontinuous functions include any $F$
  whose set of points of discontinuity (or undefined points)
  $D\subset X$ has measure $0$. Simply note that the lower
  semicontinuous hull $\til{F}$ of $F$ defined by
  $$\til{F}(x):=\begin{cases}F(x) & F(x) \text{ exists and }F(x)\leq
    \ds{\liminf_{y\to x}}F(y)\\ \ds{\liminf_{y\to x}F(y)} &
    \ds{\liminf_{y\to x}}F(y)<\infty \text{ and }x\not\in
    \op{Dom}(F) \text{ or } F(x)>\ds{\liminf_{y\to x}}F(y) \\ 1 &
    \text{otherwise} \end{cases}$$
  agrees with $F$ almost everywhere.

  One cannot remove the assumption that $\nu(\cup_{C>0} B_C)=1$, even
  if $\cup_{C>0} B_C$ is assumed to be dense, since an $L^1$ function
  with mass outside $\cup_{C>0} B_C$ could never be approximated in $L^1$.
  (Also, see the examples below.)

  Lastly, the difficulty in proving the theorem reflects a certain
  balancing act captured by the spike conditions. While it may be
  possible to weaken these, we briefly mention why these conditions
  are qualitatively necessary. Since there can be no cancellation in
  the sum, it is evident that the positive basis functions must
  contain subsequences which, when suitably normalized, converge to
  Dirac distributions based at almost every point. However the shape
  of the $f_\alpha$ are further constrained. If they all decay too
  quickly then any countable subset cannot generally approximate on a
  full measure set. If the $f_\alpha$ decay too slowly then their
  tails stack up too quickly far away from their maximum. This too
  prevents generic approximation. In fact, it is somewhat surprising
  that the spike conditions happen to be satisfied for the
  Radon-Nikodym derivatives of the geometric measures in our main
  application (see Section \ref{sec:derivs_are_spikes}).

\end{remarks}

\begin{example}
  Uniformly bounded upper semicontinuous functions need not be almost
  lower semicontinuous. For instance, let $X=[0,1]$ and define $F=F_C$
  to take value $2$ on a Cantor set $C$ of positive Lebesgue measure,
  and have value $1$ elsewhere.  Any lower continuous approximation
  must lie below $1$, otherwise all closer approximations are greater
  than $F$ by a fixed positive amount on an open set, a contradiction.
  Either way, this prevents the limit from converging in $L^1$ to
  $F_C$. Nevertheless, $F_C$ is upper semicontinuous since it is the
  pointwise limit of a nonincreasing sequence of continuous
  approximations to the step functions which take the value $2$ on the
  set of intervals representing the $n$th stage in the construction of
  the Cantor set and take the value $1$ elsewhere. The same argument
  shows that the conclusion of the theorem holds for many functions
  with a positive measure set of discontinuities, for instance,
  $F=3-F_C$.

\end{example}

\begin{example}
  In this example we present the simplest case to which we will apply
  the theorem. It also motivates the definition of a spike and further
  indicates why the hypothesis on the balls $B(b_i,r_i)$ is necessary.
  Consider $\mathbb{H}^2$ in the disk model. Let $\group$ be a
  discrete group of $\op{Isom}(\mathbb{H}^2)$ and for each $\gamma\in
  \group$ let $(d,\theta)=(d_\gamma,\theta_\gamma)$ represent the
  orbit point $\gamma\cdot 0$ where $d$ is its hyperbolic distance
  from $0$ and $\theta$ is its angle from the real axis. Then if
  $\phi$ is the angle coordinate on $\partial \mathbb{H}^2=S^1$, we
  may write the Poisson kernel based at $\gamma\cdot 0$ normalized to
  have maximum $1$ as
  $$f_{\gamma}(\phi)=\frac{1}{e^d\,\left( \cosh (d) - \cos (\theta -
      \phi )\,\sinh (d) \right) }.$$
  If $r_\gamma$ represents the
  radius on $S^1$ of the points where $f_\gamma>1/C$, then
  $r_\gamma=\arccos \left(C - (C-1)\coth (d)\right)$.  We
  will show in Section \ref{sec:derivs_are_spikes} that if $\nu$ is
  the Lebesgue probability measure on $S^1$, then the tuples
  $(f_\gamma,\theta_\gamma,r_\gamma,1,1,C)$  for all $\gamma\in \group$
  and $C=2,3,\dots$ form a family of continuous unit $\nu$-spikes. We
  will see that for the family of balls $B(\theta_\gamma,r_\gamma)$
  and sets $B_C$ as in the theorem, the radial limit set of $\group$
  coincides with $\cup_{C>1} B_C$.  However, it is easy to construct
  examples of $\group$ such that their radial limit set has measure
  zero and hence $\lim_{C\to\infty}\nu(B_C)=0$. For such a group
  $\group$, any function $F>0$ approximated by any positive sums of
  $f_\gamma$ must have more than $1/2$ of its $L^1$ norm concentrated
  on $\cup_{C>1} B_C$ since each term in the sum does. Of course one
  can provide a work around for this obstruction by simply restricting
  the measure to the limit set. This can be made to work so long as
  the radial limit set has full measure in the limit set.

  This case of the theorem for uniform approximations of continuous
  functions on the circle by Poisson kernels was first proved by
  Hayman and Lyons (\cite{HaymanLyons90}, see also
  \cite{BonsallWalsh89}) using the theory of harmonic functions.
  That result was later extended to Euclidean domains in \cite{Gardiner96}.
\end{example}

\begin{remark}
  An analysis of the proof of the next proposition shows that while
  the condition of the $\nu$-spike is probably not absolutely
  necessary for a positive basis, it is a very natural condition
  which, up to small possible improvements, is necessary for the
  intuitive approach we take.
\end{remark}

The proof of Theorem
\ref{thm:basis} will require the following proposition.

\medskip

\begin{proposition} \label{prop:subfunction}
  Let $(X, d, \nu)$ be a probability metric space such that $\nu$ has
  $(Q,\theta)$-decay. Assume $\set{(f_i(x), r_i, b_i, Q, \theta,C_i) }$ is a
  countable family of unit $\nu$-spikes on $X$ with bounded
  doubling constants and such that $C_i\leq C$. Let $Y\subset X$ be a
  set weakly covered by $\set{B(b_i, r_i)}$ with finite Lebesgue number.

For any positive function $F$ bounded away from $0$ on $X$ and bounded
from above on $Y$, set
$$t = \left(\frac{\sup_{z\in Y}(F(z))}{\inf_{z\in
      X}(F(z))}\right)^{\frac{1}{Q}}+1\geq 2.$$
Suppose for a given
$s>1$ there is a $\delta >0$ such that
$$s\geq \sup\set{\left. \frac{F(y)}{F(x)} \right| x\in X \text{ and }y\in Y \text{
    with } d(x,y)\leq \delta}$$
and $r_{i} \leq \delta/t$ for
all $i\in\nn$. Then,

\begin{enumerate}
\item[a.] There exists a constant $0<L_{\nu}<1$,  which
  do not depend on $F$, and a function $h=\sum_{i=1}^{\infty} \la_{i}
  f_{i}$ with each $\la_i \geq 0$ and only a finite number of the
  $\la_i$ not equal to $0$, such that for every $x\in X$
$$h(x)\leq  F(x)$$ and
$$\frac{L_{\nu}}{C^2\,s^2}\norm{F(x)}_{L_1(Y)}\leq
\norm{h(x)}_{L_1(Y)}.$$

\item[b.] If in addition $Y$ is (weakly) covered by a finite number of
  the $B(b_i,r_i)$, then
  $$\frac{L_{\nu}}{C^2\,s^2}F(x) \leq h(x),$$
  for ($\nu$-almost) every $x \in Y$.
\end{enumerate}
\end{proposition}

\begin{proof}[Proof of Proposition]

  In what follows, we let $T_\nu$ be the bound on the doubling
  constant. We let $B$ be the Lebesgue number of the cover, and we let
  $D_\nu$ be decay constant of $\nu$. Without loss of generality we
  also assume that $r_{i}$ is nonincreasing. Now we inductively build
  functions $g_n(x)$ as follows.

$g_0(x)=0$.
$$g_n(x)=\begin{cases}
g_{n-1}(x) & \text{ if }\quad g_{n-1}(b_{n})\geq F(b_{n})\\
g_{n-1}(x) + f_{n}(x)(F(b_{n})- g_{n-1}(b_{n})) & \text{ if }
\quad g_{n-1}(b_{n})< F(b_{n})\end{cases}$$

 This construction yields
a sequence $\la_1, \la_2, \dots$ with
$$g_n(x) =\sum_{i=1}^{n}
\la_i f_{i}(x),$$ and $0\leq \la_i \leq F(b_{n})$. Hence, $g_n(x)
\in V_{+}$.
\medskip

\begin{lemma}\label{lem:gypgeqgy}
If $y^\pr \notin \cup_{i=1}^n B(b_{i}, r_{i})$  and $d(y, y^{\pr})
\leq r_{n}$, then $g_n(y^{\pr})\leq Cg_n(y)$.
\end{lemma}

\begin{proof} Since $d(y, y^{\pr}) \leq r_{n}\leq r_{i}$ for all
$i\leq n$, by Remark above $f_{i}(y^{\pr}) \leq Cf_{i}(y)$ for all
$i=1, \dots, n$. Since $g_n(x)=\sum_{i=1}^n \la_i f_{i}(x)$ with
nonnegative coefficients, we obtain that $g_n(y^{\pr})\leq C
g_n(y)$.
\end{proof}

\begin{lemma}\label{lem:dleqrn} $g_n(y^{\pr})\leq Cg_n(y)+B\,s\,
  F(y^{\pr})$, for all $y, y^\pr \in X$ and $d(y, y^{\pr}) \leq
  r_{n}$.
\end{lemma}

\begin{proof} By the previous lemma we may assume that $y^{\pr} \in
  B(b_{k}, r_{k})$, for some $k\leq n$. By construction
  $g_n(x)=\sum_{i=1}^n \la_i f_{i}(x)$, with $0\leq \la_i\leq
  F(b_{i})$. Thus for each $i \not= k$, we have $f_{i}(y^{\pr}) \leq
  Cf_{i}(y)$.

  For $i=k$, we have that $\la_k \leq F(b_{k})$.  However, since
  $d(y^{\pr}, b_{k}) \leq r_{k}\leq \delta$, we have that
  $F(b_{k})\leq sF(y^{\pr})$. Since $B$ is the Lebesgue number of the
  weak cover, $\#\set{k\,|\,y^{\pr} \in B(b_{k}, r_{k})} \leq B$. Thus
  we conclude that $g_n(y^{\pr})\leq Cg_n(y)+B\,s\,F(y^{\pr})$.
\end{proof}

\begin{lemma}\label{lem:gnleqF} Let
  $L^{\pr}=s\,C+s\,B+CD_{\nu}B+s\,BT_{\nu}+s\,BCD_{\nu}2^{Q}$. Then
  $$g_N(y) \leq L^{\pr}F(y)$$
  for all $N\in\mathbb{N}$ and $y \in X$.
\end{lemma}

\begin{proof} Assume that $n$ is the smallest integer such that
  $g_{n}(y) \geq (s\,C+B\,s)F(y)$ for any $y \in X$. If there is no
  such $n$ then set $L'=(s\,C+B\,s)$. Otherwise, we have $g_{n-1}(y) <
  (s\,C+B\,s)F(y)$ for all $y\in X$.

  By Lemma \ref{lem:dleqrn}, $g_k(y)\leq Cg_k(y^{\pr}) +B\,s\,F(y)$
  for all $k\in\mathbb{N}$ and $y^\pr\in X$ such that $d(y,y^\pr) \leq
  r_{k}\leq \delta$.  Therefore, we obtain that $g_n(y^{\pr}) \geq
  s\,F(y)\geq F(y^\pr)$.
Thus if $d(y,b_{{k}}) \leq r_{n}$ and $k
>n$, by construction we have $\la_{k}=0$.

Now since $\lambda_i f_i(b_i)\leq F(b_i)$ and by the property of
spikes we have,
\begin{align*}
  \sum_{\substack{i\geq n,\\ \,\delta \geq d(y,\, b_{i}), \\
      d(y, B(b_{i}, \,r_{i}))\geq \frac{r_{n}}{2}}} &\la_i f_{i}(y)
  \leq  \sum_{\substack{i\geq n,\\\delta \geq d(y,\, b_{i}), \\ d(y,\,
      B(b_{i}, r_{i}))\geq \frac{r_{n}}{2}}} F(b_{i})r_i^Q
  \int_{B(b_{i}, r_{i})}\frac{C}{d(y,x)^{Q+\theta}} d\nu(x),\\
  \intertext{and since $r_i\leq r_n$ and $F(b_i)\leq s F(y)$ we have,}
  &\leq s\,F(y) r_n^Q
  B\int_{X-B(y,\frac{r_{n}}{2})}\frac{C}{d(y,x)^{Q+\theta}} d\nu(x),\\
\intertext{and by $(Q,\theta)$-decay this becomes,}
  &\leq s\,F(y)BCD_{\nu} 2^{Q}
\end{align*}

Now we consider $d(y,B(b_{i}, r_{i})) \leq \frac{r_{n}}{2}$ for
$i\geq n$. Therefore $B(b_{i}, r_{i}) \ssu B(y, 3/2 r_{n})$.

In case $r_{i}\leq r_{n}/2$, then $d(y, b_{i})\leq r_{n}$. As we
observe above,  in this case $\la_i=0$

So if $\la_i\not=0$ and $d(y,B(b_{i}, r_{i})) \leq
\frac{r_{n}}{2}$ for $i \geq n$, we obtain that $r_{i}\geq
r_{n}/2$. Since $d(y, r_{i}) \leq r_{n}$ we have that $B(y, 3/2
r_{n}) \ssu B(b_{i}, 5 r_{i})$.

\begin{align*}
 & \sum_{\substack{r_{i} \geq \frac{r_{n}}{2}, \\
d(y, B(b_{i}, \,r_{i}))\leq \frac{r_{n}}{2}}} \nu(B(b_{i},
r_{i})) \leq B \,\nu\left(B\left(y, \frac{3r_{n}}{2}\right)\right) \leq\\
& \leq \frac{B}{\# \set{i:r_{i} \geq \frac{r_{n}}{2}, d(y, B(b_{i},
\, r_{i})) \leq \frac{r_{n}}{2}}} \sum_{\substack{r_{i} \geq
\frac{r_{n}}{2},\\ d(y, B(b_{i}, \,r_{i}))\leq
\frac{r_{n}}{2}}} \nu(B(b_{i}, 5 r_{i})) \leq \\
&\leq \frac{B T_{\nu}}{\# \set{i:r_{i} \geq \frac{r_{n}}{2},d(y,
B(b_{i}, \, r_{i}))
\leq \frac{r_{n}}{2}}}  \sum_{\substack{r_{i} \geq \frac{r_{n}}{2}, \\
d(y, B(b_{i}, \, r_{i}))\leq \frac{r_{n}}{2}}} \nu(B(b_{i},
r_{i})).
\end{align*}

 Therefore
$$\# \set{i:r_{i} \geq \frac{r_{n}}{2},d(y,
B(b_{i}, \, r_{i})) \leq \frac{r_{n}}{2}} \leq T_{\nu}\, B.$$

Since $\la_i \leq F(b_{i})$ and $f_{i}(y)\leq 1$, we obtain
$$\sum_{\substack{i\geq n,\\ \,\delta \geq d(y,\, b_{i}), \\
d(y, B(b_{i}, \,r_{i}))\leq
\frac{r_{n}}{2}}} \la_i f_{i}(y) \leq \sum_{\substack{i\geq n,\\ \,\delta \geq d(y,\, b_{i}), \\
d(y, B(b_{i}, \,r_{i}))\leq \frac{r_{n}}{2}}} s\,F(y) \leq
s\,F(y)T_{\nu}B$$

Also observe that
\begin{align*}
&\sum_{\substack{i\geq n,\\ \,\delta \leq d(y,\, b_{i})}} \la_i
f_{i}(y) \leq \sum_{\substack{i\geq n,\\ \,\delta \leq d(y,\,
b_{i})}} F(b_{i}) r_{i}^Q \int_{B(b_{i}, r_{i})}
\frac{1}{d(y,x)^{Q+\theta}} d\nu(x) \leq \\
\leq & \sup_{z}(F(z)) r_{n}^Q B \int_{X-B(y,\delta-r_{n})}
\frac{C}{d(y,x)^{Q+\theta}} d\nu(x) \leq
\\ \leq & \frac{B C D_{\nu}\sup_{z}(F(z))r_{n}^Q }{(\delta-r_{n})^Q} \leq
\frac{B C D_{\nu}\sup_{z}(F(z))(\delta/t)^Q}{(\delta-\delta/t)^Q}
\leq
\\& \leq CD_{\nu} B \sup_{z}(F(z)) \frac{1}{(t-1)^{Q}}
\end{align*}

Now recall that $$t =
\left(\frac{\sup_z(F(z))}{\inf_z(F(z))}\right)^{\frac{1}{Q}}+1.$$
So we obtain that

$$\sum_{\substack{i\geq n,\\ \,\delta \leq d(y,\,
b_{i})}} \la_i f_{i}(y) \leq CD_{\nu}B\,F(y).$$

So for these choices we obtain
\begin{align*}
  g_N(y)&=\sum_{i=1}^{N}\la_i f_{i}(y)\\
  &= g_{n-1}(y)+ \sum_{\substack{N\geq i\geq n,\\ \,\delta \leq d(y,\,
      b_{i})}} \la_i f_{i}(y)+\sum_{\substack{N\geq i\geq n,\\
      \,\delta \geq d(y,\, b_{i}), \\
      d(y, B(b_{i}, \,r_{i}))\leq
      \frac{r_{n}}{2}}} \la_i f_{i}(y)+\sum_{\substack{N\geq i \geq n,\\
      \,\delta \geq d(y,\, b_{i}), \\
      d(y, B(b_{i}, \,r_{i}))\geq
      \frac{r_{n}}{2}}} \la_i f_{i}(y)  \\
  &\leq  (s\,C+s\,B)F(y)+
  CD_{\nu}BF(y)+ s\,F(y)B T_{\nu}+F(y)s\,BCD_{\nu}2^{Q} \\
  &= F(y) (s\,C+s\,B + CD_{\nu}B+s\,BT_{\nu}+s\,BCD_{\nu}2^{Q})
\end{align*}

Let $L^{\pr}=s\,C+s\,B+
CD_{\nu}B+s\,BT_{\nu}+s\,BCD_{\nu}2^{Q}>3$.
\end{proof}

Set $$3 L_{\nu} = \frac{1}{1+B+ D_{\nu}B+BT_{\nu}+
BD_{\nu}2^{Q}}.$$ Since $C> 1$ and $s \geq 1$ by Lemma
\ref{lem:gnleqF} we have
$$g_N(y) < \frac{C\,s}{3 L_{\nu}} F(y),$$ for all $y \in X$.
\medskip

Since $g_N(y)$ is bounded nondecreasing sequence, $\lim_{N\to
\infty}g_N(y) =g(y)$ is well-defined.

\begin{lemma}\label{lem:ggreaterF}  for every $y \in
\cup_{n=1}^{\infty}B(b_{n},r_{n})$,
$$g(y)\geq \frac{F(y)}{Cs}.$$
\end{lemma}

\begin{proof} Assume that $n$ is the smallest integer such that $y
\in B(b_{n}, r_{n})$. By the spike properties, we have
$$f_{n}(y)\geq
\frac{1}{C}=\frac{f_{n}(b_{n})}{C}.$$ By Lemma \ref{lem:gypgeqgy},
$$g_{n-1}(y) \geq \frac{1}{C}g_{n-1}(b_{n}).$$

So we obtain that $$g_n(y)\geq \frac{g_{n}(b_{n})}{C}.$$ Again by
construction $g_{n}(b_{n}) \geq F(b_{n})$. Also as $d(y,
b_{n})\leq r_{n}\leq \delta$, we get $s F(y)\geq F(b_{n})$.
Therefore,
$$g_n(y)\geq \frac{g_{n}(b_{n})}{C} \geq
\frac{F(b_{n})}{C} \geq \frac{F(y)}{sC}.$$

So we conclude that $$g(y)\geq g_n(y)\geq \frac{F(y)}{sC}.$$
\end{proof}

Now define
$$h_n(y)=\frac{3 L_{\nu} g_n(y)}{C\,s}\in V_{+}$$ and $\lim_{n \to
\infty}h_n(y) =h(y)$. By Lemma \ref{lem:gnleqF}, $h(y)\leq \beta F(y)$
where $\beta=\frac{3 L_{\nu}}{C\,s}L^{\pr}<1$. Note that $\beta$ is
bounded below by $3 L_\nu (1+BD_\nu 2^Q)$ independent of $C$ and $s$.
By Lemma \ref{lem:ggreaterF},
$$h(y)\geq \frac{3 L_{\nu}}{C^2\,s^2}F(y)$$ for $\nu$-almost every point.

So $$\|h(x)\| \geq \frac{3L_{\nu} \|F(x)\|_{L^1(Y)}}{C^2\,s^2}
> \frac{ L_{\nu}}{C^2\,s^2} \|F(x)\|_{L^1(Y)}.$$

Since $\lim_{n \to \infty}\|h_n(x)\|= \|h(x)\|$ there exists $n$
with required property.  This finishes part a).

\medskip For part b), it is easy to observe that if the cover is
finite, there exists $N$ such that for all $n>N$, $h_n(y)=h(y)$. If
$F$ is not locally constant then $t>1$ and therefore the conclusion of
Lemma \ref{lem:ggreaterF} is true on a uniform neighborhood of
$\cup_{n=1}^{\infty}B(b_{n},r_{n})$ which includes all of $Y$ since
a full measure subset was assumed to be dense.
This finishes the proof of Proposition \ref{prop:subfunction}.
\end{proof}

\medskip
\begin{proof}[Proof of Theorem \ref{thm:basis}]
  We begin by proving the first part of the theorem for the special
  case of uniformly positive continuous functions in $L^1$. We will
  find a universal countable subfamily of spikes indexed by
  $\set{\alpha_i}_{i=1}^\infty \subset \mathcal{A}$ so that for every
  such function $F$ we will inductively build a sequence
  $\set{h_i(x)}_{i=0}^{\infty} \ssu V_{+},$ where
  $V_{+}=\set{\sum_{i=1}^{n} a_{i} f_{\alpha_i}(x)\,:\, a_i\geq
  0,\,\alpha_i\in\mathcal{A},\, n \in \nn}$, such that for some
  $0<\gamma_n<1$,
\begin{align}\label{eq:*}\tag{*}
\begin{split}
&0< \sum_{i=0}^{n} h_i(x) \leq \gamma_n F(x)\quad
\text{ for all $x\in X$ and}\\
&\left\|F - \sum_{i=0}^{n}h_i\right\|_1\to 0 \hbox{ as } n \to \infty.
\end{split}
\end{align}
Once this is achieved we can simply note that the sums
$\displaystyle{\sum_{i=1}^{n} h_i}$ live in $V_{+}$ and converge to
$F$ in $L^1$ as desired.

Let $\set{C_n}$ denote any sequence tending to $\infty$. By
hypothesis, the sequence $B_{C_n}$ exhaust a full measure subset of
$X$. Since $\nu$ is a Radon measure, we may enclose the complement
of $B_{C_n}$ by an open set $O_n$ with approximately the same
measure. Since the complement of $O_n$ is a closed subset of
$B_{C_n}$ we may find compact subsets $Y_n\subset B_{C_n}$ such that
$\lim_{n\to\infty}\nu(Y_n)=1$.

To prove \eqref{eq:*}, first set $h_0(x)=0$. For the inductive
step, assume we found $h_0(x), h_1(x), \dots, h_n(x) \in V_{+}$
such that
$$0< \sum_{i=0}^{n} h_i(x) \leq \gamma_n F(x).$$

Set $R_n(x)=F(x) - \sum_{i=0}^{n}h_i(x)$. Observe that $R_n(x)$ is a
positive continuous function. Since $\gamma_n<1$ and $F$ is
uniformly positive, $R_n$ is also uniformly positive. Since $F$ and
each $h_i$ are uniformly continuous on $Y_n$, so is $R_n(x)$.  Hence
for $Y=Y_n$ and any fixed $s >1$ there exists a $\delta>0$ and
$t<\infty$, both depending on $n$, that satisfies the first
condition of Proposition \ref{prop:subfunction}.  Now by Theorem
2.8.7 in \cite{Federer69}, for any $\eps>0$, there is a Vitaly cover
of $B_{C_{n}}=\cap_{r>0} B_{C_{n}}(r)$ by a countable family of
balls
$\set{B(b_{\alpha_{i,\eps,n}},r_{\alpha_{i,\eps,n}})}_{i=1}^\infty$
with $r_{\alpha_{i,\eps,n}}<\eps$ and $\alpha_{i,\eps,n}\subset
\mathcal{A}$ for all $i\in\nn$.  In particular, we have a weak cover
of $Y_n$ by disjoint balls chosen from the family
$\set{B(b_{\alpha_{i,\delta/t,n}},r_{\alpha_{i,\delta/t,n}})}_{i=1}^\infty$
with all $r_{\alpha_{i,\delta/t,n}}<\delta/t$.  Moreover, by
definition, this weak cover has Lebesgue constant $B=1$. We may then
apply Proposition \ref{prop:subfunction} to the corresponding
countable family of $\nu$-spikes.

Therefore,for any $\beta<1$ there exists $h_{n+1}(x) \in V_{+}$ such
that
$$h_{n+1}(x) \leq \beta R_n(x)\quad \text{ for all }\quad x\in X$$
and
$$\norm{h_{n+1}(x)}_{L^1(Y_n)} \geq \frac{\beta L_{\nu}}{C_n^2\,s^2}\norm{R_n}_{L^1(Y_n)}.$$

Recall
$$R_n(x)=F(x)-\sum_{i=1}^{n} h_i(x)\geq (1-\gamma_n)F(x),$$
so it follows that
$$F(x)-\sum_{i=1}^{n+1} h_i(x)=R_{n}(x)-h_{n+1}(x)\geq
(1-\beta)R_n(x)\geq(1-\beta)(1-\gamma_n)F(x).$$
Therefore setting $\gamma_{n+1}=1-(1-\gamma_n)(1-\beta)<1$, we have
$$\sum_{i=1}^{n+1} h_i(x)\leq \gamma_{n+1}F(x).$$

Using the above estimates we have,
\numeq{
\label{eq:R_recur}
\norm{R_{n+1}}_{L^1(X)}=&\norm{R_n -
  h_{n+1}}_{L^1(X)} =
\norm{R_n}_{L^1(X)}-\norm{h_{n+1}}_{L^1(X)} \\
\leq&
\norm{R_n}_{L^1(X)}-\frac{L_{\nu}}{C_n^2s^2}\,\norm{R_n}_{L^1(Y_n)}\\
\leq & \left(1-\frac{L_{\nu}}{C_n^2\,s^2}\right)\norm{R_n}_{L^1(X)}+
\frac{L_{\nu}}{C_n^2\,s^2}\norm{F}_{L^{1}(X\setminus Y_n)} } Recall
$s$ was fixed independent of $n$.  Moreover since $F\in L^1(X,\nu)$,
$\norm{F}_{L^{1}(X\setminus Y_n)}$ tends to $0$ as $n\to \infty$.
Recall that $\frac{1}{T_\nu}$, and hence $L_\nu$, may tend to $0$ as
$C_n$ increases. However, we are free to choose how quickly $C_n$
tends to $\infty$. Therefore we choose the sequence
$\set{C_n}_{n=1}^{\infty}$ tending to $\infty$ sufficiently slowly
so that $\lim_{k\to \infty} \prod_{n=1}^k
\left(1-\frac{L_{\nu}}{C_n^2\,s^2}\right)=0$. In other words, we
force the sum $\sum_{n=1}^\infty \frac{L_{\nu}}{C_n^2}$ to
  diverge. The proof is finished by the next lemma, which we apply to
  the sequences $\delta_n=\frac{L_{\nu}}{C_n^2\,s^2}$ and
  $\eps_n=\norm{F}_{L^{1}(X\setminus Y_n)}$.

\begin{lemma}\label{lem:sequence_decay}
  Let $0\leq \delta_n \leq1$ for all $n$ and $\lim_{n \to \infty}
  \ep_n = 0$.  Let $\set{a_n}_{n=1}^{\infty}$ be a sequence of
  non-negative numbers, such that $$a_{n+1} \leq (1- \delta_n) a_n +
  \delta_n\ep_n.$$

Denote by
$\Delta^n_m=  \prod_{k=m+1}^n (1-\delta_k)$, for $0\leq m\leq n$ and
$\Delta_{-1}^n=0$ and $\ep_0=a_1$.
Then $$a_{n+1} \leq \sum_{k=0}^{n}(\Delta^n_{k}-\Delta^n_{k-1})\ep_k.$$
In particular  if $\lim_{n \to \infty}\ep_n = 0$ and
$\lim_{n\to \infty} \Delta^n_k =
  0$, for all $k \geq -1$ then $\lim_{n\to \infty} a_n = 0$
\end{lemma}

\begin{proof}
Proof is by induction.  For $n=0$ we have
$a_1 \leq (\Delta_0^0-\Delta_{-1}^0)\ep_0 = a_1$.  Assume we proved this for $k<n$. Now for $k=n$ we have
$$a_{n+1}\leq (1-\delta_n)a_n + \delta_n \ep_n \leq (1-\delta_n) \sum_{k=0}^{n-1}(\Delta^{n-1}_{k}-\Delta^{n-1}_{k-1})\ep_k + (1-(1-\delta_n)) \ep_n.$$
Since $(1-\delta_n) \Delta^{n-1}_{k} = \Delta^{n}_{k}$ and $\Delta_n^n=1$ and
$\Delta_{n-1}^n= 1-\delta_n$, we have
$$a_{n+1}\leq \sum_{k=0}^{n-1}(\Delta^{n}_{k}-\Delta^{n}_{k-1})\ep_k + (\Delta_n^n-\Delta^n_{n-1}) \ep_n =\sum_{k=0}^{n}(\Delta^{n}_{k}-\Delta^{n}_{k-1})\ep_k.$$
This proves the formula.

Observe as $0\leq \delta_k \leq 1$ we have that $\Delta_{k-1}^n \leq \Delta_{k}^n,$ and $\Delta_{k}^{n+1} \leq \Delta_{k}^n,$ for $0\leq k \leq n$.
For $\ep >0$ there exits $N$ such that for all $n >N$,
$\ep_n \leq \ep$.

So $$a_{n} \leq \sum_{k=0}^{N-1}(\Delta^{n}_{k}-\Delta^{n}_{k-1})\ep_k +
\sum_{k=N}^{n}(\Delta^{n}_{k}-\Delta^{n}_{k-1})\ep_k \leq
\sum_{k=0}^{N-1}(\Delta^{n}_{k}-\Delta^{n}_{k-1})\ep_k +
\sum_{k=N}^{n}(\Delta^{n}_{k}-\Delta^{n}_{k-1})\ep $$
Since the second sum telescopes we obtain
$$a_n \leq \sum_{k=0}^{N-1}(\Delta^{n}_{k}-\Delta^{n}_{k-1})\ep_k  + (\Delta^n_n - \Delta^n_{N-1})\ep.$$
Sending $n \to \infty$ and using the fact that $\Delta_n^n=1$ and $\lim_{n \to \infty} \Delta_k^n =0$ for every $k \geq -1$, we obtain that $0\leq \lim_{n \to \infty}a_n \leq \ep$. Sending, $\ep \to 0$, we complete the proof.
\end{proof}

Observe that all of the families of spikes to which we applied
Proposition \ref{prop:subfunction} could have been chosen from the
countable family corresponding to the indices
$\set{\alpha_{i,\frac1j,k}}_{i,j,k=1}^\infty\subset \mathcal{A}$.
This family is universal in that it does not depend on $F$, but only
on a particular choice of Vitaly covers for the $B_{C_n}$. We
henceforth re-index this family as
$\set{\alpha_i}_{i=1}^\infty\subset \mathcal{A}$. The corresponding
family of spikes is therefore the countable subfamily given by the
theorem.

The second statement of the theorem follows by replacing the $L^1$
norm by absolute value and using the pointwise estimate given in part
b) of the proposition. If $F$ is continuous the pointwise estimate
gives uniform convergence on compact sets since the limit of partial
sums is nondecreasing.

To prove the theorem for an $L^1$ almost lower semicontinuous
function $F$, we first chose a nondecreasing sequence of continuous
approximations $F_1\leq F_2\leq \dots \leq F$. We may assume
$F_j<F_{j+1}$ for all $j>0$ by replacing $F_j$ by $F_j-\eps_j$ for a
sequence of sufficiently small $\eps_j>0$ which tend to $0$. Setting
$F_0=0$ we may write
$$F=\sum_{j= 1}^\infty F_{j}-F_{j-1}.$$

Since $F_{j}-F_{j-1}$ is bounded, continuous and uniformly positive, it
satisfies the hypotheses of the theorem and we have shown that
$F_{j}-F_{j-1}=\sum_{i=1}^{\infty}\lambda_{\alpha_{i,j}}
f_{\alpha_{i,j}}.$ In particular,
$$\norm{\sum_{j=1}^{\infty}\left(\sum_{i=1}^\infty
    \lambda_{\alpha_{i,j}}f_{\alpha_{i,j}}\right)}_1=\norm{\sum_{j= 1}^\infty
  F_{j}-F_{j-1}}_1=\norm{F}_1<\infty.$$
However since the $f_{\alpha_{i,j}}$ and the
$\lambda_{\alpha_{i,j}}$ are all positive, we are free to rearrange
this sum. If we re-index the countable family
$\set{\alpha_{i,j}}_{i,j=1}^{\infty}$ into our universal family of
indices, $\set{\alpha_{i}}_{i=1}^{\infty}$, then we obtain
$$0=\norm{F-\sum_{j=1}^{\infty}\left(\sum_{i=1}^\infty
    \lambda_{\alpha_{i,j}}f_{\alpha_{i,j}}\right)}_1=\norm{F-\sum_{i=1}^{\infty}\lambda_{\alpha_i}f_{\alpha_i}}_1,$$
where
$$\lambda_{\alpha_i}=\begin{cases} 0 & \set{j
  :\alpha_i=\alpha_{i,j}}=\emptyset \\ \sum_{\set{j :\alpha_i=\alpha_{i,j}}}
\lambda_{\alpha_{i,j}} & \set{j
  :\alpha_i=\alpha_{i,j}}\neq \emptyset.
\end{cases}$$
The $\lambda_{\alpha_i}$ are the finite
nonnegative coefficients guaranteed by the theorem.

Again the second statement of the theorem follows similarly by
replacing the $L^1$ norm by the absolute value and working pointwise
in the obvious way.

\end{proof}

\section{Existence of the first moment for $Q$-spikes}\label{sec:first_moment}

In this section we would like to strengthen the result of the
Theorem \ref{thm:basis}, with some extra conditions. Recall that for
lower semicontinuous functions we had
$\sum_{i=1}^{\infty} \lambda_i \norm{f_i}_1=\norm{F}_{1}<\infty$. For
our main applications we will need slightly better convergence
properties. Namely,

\begin{theorem}[Finite Moment and Entropy Theorem]\label{thm:first_moment}
  Assume that we are in the setting of Theorem \ref{thm:basis}. For
  convenience, assume the local doubling constants $T_{\nu,C}$
  are bounded from above by $T_\nu<\infty$. In
  addition, assume that the family of unit $\nu$-spikes $\set{(f_\alpha,
  r_\alpha, b_\alpha,Q,\theta, C_\alpha)}_{\alpha\in \mathcal{A}}$ is a
  family of $Q$-spikes.

  Assume that for every $C,\epsilon>0$ there is a finite Lebesgue subcover
  cover $\set{ B(b_{\alpha_i}, r_{\alpha_i})}_{i=1}^{N(C,\eps)}$ of $B_C$
  such that $g(\epsilon) \leq r_{\alpha_i} \leq \epsilon$, for every $i
  \in \nn$ and some positive increasing function $g:[0,\infty]\to
  [0,\infty]$ with $g(r)<r$ for $1>r>0$. For any uniformly positive
  bounded Lipschitz function $F$ with Lipschitz constant $L$, by Theorem
  \ref{thm:basis} we may choose a positive sequence
  $\set{\lambda_{\alpha_i}}_{i=1}^{\infty}$ so that, $$F=
  \sum_{i=1}^{\infty} \la_i f_i(x)$$ with uniform convergence.
  Moreover, if $g(\ep)\geq a
  \ep^k$ for some $0<a<1$ and $k\geq 1$, then there is a constant $A>0$,
  independent of $F$, such that the following hold.
\begin{enumerate}

\item[1.] If $k=1$, $\nu(B_C) \leq 1-\frac{1}{C^{6}\log^2(C)}$, then we
may choose the constants $\la_i>0$ so
  that in addition, the first moment satisfies
  $$
  \sum_{i=1}^{\infty} \la_i \norm{f_i}_{L^1}
  \log{\frac{1}{\norm{f_i}_{L^1}}}<  A\left[1+\log\left(\frac{\sup_{z \in X}F(z)}{\inf_{z \in
        X}F(z)}\right)+\log L\right]\norm{F}_{L^\infty}< \infty
        $$
with entropy
  $$-\sum_{i=1}^{\infty}\la_i \norm{f_i}_{L^1} \log\left(\la_i
    \norm{f_i}_{L^1}\right)< A\left[1+\log\left(\frac{\sup_{z \in X}F(z)}{\inf_{z \in
        X}F(z)}\right)+\log L\right]\norm{F}_{L^\infty} < \infty.$$

\item[2.] If $k=1$ and $\nu(B_C) \leq 1-\frac{1}{C^2\log^3(C)}$, or
  $k>1$ and $\nu(B_C) \leq 1-\frac{1}{C^4\log^2(C)}$, then we may choose
  the constants $\la_i>0$ so that in addition, the first log-moment
  satisfies
  $$\sum_{i=1}^{\infty} \la_i \norm{f_i}_{L^1}
  \log\log{\frac{1}{\norm{f_i}_{L^1}}} < A\left[1+\log\left(1+\log\left(\frac{\sup_{z \in X}F(z)}{\inf_{z \in
        X}F(z)}\right)+\log L\right)\right]\norm{F}_{L^\infty}<
  \infty.$$

\item[3.] if $\nu(B_C)=1$ for some fixed $C\geq 1$ and
  $\frac{C^2}{C^2-L_\nu}> k \geq 1$, then we could choose the
  constants $\la_i\geq 0$, such that the first moment satisfies
  $$
  \sum_{i=1}^{\infty} \la_i \norm{f_i}_{L^1}
  \log{\frac{1}{\norm{f_i}_{L^1}}} < A\left[1+\log\left(\frac{\sup_{z \in X}F(z)}{\inf_{z \in
        X}F(z)}\right)+\log L\right]\norm{F}_{L^1}< \infty$$
  with
  entropy
  $$-\sum_{i=1}^{\infty}\la_i \norm{f_i}_{L^1} \log\left(\la_i
    \norm{f_i}_{L^1}\right)< A\left[1+\log\left(\frac{\sup_{z \in X}F(z)}{\inf_{z \in
        X}F(z)}\right)+\log L\right]\norm{F}_{L^1} < \infty.$$
\end{enumerate}

\end{theorem}

\begin{remarks}
  In our applications we will need to use basis functions normalized
  to have unit mass instead of unit height, i.e.
  $f_i/\norm{f_i}_{L^1}$ where $\norm{f_i}_{L^\infty}=1$. In this case
  $\lambda_i\norm{f_i}_{L^1}$ will be the coefficients of the sum
  instead of $\la_i$. Furthermore, in our applications to a
  group $\group$, $\log\left(\frac{1}{\norm{f_i}_{L^1}}\right)$ will
  roughly correspond to the word length $d_\group(\id,\gamma_i)$ in
  $\group$ and $\lambda_i\norm{f_i}_{L^1}$ is the value of the measure
  $\mu(\gamma_i)$. With these substitutions, the first moment, first
  log-moment and entropy formulas stated in the above theorem respectively take on their more recognizable form of
  $$\sum_{\gamma\in \group} \mu(\gamma) d_\group(\id,\gamma), \quad
   \sum_{\gamma\in \group} \mu(\gamma) \log  d_\group(\id,\gamma)\quad
   \text{and}\quad -\sum_{\gamma\in \group} \mu(\gamma) \log \mu(\gamma).$$

  One may observe that the entropy formula is nearly that of the first
  moment, except for the presence of an additional $\la_i$ within the
  $\log$ term. For an $L^\infty$ function the $\la_i$ are bounded
  above, but not below. Therefore, the entropy estimate implies the
  first moment but not conversely. In particular, finite entropy
  need not hold in the weakened estimate of case 2). In
\cite{Kaimanovich00} sufficient conditions for Poisson boundaries
are established including finite log-moment together with various
other possible estimates, not just finite entropy. It is for such
potential usefulness that we include this case here.

   The third case is of most interest to us, since this
incorporates the setting of a (convex-)cocompact group of isometries
of a CAT$(-1)$ space acting on $X=\Lambda$ (see Section
\ref{sec:derivs_are_Q}). However, this case still holds for any
measure $\alpha$-conformal density $\nu$ on a $\delta$-hyperbolic
space so long as the Busemann functions are Lipschitz. Also by
taking $k>1$ more general groups than convex-cocompact ones may be
considered.

  Suppose we wish to weaken our assumption on $\nu$ by allowing the
  local doubling constants $T_{\nu,C}$ for the spikes with
  $C_\alpha<C$ to depend on $C$ so that they can become unbounded.
  Then the theorem still holds in cases 1) and 2) if we everywhere
  replace $C$ by $\sqrt{T_{\nu,C}}C$. For instance in 1),
  $$\nu(B_C) \leq 1-\frac{1}{C^{6}\log^2 C} \quad \text{ becomes
  }\quad \nu(B_C) \leq 1-\frac{1}{T_{\nu,C}^3C^{6}\log^2
    (\sqrt{T_{\nu,C}}C)}.$$ However, in our main application of this
  theorem $T_{\nu}$ will be bounded above.
\end{remarks}

\begin{proof}[Proof of Theorem \ref{thm:first_moment}]
  We will examine the proof of Theorem \ref{thm:basis} paying
  closer attention to $\delta$ and $t$. In case 3) we will assume
  for the moment that $F$ is Lipschitz and prove the
  general case only at the end. As shown in Theorem
  \ref{thm:basis} we may restrict to a countable index subset of
  $\mathcal{A}$ whose corresponding family of spikes has the same
  properties. Therefore we will assume without loss of generality
  that $\mathcal{A}$ is countable.

  Again the construction is by induction. By our assumption
  $T_{\nu}<\infty$, it follows that the constant $L_\nu<1$ from
  Proposition \ref{prop:subfunction} is universally bounded from
  below. We will henceforth assume $L_\nu$ is this lower bound. Fix
  $s>2$ and recall that by the considerations in Theorem
  \ref{thm:basis}, we may assume $B_{C_n}$ is compact, otherwise we
  replace it by a compact set approximating it in measure. For case
  1) we must take care to choose a sequence $C_n$ tending to
  infinity sufficiently slowly to obtain convergence in $L^1$, but
  no slower than necessary as there will be a trade-off with the
  speed of convergence of $\nu(B_{C_n})$ to $1$. In case 1) choose
  the sequence $C_n$ to be $C_n =\max\set{1, \sqrt{\frac{L_{\nu}}{5
  s^2}n} }$. It will turn out that this choice is roughly optimal.

Recall that $R_N(x)=
F(x)-\sum_{n=1}^N h_n(x)$ where we build $h_n(x)$ by induction.

Assume that we have built $R_{N-1}(x)$.
There exists $\delta_N>0$ such that for
every $x \in X$ and $y\in B_{C_N}$ such that $d(x,y) \leq \delta_N$ we have
$$ s\geq \frac{R_{N-1}(y)}{R_{N-1}(x)}.$$

Set $\eps_0=1$ and recursively define $$\eps_i =
\min\left(\frac{\delta_i}{t_i}, \eps_{i-1}\right).$$
We may assume we have found $\eps_{N-1}$ and now we
show that a sufficient choice of $\delta_N$ is
$$\delta_N =\min\set{\frac{(s-1)\inf_{x\in X}R_{N-1}(x)}{\sup_{y\in
X}D_{g(\eps_{N-1})}(R_{N-1})(y)},g(\eps_{N-1})}.$$

Indeed, if $x\in X$ and $y\in B_{C_N}$ are such that $d(x,y)\leq
\delta_N\leq g(\eps_{N-1})$ then
\begin{align*}
s&\geq\frac{\delta_N\sup_{z\in
B_{C_N}}D_{g(\eps_{N-1})} R_{N-1}(z)}{\inf_{z\in X}(R_{N-1}(z))}+1\\
&\geq \frac{\delta_N}{\inf_{z\in X}(R_{N-1}(z))}\sup_{\substack{z\in
B_{C_N}\\ 0<d(x,z)\leq g(\eps_{N-1})}}\frac{\abs{R_{N-1}(z)-R_{N-1}(x)}}{d(x,z)} +1 \\
&\geq \frac{R_{N-1}(y)-R_{N-1}(x)}{R_{N-1}(x)}+1=\frac{R_{N-1}(y)}{R_{N-1}(x)}.
\end{align*}

Also recall that $$t_N=\left(\frac{\sup_{y \in
B_{C_N}}(R_{N-1}(y))}{\inf_{x\in
X}(R_{N-1}(x))}\right)^{\frac{1}{Q}}+1.$$

Since we can choose a finite Besecovitch cover
$\set{B(b_{\alpha_i^{(N)}}, r_{\alpha_i^{(N)}})}_{i=1}^{k_N}$ with
$g(\eps_N) \leq r_{\alpha_i^{(N)}}\leq \eps_N$, we can apply
Proposition \ref{prop:subfunction} b) to the function $\beta
R_{N-1}$ for a fixed constant $\beta<1$. Specifically, there exists
a $h_{N}(x) = \sum_{i=1}^{k_N} \la_i^{(N)} f_{\alpha_i^{(N)}}(x) \in
V_{+}$ such that we have
$$h_{N}(x) \leq  {\beta}R_{N-1}(x),$$ and
$$\frac{L_{\nu}\beta}{C_N^2 s^2} \norm{R_{N-1}(x)}_{B_{C_N}} \leq
\norm{h_N}_{B_{C_N}}.$$
As will become apparent, the role of the
constant $\beta$ is simply to make certain that $R_{N-1}-h_N$ is
uniformly positive. Recall that for all $N>\frac{5 s^2}{L_\nu}$, we
have $C_N=\sqrt{\frac{L_\nu}{5 s^2} N}$. From now on we fix $s\in
(1,2)$ and without loss of generality assume $\beta$ was chosen
sufficiently close to $1$ so that for all $N>\frac{5 s^2}{L_\nu}$,
$$\frac{5}{N}\geq \frac{L_{\nu}\beta}{C_N^2 s^2} \geq \frac{4}{N}.$$
(In the case 3) where $\nu(B_C)=1$ and
$\frac{1}{k}+\frac{L_{\nu}}{C^2}>1$ chose $s$ and $\beta$
sufficiently close to $1$ so that
$\frac{1}{k}+\frac{L_{\nu}\beta}{s^2 C^2}>1$.)

Therefore we obtain
$$R_{N}(x) = R_{N-1}(x) - h_{N}(x) \geq (1-\beta)R_{N-1}(x),$$
or inductively,
$$ R_N(x) \geq F(x)(1-\beta)^N.$$

Setting $Y_n=B_{C_N}$ from inequality \eqref{eq:R_recur} the fact that
$\norm{R_{N-1}}_{L^{\infty}(X)}\leq\norm{F}_{L^{\infty}(X)}$, we have
\begin{align}\label{eq:R_recurse}\norm{R_N}_{L^1(X)}
  &\leq \left(1-\frac{L_{\nu}\beta}{C_N^2
      s^2}\right)\norm{R_{N-1}}_{L^1(X)} + \frac{L_{\nu}\beta}{C_N^2
    s^2}\norm{F}_{L^{\infty}(X)}\nu(X-B_{C_N}).
\end{align}

By Lemma \ref{lem:sequence_decay}, $R_N \to 0$ since
$\ds{\sum_{i=1}^\infty \frac{1}{C_i^2} = \infty}$ and $\nu(B_{C_N})
\to 1.$ Recall that for each $a>0$, $\lim_{N\to
  \infty}N^a\prod_{i=\lbrack a\rbrack+1}^N\left(1-\frac{a}{i}\right)$
is a positive number. We may solve the recursive inequality
\eqref{eq:R_recurse}  for
$\norm{R_N}_{L^1(X)}$ given above (see Lemma \ref{lem:sequence_decay}). For some constants $C^{\pr}$ and
$C^{\pr\pr}$ we obtain
\begin{align*}
  \norm{R_{N+1}}_{L^1(X)} &\leq \norm{F}_{L^\infty(X)} \sum_{k=1}^N
  \frac{L_{\nu}\beta\,\nu(X-B_{C_k}) }{C_k^2 s^2}\left[\prod_{i=k+1}^{N}
    \left(1-\frac{L_{\nu}\beta}{C_i^2 s^2}\right)\right]\\
&\leq \norm{F}_{L^\infty(X)} \sum_{k=1}^N
  \frac{C^{\pr\pr}\,\nu(X-B_{C_k})}{k} \left[\prod_{i=k+1}^{N}
    \left(1-\frac{4}{i} \right)\right]\\
&\leq \norm{F}_{L^\infty(X)} \sum_{k=1}^N
  \frac{C^{\pr}\,\nu(X-B_{C_k})}{k} \left(\frac{k}{N}\right)^4.
\end{align*}
Therefore in case 1), the sum becomes $\ds{\norm{F}_{L^\infty(X)}
  N^{-4}\sum_{k=1}^N \frac{C^{\pr}}{\log^2(k)}.}$ It is straightforward
  to bound this sum from above we bound this sum
  from above, say by by integration. After relabelling $C^{\pr}$
  we may obtain a bound of
  $$\ds{\norm{R_{N+1}}_{L^1(X)}\leq
    \frac{C^{\pr}\norm{F}_{L^\infty(X)}}{N^{3}\log^2(N)}}.$$

In case 2), where in either case $\nu(B_C)\leq
\left(1-\frac{k}{\log(k)C^4\log^2(C)+C^2\log^3(C)}\right)$, the same
estimate yields $$\norm{R_{N+1}}_{L^1(X)} \leq
\frac{C^{\pr}\norm{F}_{L^\infty(X)}}{\log(k)N^{2}\log^2(N)+N\log^3(N)}.$$
In the case 3),
where $\nu(B_C)=1$, we have
$$\norm{R_{N+1}}_{L^1} \leq
\norm{F}_{L^1}\left(1-\frac{L_{\nu}\beta}{C^2s^2}\right)^N.$$
 Therefore, we obtain that
$$F(x) = \sum_{N=1}^{\infty}h_N(x)=\sum_{N=1}^{\infty}
\sum_{i=1}^{k_N}\la_i^{(N)} f_{\alpha_i^{(N)}}(x)$$ with convergence
in $L^1$-norm.
Also $$ \sum_{i=1}^{k_N} \la_i^{(N)}
\norm{f_{\alpha_i^{(N)}}}_{L^1}=\norm{h_N}_{L^1} \leq \beta\norm{R_{N-1}}_{L^1}.$$

Therefore, after renaming $C^{\pr}$ we have
\begin{gather}\label{eq:h_est}\tag{**}
\sum_{i=1}^{k_N} \la_i^{(N)}
\norm{f_{\alpha_i^{(N)}}}_{L^1(X)}\leq \begin{cases}
\frac{C^{\pr}\norm{F}_{L^\infty(X)}}{N^{3}\log^2(N)} & \text{case 1)} \\
\frac{C^{\pr}\norm{F}_{L^\infty(X)}}{\log(k)N^{2}\log^2(N)+N\log^3(N)} & \text{case 2)} \\
C^{\pr}\norm{F}_{L^1(X)}\left(1-\frac{L_{\nu}\beta}{C^2s^2}\right)^N & \text{case 3)}
\end{cases}
\end{gather}
Now we would like to connect $\eps_{N+1}$ to $\eps_N$. (it
might seem that it is easier to connect $\eps_N$ and $\eps_{N-1}$, but
the former relation is easier to write).

Since $h_N(x) \in V_{+}$ and $r_{\alpha_i^{(N)}}\geq g(\eps_N)$ and
$h_N(x)$ is composed of $Q$-spikes we have,
$$D_{g( \eps_N)} R_N(x) \leq D_{g( \eps_{N})} R_{N-1}(x) + D_{g( \eps_N)} h_N(x)$$

So
\begin{align*}
D_{g(\eps_N)} h_N(x) & \leq \sum_{i=1}^{k_N}
\la_i^{(N)} D_{g(\eps_N)} f_{\alpha_i^{(N)}}(x) \leq \sum_{i=1}^{k_N}
\la_i^{(N)} D_{r_{\alpha_i^{(N)}}}
f_{\alpha_i^{(N)}}(x)\\
&\leq \sum_{i=1}^{k_N} \la_i^{(N)} \frac{C_N
  f_{\alpha_i^{(N)}}(x)}{r_{\alpha_i^{(N)}}} \leq
\frac{C_N}{g(\eps_N)}h_N(x)\\
&\leq \frac{C_N}{g(\eps_N)}\, \beta R_{N-1}(x).
\end{align*}

In particular, since $g(\eps_N)\leq g(\eps_{N-1})$,
\begin{align*}
  \sup_{z\in X}(D_{g(\eps_N)} R_{N}(z)) &\leq \sup_{z\in
    X}(D_{g(\eps_N)} R_{N-1}(z))+ \sup_{z\in
    X} D_{g( \eps_{N})} h_N(z) \\
  &\leq \sup_{z\in X}(D_{g(\eps_{N-1})} R_{N-1}(z))+ \sup_{z\in
    X} D_{g(
    \eps_{N})} h_N(z), \\
  \intertext{\indent and recalling that $g(\eps_N)\leq \eps_N \leq
    \frac{\delta_N}{t_N}$ and $$\sup_{z\in X}(
    D_{g(\eps_{N-1})}R_{N-1}(z))\leq \frac{(s-1)\inf_{z\in
        X}(R_{N-1}(z))}{\delta_N}\leq \frac{(s-1)\inf_{z\in
        X}(R_{N-1}(z))}{t_N\eps_N}$$
    the previous inequality becomes, for $N\geq 1$,}
 \sup_{z\in X}(D_{g(\eps_N)} R_{N}(z)) & \leq \frac{(s-1)\inf_{z \in X}(R_{N-1}(z))}{t_N \eps_N} +
  \frac{C_N\, \beta}{g(\eps_N)}\sup_{z\in
    X}(R_{N-1}(z))\\
  & \leq \left(\frac{s-1}{t_N}+C_N\, \beta\right) \frac{\sup_{z\in
      X}(R_{N-1}(z))}{g(\eps_N)} \\
&\leq (s-1) C_N \frac{\sup_{z\in
      X}(R_{N-1}(z))}{g(\eps_N)},
\end{align*}
since $t_N\geq 2$, $\beta<1$ and $s>2$.

Therefore using the above estimate, we have
$$
\delta_{N+1} =\min\set{\frac{(s-1)\inf_{z\in
      X}(R_{N}(z))}{\sup_{z\in X}(
    D_{g(\eps_{N})}R_{N}(z))},g(\eps_{N})} \geq \frac{\inf_{z\in
    X}(R_{N}(z))}{C_N\sup_{z\in X}(R_{N-1}(z))} g(\eps_N) $$

Therefore we obtain,
$$\frac{\delta_{N+1}}{t_{N+1}} \geq \frac{\inf_{z\in X}(R_{N}(z))}{C_N\sup_{z\in X}(R_{N-1}(z))t_{N+1}} g(\eps_N).$$

Now we have
$$t_{N+1}-1 =\left(\frac{\sup_{z\in B_{C_{N+1}}}(R_N(z))}{\inf_{z
\in X} (R_N(z))}\right)^{\frac{1}{Q}}\leq
\left(\frac{1}
{(1-\beta)^N}
\frac{\sup_{z\in X}F(z)}{\inf_{z \in
X}F(z)}\right)^{\frac{1}{Q}},$$ and since $t_N\geq 2$ we also
have
$$t_N \leq 2(t_N-1).$$

Therefore we obtain
$$\frac{\delta_{N+1}}{t_{N+1}} \geq  \frac{(1-\beta)^{N(1+\frac{1}{Q})}}{2(1-\beta)C_N} g(\eps_N) \left(\frac{\inf_{z \in X}F(z)}{\sup_{z \in X}F(z)}\right)^{N(1+ \frac{1}{Q})}.$$
Since $\eps_N> (1-\beta)^{N(1+\frac{1}{Q})}
g(\eps_N)$, we have
\begin{align*}
  \eps_{N+1}=\min\left(\frac{\delta_{N+1}}{t_{N+1}}, \eps_N\right)
  &\geq \frac{(1-\beta)^{N(1+\frac{1}{Q})}}{2 C_N}
  g(\eps_N)\left(\frac{\inf_{z \in X}F(z)}{\sup_{z \in
        X}F(z)}\right)^{(1+ \frac{1}{Q})} \\
  &\geq K^{\pr}\left((1-\beta)^{N}\frac{\inf_{z \in X}F(z)}{\sup_{z \in
        X}F(z)}\right)^{2} g(\eps_N),
\end{align*}
for some constant $K^{\pr}<1$ which is independent of $N$ since $C_N$
is dominated by the additional exponential term. This is still true in
case 3) of the theorem where we may consider $C_N$ to be eventually
constant.

\begin{lemma}\label{lem:eps_estimate} Assume  $g(\ep) \geq a \ep^k$ for some $1>a>0$ and
  $k\geq 1$.
\begin{enumerate}
\item if $k=1$, then there exists $\la_0 >0$, independent of
  $F$, such that for all $N$ we have
  $$\eps_{N+1} \geq \frac{e^{-\la N^2}}{\sup_{z\in X}D_1 F(z)}\quad\text{where}\quad
  \la=\la_0\left[1+\log\left(\frac{\sup_{z \in X}F(z)}{\inf_{z \in
        X}F(z)}\right)\right].$$
\item if $k>1$, then there exists $\la_0 >0$,
  independent of $F$, such that for all $N$ we have
  $$\eps_{N+1} \geq e^{-\la k^N}\quad\text{where}\quad
  \la=\la_0\left[1+\log\left(\frac{\sup_{z \in X}F(z)}{\inf_{z \in
        X}F(z)}\right)+\log\(\sup_{z\in X}D_1 F(z)\)\right].$$
\end{enumerate}
\end{lemma}

\begin{proof}
Let $a_N= \log(\eps_N)$. Making a gross underestimate, the inequality
$$\eps_{N+1} \geq K^{\pr} {\left((1-\beta)\frac{\inf_{z \in
        X}F(z)}{\sup_{z \in X}F(z)}\right)^{2N}} g(\eps_N)$$
implies for $N\geq 1$,
$$a_{N+1} \geq 2N\log{\left((1-\beta)\frac{\inf_{z \in
        X}F(z)}{\sup_{z \in X}F(z)}\right)} + k a_N + \log{(K^{\pr}\, a)}.$$
So inductively we obtain
$$ a_{N+1} \geq 2\log{\left((1-\beta)\frac{\inf_{z \in
        X}F(z)}{\sup_{z \in X}F(z)}\right)} \sum_{i=1}^{N}(N-i) k^{i} + k^{N} a_1 +\log{(K^{\pr}\, a)} \sum_{i=1}^{N} k^i .$$

Recall that $R_0=F$ so
$\eps_1=\frac{\delta_1}{t_1}>\frac{\left(\frac{\inf_{z \in
        X}F(z)}{2\sup_{z \in X}F(z)}\right)^2 }{\sup_{z\in
    X}D_1 F(z)}$. If $k=1$, then since $1-\beta<1$ and
$K^{\pr}a<1$ we have that
$$\aligned a_{N+1} &\geq 2\log{\left((1-\beta)\frac{\inf_{z \in
        X}F(z)}{\sup_{z \in X}F(z)}\right)} \sum_{i=1}^{N}(N-i)
+a_1+N\log{(K^{\pr}\, a)} \\ & \geq 2\log{\left((1-\beta)\frac{\inf_{z \in
        X}F(z)}{\sup_{z \in X}F(z)}\right)} \frac{N^2+N}{2}
 +a_1+N\log{(K^{\pr}\, a)}\\
& \geq -\la_0\left[1+\log\left(\frac{\sup_{z \in X}F(z)}{\inf_{z \in
        X}F(z)}\right)\right] N^2 -\log\sup_{z\in X}D_1 F(z),
\endaligned$$
for some
$\la_0 >0$.  This proves 1).

If $k>1$, then
$$
\aligned a_{N+1} &\geq 2\log{\left((1-\beta)\frac{\inf_{z \in
        X}F(z)}{\sup_{z \in X}F(z)}\right)} \sum_{i=0}^{N}(N-i) k^{i} +k^Na_1+\log{(K^{\pr}\, a)} \sum_{i=0}^{N} k^i  \\
& = k^N\left( 2\log{\left((1-\beta)\frac{\inf_{z \in
        X}F(z)}{\sup_{z \in X}F(z)}\right)} \sum_{i=0}^{N}i k^{-i}
+a_1+\log{(K^{\pr}\, a)} \sum_{i=0}^{N} k^{i-N}\right) \\
& \geq -\la_0\left[1+\log\left(\frac{\sup_{z \in X}F(z)}{\inf_{z \in
        X}F(z)}\right)+\log\(\sup_{z\in X}D_1 F(z)\)\right] k^N, \endaligned$$
for some $\la_0 >0$.

This proves the lemma.
\end{proof}

\medskip
Now we continue with the proof of Theorem \ref{thm:first_moment}. By
overestimating the first case, we can incorporate both cases of Lemma
\ref{lem:eps_estimate} into one estimate. Namely, there exists a $\la>0$
given by the second part of the lemma such that
$$\eps_{N+1}\geq e^{-\la N^2 k^N}.$$

Now from the definition of a spike, we have that for
all $N$ and $i = 1, \dots, k_N$, if $K=\frac{L_{\nu}\beta}{5 s^2}$, then
$$\frac{\norm{f_{\alpha_i^{(N)}}}_{L^1}}{\norm{f_{\alpha_i^{(N)}}}_{L^{\infty}}} \geq \frac{\nu(B(b_{\alpha_i^{(N)}},r_{\alpha_i^{(N)}}))}{C_N}\geq
\frac{r_{\alpha_i^{(N)}}^Q}{C_N^2} \geq \frac{e^{-\la Q(N-1)^2
k^{(N-1)}}}{K\,N}$$

In case 3), this becomes
$$\frac{\norm{f_{\alpha_i^{(N)}}}_{L^1}}{\norm{f_{\alpha_i^{(N)}}}_{
    L^{\infty}}}\geq \frac{e^{-\la Q (N-1)^2 k^{(N-1)}}}{C^2}.$$

Now we use the estimate \eqref{eq:h_est} and the above estimates to finish
each case. We set the constant $L$ to be the global Lipschitz constant for
$F$.

Case 1). Using the fact that spikes are unit, our construction yields
$$\aligned \sum_{N=1}^{\infty}
\sum_{i=1}^{k_N} &\la_i^{(N)} \norm{f_{\alpha_i^{(N)}}}_{L^1}
\log{\left(\frac{1}{\norm{f_{\alpha_i^{(N)}}}_{L^1} }\right)}  \\
& \leq \sum_{N=1}^{\infty}\sum_{i=1}^{k_N}\la_i^{(n)}
\norm{f_{\alpha_i^{(N)}}}_{L^1} \left(\la Q  (N-1)^2 +2 \log(K N)\right)
\\& \leq
\sum_{N=1}^{\infty} \frac{C^{\pr}\norm{F}_{L^\infty(X)}}{N^{3}\log^2(N)} \left(\la Q
(N-1)^2 +2 \log(K N)\right) \\
& < A\left[1+\log\left(\frac{\sup_{z \in X}F(z)}{\inf_{z \in
        X}F(z)}\right)+\log L\right]\norm{F}_{L^\infty(X)} <\infty.
\endaligned
$$

Case 2). Using that $\log(x+y) \leq \log(x) + y$ for $x\geq 1$ and
$y>0$, for a possibly different constant $C^{\pr}$ we obtain,
$$\aligned \sum_{N=1}^{\infty}
\sum_{i=1}^{k_N}&\la_i^{(N)} \norm{f_{\alpha_i^{(N)}}}_{L^1}
\log\log{\left(\frac{1}{\norm{f_{\alpha_i^{(N)}}}_{L^1} }\right)}  \\
& \leq \sum_{N=1}^{\infty}\sum_{i=1}^{k_N}\la_i^{(n)}
\norm{f_{\alpha_i^{(N)}}}_{L^1} \log\left(\la Q  (N-1)^2 k^{N-1} +2 \log(K N)\right)
\\& \leq
\sum_{N=1}^{\infty} \frac{C^{\pr}\norm{F}_{L^\infty(X)}}{\log(k)
  N^2\log^2(N)+N\log^3(N)}\left(N\log(k)+\log(\lambda N)\right) \\
& < A\left[1+\log\left(1+\log\left(\frac{\sup_{z \in X}F(z)}{\inf_{z \in
        X}F(z)}\right)+\log L\right)\right]\norm{F}_{L^\infty(X)}<\infty.
\endaligned
$$

Case 3). The construction yields
$$\aligned \sum_{N=1}^{\infty}
\sum_{i=1}^{k_N}&\la_i^{(N)} \norm{f_{\alpha_i^{(N)}}}_{L^1}
\log{\left(\frac{1}{\norm{f_{\alpha_i^{(N)}}}_{L^1} }\right)} \\
& \leq \sum_{N=1}^{\infty}\sum_{i=1}^{k_N}\la_i^{(n)}
\norm{f_{\alpha_i^{(N)}}}_{L^1} (\la Q  (N-1)^2 k^{N-1} +2 \log(C))
\\& \leq
\sum_{N=1}^{\infty}C^{\pr}\norm{F}_{L^1}\left(1-\frac{L_{\nu}\beta}{C^2s^2}\right)^N (\la Q
(N-1)^2k^{N-1} +2 \log(C)) \\
&< A\left[1+\log\left(\frac{\sup_{z \in X}F(z)}{\inf_{z \in
        X}F(z)}\right)+\log L\right]\norm{F}_{L^1}<\infty.
\endaligned
$$

Since $\frac{1}{k} + \frac{L_{\nu}\beta}{C^2s^2}>1$ we have that
$$\left(1-\frac{L_{\nu}\beta}{C^2 s^2}\right) k <1,$$
and in this case
the sum converges. Note that in each case the constant $A$ does not
depend on $F$.

For the bound on the entropy we recall that $k_N$ is the number of
elements in the Besicovitch cover with radii on the scale of
$\epsilon_N$. By passing to a subcover, we may assume $k_N$ has been
chosen minimally. Since the cover is Besicovitch and by the condition
on $\nu$ in the definition of a $Q$-spike, we have
$$k_N\leq \frac{K}{\ep_N^{Q}}\leq K L e^{\la Q (N-1)^2 k^{N-1}}$$
for some constant $K$ independent of $N$.

Now as before we break the computation into cases. The
convexity of $-\log$ implies for case 1),
\begin{align*}
  -\sum_{N=1}^{\infty} \sum_{i=1}^{k_N} &\la_i^{(N)}
  \norm{f_{\alpha_i^{(N)}}}_{L^1} \log \left( \la_i^{(N)}
    \norm{f_{\alpha_i^{(N)}}}_{L^1} \right) \\
  & \leq - \sum_{N=1}^{\infty} \left(\sum_{i=1}^{k_N}\la_i^{(N)}
  \norm{f_{\alpha_i^{(N)}}}_{L^1}\right) \log \left(
    \frac{1}{k_N}\sum_{i=1}^{k_N}\la_i^{(N)}
    \norm{f_{\alpha_i^{(N)}}}_{L^1} \right) \\
  & \leq -\sum_{N=1}^{\infty}
  \frac{C^{\pr}\norm{F}_{L^\infty(X)}}{N^{3}\log^2(N)}
  \log\left(\frac{C^{\pr}\norm{F}_{L^\infty(X)}}{k_N N^{3}\log^2(N)} \right)\\
  & \leq \sum_{N=1}^{\infty}
  \frac{C^{\pr\pr}\norm{F}_{L^\infty(X)}}{N^{3}\log^2(N)} \left(\la Q N^2 + 4\log(N) - \log \norm{F}_{L^\infty(X)}\right)\\
  & < A\left[1+\log\left(\frac{\sup_{z \in X}F(z)}{\inf_{z \in
        X}F(z)}\right)+\log L\right]\norm{F}_{L^\infty(X)}  < \infty.
\end{align*}
Repeating this procedure for case 3) using the estimate in
\eqref{eq:h_est} is similar. This completes the proof of the
theorem.
\end{proof}

\begin{remarks}
If we wish to extend the theorem to a more general class of $F$, then we
can approximate a bounded uniformly positive lower semicontinuous $F\in
L^1(X,\nu)$ by an increasing sequence of Lipschitz (see
\cite{Miculescu01}) functions and use an argument similar to that in the
proof of Theorem \ref{thm:basis}. Because the estimates are in terms of
$\norm{F}_{L^\infty}$ in case 1) and 2), this method fails. In fact, it is
easy to construct a continuous function $F$ on the interval $[0,1]$ for
which any increasing sequence of Lipschitz approximations $F_j$ satisfies
$\sum_{j=1}^{\infty}\norm{F_j-F_{j-1}}_\infty=\infty$.

However, in case 3) this approach can work. One explicit way of doing this
is to consider the metric space version of the quadratic Moreau-Yosida
infimal convolution approximation. Namely, let $F_s$ be defined by
$$F_s(x):=\inf_{y\in X}\set{ F(y)+s d(x,y)^2}.$$
For lower semicontinuous $F$, the functions $F_s$ satisfy
\begin{enumerate}
\item For all $s\geq 0$, $F_s$ is Lipschitz with Lipschitz constant $s$.
\item For all $t\geq s\geq 0$, we have $F_s\leq F_t \leq F$ and
  $\lim_{s\to\infty} F_s(x)= F(x).$
\item For all $s\geq 0$, we have $\inf_{x\in X} F_s(x)=\inf_{x\in X}
  F(x)=H_0.$
\end{enumerate}
Items 2) and 3) are routine to verify. For the first statement, see
Chapter 5 of \cite{ClarkeLedyaevSternWolenski98} or
\cite{ClarkeLedyaevWolenski95} and observe that the proof only
depends on the triangle inequality and not on the distance being a
norm.

If $F$ arises as a strictly increasing pointwise limit of such
uniformly positive bounded Lipschitz functions $F_i$ with Lipschitz
constants $L_i$ and $$\sum_{i=1}^\infty \log(L_i)
\norm{F_i-F_{i-1}}<\infty,$$ then it is not hard to show that the
conclusion of case 3) of the theorem still holds, except for the
quantitative bounds. However this condition most likely fails for
continuous functions whose pointwise modulus of continuity is
greater than $-\frac{1}{\log}$ on a positive measure set, but we did
not check this.

\end{remarks}

\section{Radon-Nikodym derivatives are spikes for $\delta$-hyperbolic
spaces.}

\label{sec:derivs_are_spikes} Now we recall the notation of the
first section. For a Gromov $\delta$-hyperbolic space $(X,d)$ and
$\alpha>0$, we would like to find an $\alpha$-quasiconformal
measure, but we do not know that any exist. If we let
$\frac{\alpha}{\ep}$ be the Hausdorff dimension of the metric
$\delta_p^\eps$ (or equivalently, the bilipschitz quasimetric
$d_p^{\ep}$) and let $m_p$ be the corresponding
$\frac{\alpha}{\eps}$-dimensional Hausdorff measure, then $m_p$ is a
reasonable candidate for such a measure.  In the case that $X$ is
geodesic and $\Gamma$ acts quasiconvex cocompactly on $X$, Coornaert
showed (Proposition 7.5 of \cite{Coornaert93}), based on an argument
of Sullivan \cite{Sullivan79}, that $m_p$ is $\alpha$-quasiconformal
and any two $\alpha$-quasiconformal measures are equivalent with a
bounded uniformly positive Radon-Nikodym derivative. Nevertheless,
even in this setting, it is not clear that there exist continuous
$\alpha$-quasiconformal measures.

In this section we will prove Theorem \ref{thm:main1}. For this
purpose we now allow $\alpha$ to be any positive number. However, we
have already seen in Sections \ref{sec:delta_hyp} and
\ref{sec:metric_measure} that, at least for the somewhat regular
spaces, we must assume $\alpha$ to be sufficiently large in order
for $\alpha$-quasiconformal measures even to exist.

\begin{proposition}\label{prop:spike}
  Let $X$ be a $\delta$-hyperbolic metric space which is Gromov product
  bounded with respect to $\group$ and some $p\in X$. For any $\alpha>0$, let
  $\nu$ be an $\alpha$-quasiconformal measure which
  does not consist of a single atom.
  Then there exist constants $\beta\geq 1$ and $D_0>0$, independent
  of $p$, such that for all $D\ge D_0$, all $\eps>0$ and all $\ga\in
  \group$,
  $$\left(e^{2\alpha( x\cdot \ga^{-1} p)_p},e^{-\ep (U_{p,\ga}-D)},
    z_{p,\ga}^+, \frac{\alpha}{\ep},\frac{\alpha}{\ep}, \beta e^{\alpha(10\delta+2D)}\right)$$
  is a $\nu$-spike on $\D X$ with respect
  to the quasimetric $d_p^\eps$.
\end{proposition}

\medskip

\begin{proof}
We must verify the conditions of Definition \ref{def:spike}. Recall
that
$$O_p(\ga, D)= \set{ x \,: \, d_p^{\ep}( z_{p,\ga}^+, x) \leq
e^{-\ep(U_{p,\ga} -D)}}.$$ By the Shadow Lemma
(\ref{lem:Shadow_lemma}) there exist $D_0$ and $\beta\ge 1$ such
that for all $D\ge D_0$ we have
$$\nu(O_p(\ga, D)) \leq
\beta e^{-\alpha U_{p,\ga}} e^{2\alpha D}.$$
We also have
$$\frac{1}{\beta}e^{-\alpha U_{p,\ga}} \leq \nu(O_p(\ga, D)).$$
Now for all $ z,  x \in \pa X$ such that $d_p^{\ep}( z,
 x) \leq e^{\ep(U_{p,\ga} -D)}$ by Lemma \ref{lem:value_on_balls}
observe that $|( x\cdot \ga^{-1}p)_p - ( z\cdot \ga^{-1} p)_p|
\leq 2\delta +D$.  In particular
$$e^{2\alpha( x\cdot \ga^{-1}p)_p} \geq
e^{-2\alpha(2\delta+D)}e^{2\alpha( z\cdot \ga^{-1}p)_p}.$$ This implies
condition 3).

Now it is not difficult to obtain condition 1). For all $ x \in
\pa X$ such that $d_p^{\ep}( z_{p,\ga}^+  x) \leq
e^{\ep(U_{p,\ga} -D)}$, we also have that
$$e^{2\alpha( x\cdot \ga^{-1} p)_p} \geq
e^{-2\alpha(2\delta+D)}e^{2\alpha( z_{p,\ga}^+\cdot \ga^{-1}p)_p}
\ge e^{-2\alpha(4\delta+D)} e^{2\alpha U_{p,\ga}},$$ since
$( z_{p,\ga}^+\cdot \ga^{-1}p)_p \ge U_{p,\ga}-\delta$. This implies
condition 1).

\medskip

Condition 2): By definition of $U_{p,\gamma}$ it follows from the
Shadow Lemma \ref{lem:Shadow_lemma} that $\nu$ has
$(\frac{\alpha}{\eps},\frac{\alpha}{\eps})$-decay. For this case
we also assume $ x \notin O_p(\ga, D)$, so we have $(x\cdot
z^+_{p,\ga})_p \leq U_{p,\ga}-D$. On the other hand, for all $y
\in O_p(\ga, D)$ we have $(y\cdot z^+_{p,\ga} )_p \geq U_{p,\ga}-D
\geq (x\cdot z^+_{p,\ga})_p$.  Now by $\delta$-hyperbolicity we
have
$$(x\cdot y)_p \geq \min((x\cdot z^+_{p,\ga})_p, (y\cdot z^+_{p,\ga}
)_p)-\delta \geq (x\cdot z^+_{p,\ga})_p-\delta.$$

By Lemma \ref{lem:compare_to_distance} we have that for all $ x
\notin O_p(\ga, 0)$ and all $y \in O_p(\ga,D)$, we have
$$e^{\ep( x\cdot \ga^{-1}p)_p} \leq \frac{e^{3\ep
\delta}}{d_p^{\ep}( x, z_{p,\ga}^+)} \leq \frac{e^{4\ep
\delta}}{d_p^{\ep}( x, y)}.$$ Combining these estimates together,
we obtain
$$\aligned
e^{2\alpha( x\cdot \ga^{-1} p)_p} &= \frac{1}{\nu(O_p(\ga,D))}
\int_{O_p(\ga,D)} e^{2\alpha( x\cdot \ga^{-1} p)_p} d\nu(y) \\
&\leq \frac{1}{\nu(O_p(\ga,D))}\int_{O_p(\ga, D)} \frac{e^{8\alpha
\delta}}{d_p^{\ep}( x, y)^{\frac{2\alpha}{\ep}}} \nu( y)
\endaligned$$

By setting $r= e^{-\ep (U_{p,\ga}-D)}$ we obtain
$$r^{\frac{\alpha}{\eps}}e^{2\alpha( z_{p,\ga}^+\cdot \ga^{-1}
p)_p}\geq e^{\alpha (U_{p,\ga}+D)-2\alpha\delta}\geq
\frac{e^{\alpha D-2\alpha \delta}}{\beta\nu(O_p(\ga, D))}.$$ This
implies the second condition, and therefore,
$$\left(e^{2\alpha( x\cdot \ga^{-1}p)_p},e^{-\ep (U_{p,\ga}-D)},
 z_{p,\ga}, \frac{\alpha}{\ep}, \frac{\alpha}{\ep}, \max\set{\beta
e^{\alpha(10\delta -D)}, e^{\alpha(8\delta + 2D)}}\right)$$ is a
spike for every $D \ge D_0$. For simplicity we note that
$$\max\set{\beta e^{\alpha(10\delta -D)}, e^{\alpha(8\delta + 2D)}}<
\beta e^{\alpha(10\delta+2D)},$$ and so we use the right hand
expression for the spike constant.

Lastly, to show we have a $\nu$-spike, observe that the local
doubling constants are given by
$$\frac{\nu(O_p(\ga, D+\log(5)/\ep))}{\nu(O_p(\ga, D))} \leq
\beta^2 e^{2\alpha(D+\log(5)/\ep)}.$$

This proves the lemma.
\end{proof}

\begin{corollary}\label{cor:radon_spikes}

  Let $X$ be as in Proposition \ref{prop:spike} and $\nu$ a
  nonatomic $\alpha$-quasiconformal measure with constant $C$. If $\nu^{\pr}=f \nu$
  for a uniformly positive $f\in L^\infty(X,\nu)$, then there exists
  a constant $\beta\geq 1$ and $D_0>0$ such that for all $D\ge D_0$ and all $\eps>0$,
  $$\left(\frac{d(\ga_{\star} \nu^{\pr}) }{d\nu} ,e^{-\ep
      (U_{p,\ga}-D)}, z^+_{p,\ga}, \frac{\alpha}{\ep}, \frac{\alpha}{\ep},C^2 \beta e^{\alpha(10\delta+2D)}\norm{f}_{L^{\infty}}\norm{\frac1{f}}_{L^{\infty}}\right)$$
  is a $\nu$-spike on $\D X$
  with respect to the quasimetric $d_p^\eps.$

\end{corollary}

\begin{proof}
  Observe that
  $$\frac{d\ga_*\nu^{\pr}}{d\nu}(x)=f(\ga^{-1}x)\frac{d\ga_*\nu}{d\nu}(x)$$
  and that for some $C\geq 1,$
$$C^{-1}e^{2\alpha( x\cdot \ga^{-1} p)_p}\leq
e^{\alpha d(p,\ga^{-1}p)}\frac{d\ga_*\nu}{d\nu}(x)\leq C e^{2\alpha(
x\cdot \ga^{-1} p)_p}.$$
  Since $f$  is bounded from below and above we apply
  Lemma \ref{lem:spike_pinch_est} twice to the previous proposition to obtain the result.
\end{proof}

\begin{proof}[Proof of Theorem \ref{thm:main1}]
Fix $p \in X$. Corollary \ref{cor:radon_spikes} implies that the set
of 6-tuples,
$$\set{\left(\frac{d\ga_{\star}\nu}{d\nu},e^{-\ep (U_{p,\ga}-D)},
    z_{p,\ga}^+, \frac{\alpha}{\ep},\frac{\alpha}{\ep},
C(D)\right)}_{\ga \in \Ga, D\geq D_0},$$ are a family of
$\nu$-spikes where $C(D)=\beta e^{\alpha(10\delta+2D)}$. Observe
that the balls $B_{d_p^\ep}(z_{p,\ga}^+,e^{-\ep (U_{p,\ga}-D)})$ in
the quasimetric $d_p^\ep$ are exactly $O_p(\ga,D)$. Observe, that
every Radon-Nykodym derivative appears once for each $D\geq D_0$.
Noting that $U_{p,\gamma}$ is comparable to $d(p,\ga^{-1}p)$, from
the expression of Remark \ref{rem:limit_set} we can write the radial
limit set as
$$\Lambda_r(\Ga) =
\bigcup_{D\geq D_0}\bigcap_{r>0}\bigcup_{\set{\gamma\in\group
|e^{-\ep (U_{p,\ga}-D)}\leq r}} O_p(\ga, D).$$ Note that $C(D)$ and
$D$ are monotone increasing with respect to each other. Write
$D(\cdot)$ for the inverse function to $C(\cdot)$. In terms of the
notation of Theorem \ref{thm:basis}, we have
$\ds{B_{C}(r)=\bigcup_{\set{\gamma\in\group |e^{-\ep
(U_{p,\ga}-D)}\leq r}}O_p(\ga, D(C))}$.  Our assumption that the
radial limit set has full measure within the limit set implies that
$\lim_{C\to \infty} \nu(B_C)=1$. We also assumed that
$\frac{d\ga_\star\nu}{d\nu}$ is continuous for each $\ga\in\group$.
We will apply Theorem \ref{thm:basis} to this family of continuous
spikes (normalized to be unit spikes), and we take
$F=\frac{d\nu^\pr}{d\nu}$ which we assumed to be lower
semicontinuous.

Before applying the Basis Theorem, we first need to know that $\nu$
has $\left(\frac{\alpha}{\eps},\frac{\alpha}{\eps}\right)$-decay.
Since $(X,\group)$ is assumed to be quasiconvex cobounded, for any
point $z\in \pa X$, there is a $C$-quasigeodesic ray ending in $z$
for some universal constant $C$. Therefore for any $r>0$, the family
of shadows $\set{O_p(\ga,D)}_{\ga\in\group}$ has a member containing
$z$ and with $d_p^\eps$ radii bounded between $C r$ and $\frac{1}{C}
r$ for some fixed $C>1$ and any sufficiently small $r>0$. Therefore
the Shadow Lemma \ref{lem:Shadow_lemma} implies that $\nu$ is upper
and lower $\frac{\alpha}{\ep}$-regular with respect to $d_p^\eps$.
Hence Lemma \ref{lem:regular_has_decay} implies the decay property.
Incidentally, the lower regularity gives us the strong doubling
property for $\nu$, though we did not need this.

Similarly, when $\nu=\nu_p$ belongs to a bounded density
$\set{\nu_q}_{q\in X}$, case b') of the Shadow Lemma together with
the assumption on $X$ implies the upper regularity for shadows
$\set{O_p(q,D)}_{q\in X}$ which form a cover of $\Lambda_r$ with
radii bounded between $C r$ and $\frac{1}{C} r$ for some fixed $C>1$
and any sufficiently small $r>0$. This implies the upper Ahlfors
regularity for $\nu$ on $\Lambda_r$ without the quasiconvex
cobounded assumption on $\group$. In particular, $\nu$ has
$\(\frac{\alpha}{\eps},\frac{\alpha}{\eps}\)$ decay.

Theorem \ref{thm:main1} now follows from the Basis Theorem by taking
$\mu(\gamma)=\lambda_{\gamma}$, since for every measurable set $E$
we have
$$\mu \conv \nu(E)=\int_E \sum_{\gamma\in \group}
\mu(\ga)\frac{d\ga_\star\nu}{d\nu} d\nu=\int_E
\frac{d\nu^{\pr}}{d\nu}d\nu=\nu^{\pr}(E).$$
\end{proof}

\begin{proof}[Proof of Corollary \ref{cor:mu_boundary}]
By the proximality criterion established in \cite{Furstenberg73} we need
to show that for each $x\in \La$, there is a sequence $\(\ga_i\)\subset
\group$ such that ${\ga_i}_*\nu$ converges weakly to an atomic measure (of
any positive weight) at $x$. Since $x\in \La$ there is a sequence
$\(\ga_i\)$ such that for the base point $p\in X$, $\ga_ip$ tends to $x$.
In particular, $(\ga^{-1}_ip\cdot x)_p$ and $d(p,\ga_ip)$ both tend to
$\infty$ as $i$ does. Since $\nu$ is bounded $\alpha$-quasiconformal,
$\int_{\Lambda} d{\ga_i}_*\nu \geq C$ independently of $i$. On the other
hand for any open neighborhood $\mathcal{O}\subset \Lambda$ of $x$, we
have
$$\int_{\Lambda\setminus
\mathcal{O}}d{\ga_i}_*\nu=\int_{\Lambda\setminus
\mathcal{O}}\frac{d{\ga_i}_*\nu}{d\nu}(z)d\nu(z)\leq
\int_{\Lambda\setminus \mathcal{O}}Ce^{-\alpha
\rho_{p,z}(\ga_i^{-1}p)}d\nu(z)$$ which vanishes as $i$ tends to
$\infty$ since $(z\cdot\ga_i^{-1}p)_p$ remains bounded for $z\not\in
\mathcal{O}$ while $d(p,\ga_i^{-1}p)$ becomes uniformly large.
Therefore,  $\set{\gamma_*\nu}_{\ga\in\group}$ forms a proximal
family on $\Lambda$.
\end{proof}
\section{Radon-Nikodym derivatives as $Q$-spikes}\label{sec:derivs_are_Q}

In this section we prove Theorem \ref{thm:main2}. We assume that
$(H,d)$ is a CAT(-1) space and that $\nu$ is a Lipschitz
$\alpha$-quasiconformal measure on $\Lambda\subset\pa H$ for any
$\alpha>0$ for which one exists. In order to avoid treating separate
cases in the analysis below, we will also assume $\ep \leq \alpha$.

For each $\gamma\in \group$ we let
$f_\gamma(z)=\frac{d\gamma_{\star}\nu}{d\nu}(z).$ Recall that
$U_{p,\ga}=d(p,\ga^{-1}p)$ and set $z_{p,\ga}^+$ to be the end of
a geodesic starting at $p$ and passing through $\ga^{-1}p$. In
this section we will show the following proposition.
\begin{proposition} \label{prop:Q-spike} For each $p\in
H$ and $\gamma\in \group$, the tuple
$$\left(e^{2\alpha(x\cdot \ga^{-1} p)_p},e^{-\ep (d(p, \ga^{-1} p)-D)}, z^+_{p,\ga},\frac{\alpha}{\ep},\frac{\alpha}{\ep},
C(\alpha,\eps, D)\right)$$ is a $Q$-spike for
$Q=\frac{\alpha}{\ep}$, where $C(\alpha,\eps, D)$ only depends on
$\alpha$, $\eps$ and $D$.
\end{proposition}

We will need the following lemmas.
\begin{lemma} \label{lem:der_estimate}
Let $\rho,\rho^{\pr}$ be two metrics on $\pa H$. Let $x \in \pa H$.
Define $f( z)= \frac{\rho^{\pr}( x,   z)}{\rho(x,
 z)}$.

Then for any $ z,  y \not=  x$ we have

$$\frac{f( z)-f( y)}{\rho( z,  y)} \leq
\frac{f( z)}{\rho( x,  y)}+ \frac{\rho^{\pr}(z,
 y)}{\rho( z,  y)} \frac{1}{\rho( x,  y)}.$$ \end{lemma}

\medskip
\begin{proof}
$$\aligned
\frac{|f( z)-f( y)|}{\rho( z,  y)} &= \frac{|\frac{\rho^{\pr}(x,
z)}{\rho( x,  z)}-\frac{\rho^{\pr}( x,  y)|}{\rho( x,
 y)}}{\rho( z,  y)}=\frac{|\rho^{\pr}( x,   z)\rho( x,
 y)-\rho^{\pr}(x,   y)\rho( x,
 z)|}{\rho( z,  y)\rho( x,  z)\rho( x,  y)}=\\
&=\frac{\rho^{\pr}(  x,   z)|\rho( x,  y)-\rho( x,  z)|+|\rho^{\pr}(
 x,   z)- \rho^{\pr}( x,   y)|\rho( x,  z)}{\rho( z,
 y)\rho( x,  z)\rho( x,  y)}=\\
= &\frac{f( z)}{\rho( x,  y)} \frac{|\rho( x,  y)-\rho( x,
 z)|}{\rho( z,  y)} + \frac{|\rho^{\pr}( x,   z)- \rho^{\pr}( x,
  y)|}{\rho( z,  y)}\frac{1}{\rho( x,  y)}
\endaligned$$

By triangle inequality we have $|\rho( x,  y)-\rho( x,  z) |\leq
\rho( y, z)$ and $|\rho^{\pr}(  x,   z)- \rho^{\pr}( x,
 y)|\leq \rho^{\pr}( z,   y)$

Substituting we obtained required inequality.
\end{proof}

Now we apply this lemma to our Radon-Nikodym derivatives. The
following can be viewed as a quantitative version of the Lipschitz
property for Busemann functions established by Bourdon
\cite{Bourdon96}. Perhaps this can be derived directly from that
result.

\begin{lemma}\label{lem:deriv_est}
For $p,q\in H$, fix $\ep$ such that $d_p^{\ep}$ and $d_q^{\ep}$ are
metrics. Then
$$\frac{|e^{\ep(z \cdot q)_p}-e^{\ep(y \cdot q)_p}|}{d_p^{\ep}(y,z)}
\leq 4 e^{\ep d(p,q)} \frac{e^{\ep (z\cdot
q)_p}}{\op{diam}(\Lambda)},$$ for any $y,z \in \pa H$.
\end{lemma}

\begin{proof}  Fix $x \in \pa H$ such that $d_p^{\ep}(y,x)\geq \frac{\op{diam}(\Lambda)}{2}$.  Let
$$f(w)= \frac{d_q^{\ep}(w,x)}{d_p^{\ep}(w,x)}.$$

By Lemma \ref{lem:ratio_in_delta_space} we have
$$\frac{d_q^{\ep}(z,y)}{d_p^{\ep}(z,y)} \leq f(z) e^{\ep((y\cdot q)_p - (x\cdot q)_p)} \leq f(z) e^{\ep d(p,q)}$$
By Lemma \ref{lem:der_estimate} we have
$$\frac{|f(z)-f(y)|}{d_p^{\ep}(z,y)} \leq \frac{f(z)}{d_p^{\ep}(x,y)}+
\frac{d_q^{\ep}(z,y)}{d_p^{\ep}(z,y)}\frac{1}{d_p^{\ep}(x,y)} \leq
\frac{4 f(z) e^{\ep d(p,q)}}{\op{diam}(\Lambda)}.$$ Now use the fact that
$$f(w)=e^{\ep(w \cdot q)_p}e^{\ep(x \cdot q)_p},$$
to eliminate $f$ everywhere in the above formula. To finish the proof
divide everything by $e^{\ep(x \cdot q)_p}$.
\end{proof}

\begin{proof}[Proof of Proposition \ref{prop:Q-spike}]
Let $D_0$ be the constant from Lemma \ref{lem:Shadow_lemma}. Fix any
$D\geq D_0$, so that by Proposition \ref{prop:spike}
$$\left(e^{2\alpha(x\cdot \ga^{-1} p)_p},e^{-\ep (d(p, \ga^{-1}
p)-D)}, z^+_{p,\ga},\frac{\alpha}{\ep},\frac{\alpha}{\ep}, C(\alpha,
D)\right)$$ is a $\nu$-spike for $C(\alpha,D)=\beta
e^{\alpha(10\delta+2D)}$. Set $q=\ga^{-1}p$ and
$r=e^{-\ep(d(p,q)-D)}$.

Now apply the inequality $|x^k-y^k| \leq k|x-y|\max(x,y)^{k-1}$
for all $k\geq 1$ to the Lemma \ref{lem:deriv_est} to obtain
\begin{align*}
\frac{|e^{\alpha(z \cdot q)_p}-e^{\alpha(y \cdot
q)_p}|}{d_p^{\ep}(y,z)}& \leq \frac{\alpha}{\ep} \frac{|e^{\ep(z
\cdot q)_p}-e^{\ep(y \cdot q)_p}|}{d_p^{\ep}(y,z)} \max(e^{\ep(z
\cdot q)_p},e^{\ep(y \cdot q)_p})^{\frac{\alpha}{\ep}-1}\leq \\ &
\leq 4 \frac{\alpha}{\ep} e^{\ep d(p,q)} \frac{e^{\ep (z\cdot
q)_p}}{\op{diam}(\Lambda)}\max(e^{\ep(z \cdot q)_p},e^{\ep(y \cdot
q)_p})^{\frac{\alpha}{\ep}-1}.\end{align*} From the definition of
$\nu$-spike we have that $e^{\alpha(z \cdot q)_p}\leq
C(\alpha,D)e^{\alpha(y \cdot q)_p}$ for all $z,y$ such that
$d_p^{\ep}(z,y)\leq e^{-\ep(d(p,\ga^{-1}p)-D)}$. This implies that

\begin{align*} D_{r} e^{2\alpha(y \cdot q)_p} &\leq 4
\frac{\alpha}{\ep} e^{\ep d(p,q)} \frac{e^{\alpha (y\cdot
q)_p}}{\op{diam}(\Lambda)} C(\alpha,
D)^{\frac{\alpha-\ep}{\alpha}} = \\ & = \frac{4\alpha C(\alpha,
D)^{\frac{\alpha-\ep}{\alpha}}e^{\alpha D}}{\ep\op{diam}(\Lambda)}
\frac{e^{2\alpha(y \cdot q)_p}}{r}. \end{align*}

However, we assumed $\nu$ was $\alpha$-quasiconformal, so the Shadow
Lemma \ref{lem:Shadow_lemma} shows
$$\nu(B(z_{p,\ga}^+,e^{-\ep(d(p,q)-D)}))\geq \beta e^{\alpha
d(p,q)} = \beta r^{\frac{\alpha}{\ep}} e^{-\alpha D}.$$ Setting
$C(\alpha, \ep, D) = \max(\frac{4\alpha C(\alpha,
D)^{\frac{\alpha-\ep}{\alpha}}e^{\alpha
D}}{\ep\op{diam}(\Lambda)}, C(\alpha, D),\beta e^{\alpha D})$ we
finish the proof.\end{proof}

\begin{proof}[Proof of Theorem \ref{thm:main2}]  Since $\nu$ is Lipschitz
$\alpha$-quasiconformal, and in light of Lemma
\ref{lem:Lip_const_power}, each Radon-Nikodym derivative
$\frac{d\ga_\star\nu}{d\nu}$ can be expressed as
$$\frac{d\ga_\star\nu}{d\nu}(z)=R_\ga(z) e^{2\alpha(\ga^{-1}p\cdot
z)_p}e^{-\alpha d(p,\ga^{-1}p)}$$ for a Lipschitz function $R_\ga$
on $\Lambda$ satisfying $K^{-1}\leq R_\ga \leq K$ and $D_r R_\ga
\leq \frac{K}{r}$, for a constant $K\geq 1$ independent of $\ga\in
\group$ and all $r>0$. Hence by Lemma \ref{lem:func_still_Q} the
tuples
$$\left(\frac{d\ga_\star\nu}{d\nu},
e^{-\ep(d(p,\ga^{-1}p)-D)},z^+_{p,\ga},\frac{\alpha}{\ep},\frac{\alpha}{\ep},
K^2 C(\alpha,\ep, D)\right)$$ also form a family of
$\frac{\alpha}{\ep}$-spikes indexed by $\ga\in \group$ and all
$D\geq D_0$.

We assumed that $\op{diam}(X/\group)<\infty$, so we may set
$D=\max\set{\op{diam}(X/\group),D_0}$.  Hence, for every $R>0$, the
union of shadows $\bigcup_{\set{\ga | R-D\leq d(p,\ga^{-1}p)\leq
R}}O_p(\ga,D)$ covers $\La$. We assume that this was the choice of
$D$ taken in the proof of Proposition \ref{prop:Q-spike} above. In
the notation of Theorem \ref{thm:first_moment}, this means that we
can take  $g(r)= e^{-\ep D}r$ and that $B_{C}=\Lambda$ where $C=K^2
C(\alpha,\ep,D)$.

We now apply case 3) of Theorem \ref{thm:first_moment} to these
spike and set the function $F=1$. Setting $\mu(\ga)=\la_\gamma$,
then as in the proof of Theorem \ref{thm:main1} we have $\mu \conv
\nu=\nu$. Since $\mu$ also has finite first moment, we apply the
criteria of Kaimanovich to conclude that $(\Lambda,\nu)$ is a
Poisson boundary for $\mu$.
\end{proof}

\begin{remark}
We could have weakened the hypothesis on $\group$ in Theorem
\ref{thm:main2} so that $$B_C=\bigcap_{r>0}\bigcup_{\set{\ga|
e^{-\ep d(p,\ga^{-1}p)}<r}}O_p(\ga,D(C))$$ satisfies either the
assumptions of case 1) or 2) in Theorem \ref{thm:first_moment}.
However, there does not seem to be a simple intrinsic condition on
$\group$ which guarantees this.
\end{remark}

\begin{corollary}\label{cor:param_space}
In Corollary \ref{cor:mu_boundary} and Theorem \ref{thm:main2}, the
space of measures $\mu$ for which a given measure $\nu$ on $\pa X$
is stationary (respectively, a Poisson Boundary) is infinite
dimensional. Moreover, if $\eta$ is any finite Borel measure on
$\group$ then the measure $\mu$ can be chosen so that $\eta$ is
absolutely continuous with respect to $\mu$. If $\eta$ has full
support on $\group$, then $\mu$ can be chosen in the same class as
$\eta$.
\end{corollary}

\begin{proof}
if we started with any finite sum $\sum_{i=1}^k \lambda_i f_i< F$
with $f_i=\frac{d\ga_i\nu}{d\nu}$ for any choice of $\ga_i\in
\group$, then we could apply Theorem \ref{thm:first_moment} or
\ref{thm:basis} to write $F- \sum_{i=1}^k \lambda_i f_i$ in a basis
that excluded $f_1,\dots_,f_k$. In other words, for any finite set
$\set{\ga_1,\dots,\ga_k}\subset\group$, we can specify the value of
$\mu(\ga_i)$ within an interval $[0,\ep_i]$ for some $\ep_i>0$. In
particular, the space of measures for which $\nu$ is stationary, or
a Poisson boundary measure, is always infinite dimensional.

For the second statement, recall that
$f_\ga(x)=\frac{d\ga_*\nu}{d\nu}(x)$ are bounded continuous
functions on $\Lambda$. We can choose a positive function
$\eta$-measureable $w$ on $\group$ so that $w(\ga)\leq
\frac{1}{\sup_{x \in \Lambda} f_\ga(x)}$ for all $\ga \in \group$.
For each $x\in \Lambda$, $f_\ga(x)$ is continuous in $\ga$ with
respect to the topology on $\group$ inherited as a subspace of
$\op{Isom}(X)$ endowed with the compact open topology. Since $\eta$
is Borel, it follows that $w(\ga)f_\ga(x)$ is $\eta$-measurable for
all $x\in \Lambda$. By Lebesgue dominated convergence, $\int_\group
w(\ga)f_\ga d\eta(\ga)$ is continuous on all of $\Lambda$. By
dividing $w$ by a fixed constant we may assume $\int_\group
w(\ga)f_\ga d\eta(\ga)$ is less than $1$ or $\frac{d\nu^\pr}{d\nu}$,
whichever is the case.

Now we repeat the proof of the theorems by approximating a new function.
For instance, in the first theorem we approximate
$\frac{d\nu^\pr}{d\nu}-\int_\group w(\ga)f_\ga d\eta(\ga)$ in place of
$\frac{d\nu^\pr}{d\nu}$. In either case, this results in a measure
$\mu^\pr$. We finally set $\mu=\mu^\pr+ w\eta$.

For the last statement, suppose we have satisfied the hypotheses of
Corollary \ref{cor:mu_boundary}  (respectively Theorem
\ref{thm:main2}) so that we have found a subset $S\subset \group$
such that $\set{\frac{d\ga_*\nu}{d\nu}}_{\ga\in S}$ form a family of
$\nu$-spikes (resp. $Q$-spikes) satisfying the hypotheses of Theorem
\ref{thm:basis} (resp. Theorem \ref{thm:first_moment}). For each
$\ga\in S$, we wish to replace the continuous $\nu$-spikes (resp.
$Q$-spikes) corresponding to the functions $\frac{d\ga_*\nu}{d\nu}$
on $\Lambda$ with the functions $f_\ga$ given by
$$f_\ga(x)=\int_\group
\phi_\ga(\ga^\pr)\frac{d\ga^\pr_*\nu}{d\nu}(x) d\eta(\ga^\pr),$$
where $\phi_\ga$ is a yet to be chosen bounded Borel measurable
function on $\group$. We have assumed that the support of $\eta$
contains an open neighborhood of $S$. Hence for each $\ga\in S$,
there is a sequence $\(\theta_{\ga,i}\)$ of bounded continuous
functions on $\group$ such that $\theta_{\ga,i}\eta$ converges
weakly to the unit atomic measure at $\ga$. Since the
$\frac{d\ga_*\nu}{d\nu}(x)$ are continuous $\nu$-spikes (resp.
$Q$-spikes for $Q=\frac{\alpha}{\eps}$), the Lebesgue dominated
convergence theorem implies that if $\phi_\ga$ is chosen to be
$\theta_{\ga,i}$ for some sufficiently large $i=i(\ga)$, then
$f_\ga$ will still satisfy the continuous $\nu$-spike (resp.
$Q$-spike with $Q=\frac{\alpha}{\eps}$) conditions, but with each
spike constant $C_\ga$ for $\frac{d\ga_*\nu}{d\nu}(x)$ uniformly
enlarged to, for instance, $2C_\ga$ for $f_\ga$. Applying the proofs
of Theorems \ref{thm:main1} and \ref{thm:main2} to the new spikes
$f_\ga$ yields a finite measure $\mu$ of the form
$\mu=\sum_{i=1}^\infty \la(\ga_i) \phi_{\ga_i} \eta$. Consequently,
$\mu$ is absolutely continuous with respect to $\eta$, though not
necessarily conversely. However, by the penultimate argument, we may
augment $\mu$ to be the same measure class as $\eta$.
\end{proof}

\begin{remark}
In all of our constructions of stationizing measures $\mu$, the
support must at least contain a countably infinite subset of
$\group$. However, as we will see in the next section, there are
special cases of $\group$ and $(\pa \group,\nu)$ for which the
support of $\mu$ can be finite. However, these do not arise from our
construction directly.
\end{remark}

We end this section with the proof of the last corollary of the
introduction.
\begin{proof}[Proof of Corollary \ref{cor:Lip_examp}]
We first show that the Hausdorff measures and the Patterson-Sullivan
measures, assuming they exist for the given $\group$ action on $H$, are
Lipschitz $\alpha$-quasiconformal measures. By Proposition
\ref{prop:haus_conformal} and Corollary \ref{cor:patt-sull_conformal}, we
know that these measures are $\alpha$-conformal measures. However,
Proposition \ref{lem:deriv_est} shows that the Gromov product, and
consequently Busemann functions, are Lipschitz on the boundary of CAT(-1)
spaces. Hence these measures have Lipschitz Radon-Nikodym derivatives with
respect to the $\group$ transformations.

If $m$ is either such a Hausdorff or Patterson-Sullivan measure and $\nu=f
m$ for a Lipschitz $f$, then we have
$$\frac{d\ga_*\nu}{d\nu}(x)=\frac{d\ga_*f m}{df
m}(x)=\frac{f(\ga^{-1}x)}{f(x)}\frac{d\ga_*m}{dm}.$$

We assumed $\frac{1}{K}\leq f<K $ and $D_r f<K$ for some $K\geq 1$
and all $r$ less than some fixed value $C$. On the other hand, we
have by Lemma \ref{lem:Lip_const},
$$D_r\(f\of \ga^{-1}\)(x)\leq
D_{r D_r\(\ga^{-1}\)(x)}\(f\)\(\ga^{-1}x\) D_r\(\ga^{-1}\)(x).$$ However,
we only need to estimate this for $r=e^{-\eps d(p,\ga^{-1}p)}$. Since
$D_r\(\ga^{-1}\)(x)<C e^{\eps d(p,\ga^{-1}p)}$ for some $C>0$ we obtain
$D_r\(f\of \ga^{-1}\)(x)<\frac{K C}{r}$. Hence by Lemma
\ref{lem:Lip_const} $D_r
\frac{f(\ga^{-1})}{f}<\frac{K^3 C}{r}$. Lemma \ref{lem:func_still_Q}
implies that $\frac{d\ga_*\nu}{d\nu}$ are $Q$-spikes with spike constant
$2 K^5 C$ which is independent of $\ga$. This allows us to employ the rest
of the proof of Theorem \ref{thm:main2}.
\end{proof}
\section{Application to the Free Groups}\label{sec:free_group}

Here we present an application of Theorem \ref{thm:first_moment} to
the case of a free group $\group$. It is easy to verify that the
Cayley graph of $\group$, equipped with any metric quasi-isometric
to the word metric, is a Gromov hyperbolic space. With even less
effort one may check that its Cayley graph for any set of generators
is a CAT($-\kappa$) space for any $\kappa\leq 0$.

We will prove that the (class of) Patterson-Sullivan measures is harmonic,
i.e., there exists a random walk on the free
 group such that the induced measure on the boundary is in the maximal Hausdorff class.

 \medskip
 We show that normalized Radon-Nikodym derivative satisfy
 Theorem \ref{thm:first_moment}

 Let $X$ be a Cayley graph for a free group $\group$ with standard set of
 generators (i.e. generators without any relations) and any positive weights on
 the graph. Let $d$ be the corresponding weighted path metric.

 Let $x,y \in X$ and $ z \in \pa X$.  Let $z_o \in X$ be a
 point on the segment $[x,y]$ and geodesic connecting $x$ and
 $ z$ such that the $d(x,z_o)$ achieves the maximum.  Then the
 Busemann function can be represented as,
 $$\rho_{x,  z}(y) = d(x,z_o)-d(y,z_o).$$

Therefore for $ z,  o \in \pa X$. Let $z \in X$ be a point on
the geodesic $[ z, o]$ connecting $ z$ and $ o$ which is
closest to $x$.  Then
$$d_x^{\ep}( z,  o) = e^{-\ep d(x,y)}.$$

Let $\alpha/\ep$ be the Hausdorff dimension for the metric
$d_{p}^{\ep}$ and $\nu$ be the corresponding $\alpha$-conformal
density, which coincides with the Patterson-Sullivan measure at $p$.
Then we have for $p=\id$ and all $\ga \in \group$
$$\frac{d \ga_{\star} \nu}{d\nu}( z) = e^{-\alpha \rho_{\ga^{-1}p,  z}(p)}.$$
So since $-d(p, \ga^{-1}( p)) \leq \rho_{\ga^{-1}p,  z}(p)\leq d(p,
\ga^{-1}( p))=\|\ga\|$, we have that
$$e^{-\alpha\|\ga\|}\leq \frac{d\ga_{\star} \nu}{d\nu}( z)\leq
e^{\alpha\|\ga\|}.$$

Let $f_\ga( z) =\frac{d\ga_{\star} \nu}{d\nu}( z)$.  It is not difficult
to check that $(f_\ga( z),  z_\ga, e^{-\ep\|\ga\|}, \alpha,\alpha,
e^{-1})$ is a spike, where $ z_\ga$ is the endpoint of any geodesic which
starts at $p$ and passes through $\ga^{-1} p$.

Now let us estimate $D_{e^{-\ep\|\ga\|}}f_\ga( z)$.

\begin{lemma}\label{lem:loc_const} For any $q,p\in X$, the Busemann function
  $\rho_{q,\cdot}(p)$ is locally constant on $\partial X$ taking
  on  only a finite number of values. Furthermore,
  $D_{e^{-\ep\|\ga\|}}f_\ga( z)=0$ for all $\ga\in \group$.
\end{lemma}

\begin{proof}
Let $ o, z \in \pa X$ be such that $d_{p}^\eps( o, z) < e^{-\ep
d(p,q)}$. In particular, the point $y$ lying on the geodesic
  $[o,z]$ and closest to $p$ does not lie on the segment $[p,q]$.

Let $x$ be the point of the segment $[p, q]$ closest to $y$. Now
since $y\not=x$, we have that
$$\rho_{q,  z}(p) = d(y, q)- d(y, p)= d(x,q)- d(x,p)$$ on the other hand,
$$\rho_{q,  o}(p) = d(y,q)- d(y, p)= d(x,q)- d(x,p).$$
So $\rho_{q,  z}(p)=\rho_{q,  o}(p)$. The first claim follows from
the fact that a finite number of $d_p$ balls with radius $e^{-\ep
d(p,q)}$ cover $\D X$.

For the second statement, let $ o, z \in \pa X$ be such that
$d_{p}^\eps( o, z) < e^{-\ep \norm{\ga}}$. Since
$\norm{\ga}=d(\ga^{-1}p,p)$, we have $\rho_{\ga^{-1}p,
z}(p)=\rho_{\ga^{-1}p, o}(p)$ or $f_\ga( z) = f_\ga( o)$.
\end{proof}

It is clear that the ball of the radius $e^{-\ep n}$ form a Vitally
Covering for any $n \in \nn$. Therefore, we can apply Theorem
\ref{thm:first_moment} to obtain that $$1= \sum_{\ga\in \group}
\la_\ga f_\ga( z)$$ or equivalently
$$\nu = \sum_{\ga\in \group} (\la_\ga \|f_\ga\|_{L^1}) \ga_{\star}\nu.$$
with the property that
$$-\sum_{\ga\in \group} (\la_\ga \|f_\ga\|_{L^1}) \log(\|f_\ga\|_{L^1}) = \sum_{\ga\in \group} \la_\ga
\|f_\ga\|_{L^1}(\ep \|\ga\|) <\infty.$$

Hence if $\mu(\ga)=\la_\ga \|f_\ga\|_{L^1}$, then $\mu$ is a probability
measure with finite first moment and $\nu$ is $\mu$-stationary.

In the case when the edge weights only depend on the reduced word
distance to the identity, the measure $\mu$ can be represented very
simply. If $\mu$ is constant on a single sphere about the identity,
$\norm{\ga}=k$, then it is easy to compute using Lemma
\ref{lem:loc_const} that $\mu\conv\nu=\nu$. More generally, any
measure which is constant on spheres has this property since the
average of two measures which stationize $\nu$ also stationizes
$\nu$. Although from Corollary \ref{cor:param_space} we see that
these are only the most symmetric of the measures which stationize
$\nu$.

\providecommand{\bysame}{\leavevmode\hbox to3em{\hrulefill}\thinspace}
\providecommand{\MR}{\relax\ifhmode\unskip\space\fi MR }
\providecommand{\MRhref}[2]{%
  \href{http://www.ams.org/mathscinet-getitem?mr=#1}{#2}
}
\providecommand{\href}[2]{#2}

\end{document}